\documentclass[10pt,reqno]{amsart}
\usepackage{graphicx}
\usepackage{hyperref}
\usepackage{amssymb}
\usepackage{amsthm}
\usepackage{tikz}
\usetikzlibrary{intersections}
\usetikzlibrary{positioning}
\usetikzlibrary{graphs}

\usepackage{verbatim}
\usepackage{empheq} %A lateral brace to group a set of equations with labels
\usepackage{bbm}%kongxin numbers
\usepackage{array}
\usepackage{hyperref}
\hypersetup{
    colorlinks=true,
    linkcolor=blue,
    filecolor=magenta,      
    urlcolor=cyan,
    citecolor=blue,
    pdftitle={Sharelatex Example},
}
\newcommand{\ol}{\overline}

\newcommand{\RN}[1]{%
  \textup{\uppercase\expandafter{\romannumeral#1}}%
}
\newtheorem{thm}{Theorem}[section]
\newtheorem{thmx}{Theorem}[section]

\newtheorem{defn}{Definition}[section]
\newtheorem{lem}{Lemma}[section]
\newtheorem{cor}{Corollary}[section]

\newtheorem{rem}{Remark}[section]

\numberwithin{equation}{section}

\def\E{\mathbb{E}}

\def\N{\mathbb{N}}

\def\R{\mathbb{R}}
\def\S{\mathbb{S}}
\def\T{\mathbb{T}}
\def\Z{\mathbb{Z}}

\def\BB{\mathcal{B}}

\def\HH{\mathcal{H}}

\def\LL{\mathcal{L}}
\def\SS{\mathcal{S}}

\def\MM{\mathcal{M}}

\def\andd{\,\,\text{and}\,\,}

\def\div{{\rm div}}

\def\pa{\partial}

\def\sm{\setminus}

\def\al{\alpha}
\def\be{\beta}
\def\ep{\epsilon}

\def\la{\lambda}

\def\de{\delta}

\def\ga{\gamma}

\def\vp{\varphi}

\begin{document}
%---------------------------------------------------------------%
\title[Stochastic stability of physical measures in conservative systems]{Stochastic stability of physical measures in conservative systems}

\author{Weiwei Qi}
\address{State Key Laboratory of Mathematical Sciences (SKLMS), Academy of Mathematics and Systems Science, Chinese Academy of Sciences, Beijing 100190, China}
\email{wwqi@amss.ac.cn}

\thanks{ W.Q. was partially supported by a start-up grant from Academy of Mathematics and Systems Science. 
Z.S. was partially supported by NSERC RGPIN-2018-04371 and NSERC RGPIN-2024-04938.  Y.Y. was partially supported by NSERC RGPIN-2020-04451, a faculty development grant from the University of Alberta, and a Scholarship from Jilin University.}

\author{Zhongwei Shen}
\address{Department of Mathematical and Statistical Sciences, University of Alberta, Edmonton, AB T6G 2G1, Canada}
\email{zhongwei@ualberta.ca}

\author{Yingfei Yi}
\address{Department of Mathematical and Statistical Sciences, University of Alberta, Edmonton, AB T6G 2G1, Canada, and School of Mathematics, Jilin University, Changchun 130012, PRC}
\email{yingfei@ualberta.ca}

\begin{abstract}

Given the significance of physical measures in understanding the complexity of dynamical systems as well as the noisy nature of real-world systems, investigating the stability of physical measures under noise perturbations is undoubtedly a fundamental issue in both theory and practice.

The present paper is devoted to the stochastic stability of physical measures for conservative systems on a smooth, connected, and closed Riemannian manifold. It is assumed that a conservative system admits an invariant measure with a positive and mildly regular density. Our findings affirm, in particular, that such an invariant measure has strong stochastic stability whenever it is physical, that is, for a large class of small random perturbations, the density of this invariant measure is the zero-noise limit in $L^{1}$ of the densities of unique stationary measures of corresponding randomly perturbed systems. Stochastic stability in a stronger sense is obtained under additional assumptions. Examples are constructed to demonstrate that stochastic stability could occur even if the invariant measure is non-physical. The high non-triviality of constructing such examples asserts the sharpness of the stochastic stability conclusion. Similar results are established for conservative systems on bounded domains.

Our approach to establishing stochastic stability is rooted in the analysis of Fokker-Planck equations associated with randomly perturbed systems. The crucial element in our proof is the establishment of uniform-in-noise estimates in Sobolev spaces and positive lower and upper bounds for the densities of stationary measures, which are natural yet far-reaching consequences of the conservativeness of the unperturbed system. Not only do these results have stochastic stability as immediate results, but also they allow us to achieve in particular an optimal lower bound on the exponential convergence rate to the stationary measure and readily confirm the so-called sub-exponential large deviation principle of stationary measures. A distinguishing feature of our approach is that it does not rely on uniform, non-uniform, or partial hyperbolicity assumptions, which are often required in the existing literature when investigating stochastic stability. Consequently, our study opens up a new avenue for the exploration of stochastic stability and related issues.

\end{abstract}

\subjclass[2010]{Primary 37A50,  34F05, 37H30; secondary 60J60, 60H10}

%\date{January 1, 2001 and, in revised form, June 22, 2001.}


\keywords{stochastic stability, conservative system, invariant measure, physical measure, ergodicity, stationary measure}

\maketitle

\tableofcontents
%---------------------------------------------------------------%

%\section{\bf Lotka-Volterra Model}
%\url{http://www.tiem.utk.edu/~gross/bioed/bealsmodules/competition.html} 

%-----------------------------------------------------------------%
\section{\bf Introduction}\label{sec-intro}

Invariant measures are fundamental objects in the study of statistical properties of dynamical systems and are especially powerful in the characterization of complex dynamics. Among these, physical measures \cite{MR0800052,MR1933431} -- those that can be observed -- are of notable interest and paramount relevance in numerous practical applications across various scientific and engineering disciplines. Statistical properties such as entropy, Lyapunov exponents, and mixing properties, quantified by physical measures, often provide valuable information about a system's predictability, sensitivity to initial conditions, and overall complexity. Since the pioneering work on Axiom A attractors \cite{MR0399421,MR0511655,MR0415683,MR0277003}, finding physical measures that can capture the complexity of systems has attracted a lot of attention resulting in a substantial body of literature. Interested readers are referred to surveys \cite{MR1933431,MR3078694} and references therein for earlier developments, and to \cite{MR3103173,MR3210148,MR3385635,MR3607577,MR3712997,MR3742478,MR3868016,MR3881126,MR4396671,MR4405573}, to name just a few, for more recent ones. 

Given that real-world systems are intrinsically slightly noisy, investigating the stability of dynamical systems, particularly their key dynamical properties, under small noise perturbations is of both theoretical and practical significance \cite{MR1652127,MR1015933,MR1037009}. This makes it imperative and fundamental to study the stochastic stability of physical measures, especially those related to complex dynamical behaviors, that is crucial for making reliable predictions about the long-term behaviors and helps in understanding how noises affect systems' dynamics. In the present paper, we focus on this issue for conservative systems. Besides addressing the stochastic stability issue, we achieve an optimal lower bound on the exponential convergence rate to the stationary measure as well as the sub-exponential large deviation principle of stationary measures in the setting of randomly perturbed conservative systems.

%%%%%%%%%%%%%

\subsection{Setup}

We formulate the problem and present our findings in dimension $d\geq2$; the one-dimensional case is special and treated separately in Section \ref{sec-1D}. From now on, we are going to use some special notations and direct the reader to Subsection \ref{subsec-notaion} for their precise meanings. 

Let $(M,g)$ be a $d$-dimensional smooth, connected and closed manifold endowed with a Riemannian metric $g$, and denote by ${\rm Vol}$ the Riemannian density of $M$. It is assumed that $g$ is normalized so that ${\rm Vol}(M)=1$ (see \cite[Chapter 16]{MR2954043}). Consider the following ordinary differential equation (ODE) over $M$:
\begin{equation}\label{eqn-ode}
    \dot{x}=B(x),
\end{equation}
where $B: M\to TM$ is a Lipschitz continuous vector field on $M$, ensuring the global well-posedness of \eqref{eqn-ode}. The system \eqref{eqn-ode} is assumed to be conservative (or generalized volume-preserving) in the sense that \eqref{eqn-ode} admits an invariant measure $\mu_0$ with a positive density $u_0\in W^{1,p_0}$ for some $p_0>d$ with respect to  ${\rm Vol}$. Clearly, $\div (u_0 B)=0$ ${\rm Vol}$-a.e. (see Lemma \ref{lem-new-measure-vectors}). 

Note that \eqref{eqn-ode} may admit many and even infinitely many invariant measures like $\mu_0$; typical examples are rotations on tori with rationally dependent frequencies. In which case, $\mu_0$ is a representation of these invariant measures.

%In general, the flow $\vp^t$ admits multiple invariant measures except for $\mu_0$. A fundamental issue in the theory of dynamical systems is to figure out which measure survives from noise perturbation. Since both deterministic and random noises are ubiquitous in the physical world, the purpose of the present paper is to investigate the stability of invariant measures of $\vp^t$ under general noise perturbations.

To investigate the stochastic stability of $\mu_0$ or any similar invariant measures, we examine the following small random perturbation of \eqref{eqn-ode}:
\begin{equation}\label{sde-manifold}
dX_t^{\ep}=B(X_t^{\ep})dt+\ep^2 A^{\ep}_0(X_t^{\ep})dt +\ep\sum_{i=1}^{m} A^{\ep}_{i}(X_t^{\ep})\circ dW^i_t,
\end{equation}
where $0<\ep\ll1$ is the noise intensity, $m\geq d$, $A=\left\{A^{\ep}_i, i\in \{0, 1,\dots m\},\,\,\ep\right\}$ is a collection of vector fields on $M$, $\{W^i_t\}$ are $m$ independent and standard one-dimensional Brownian motions on some probability space, and the stochastic integrals are understood in the sense of Stratonovich. The vector fields $A$ are chosen from the \emph{admissible class} defined as follows. 

% A stochastic process $X^{\ep}_t$ is a solution of \eqref{sde-manifold} if for any $f\in C^2$,
% \begin{equation*}
%     f(X^{\ep}_t)=f(X^{\ep}_s)+\int_0^t (Bf)(X^{\ep}_{s'})ds'+\ep^2\int_0^t (A^{\ep}_0f)(X^{\ep}_{s'})ds'+\ep\int_0^t (A^{\ep}_if)(X^{\ep}_{s'})\circ dW^{i}_{s'},\quad \forall s<t.
% \end{equation*}

\begin{defn}[Admissible class]\label{admissible-class}
    A collection $A$ of vector fields on $M$ is said to be in the {\em admissible class $\mathcal{A}$} if $A=\left\{A^{\ep}_{i},i\in\{0,\dots,m\},\,\,\ep\right\}$ for some $m\geq d$ and the following conditions are satisfied:
    \begin{enumerate}
\item [{\bf(A1)}] there exists $p>d$ such that $A^{\ep}_0\in L^p$, $A^{\ep}_i\in W^{1,p}$ for $i\in\{1,\dots,m\}$, and 
$$
\|A^{\ep}_0\|_p+ \max_i\|A^{\ep}_i\|_{1,p}\lesssim1;
$$

\item [{\bf(A2)}] there exists $\la>0$ such that 
\begin{equation*}
\inf_{\ep}\sum_{i=1}^m|A^{\ep}_i f|^2\geq \la |\nabla f|^2\quad {\rm Vol}\text{-a.e.}, \quad \forall f\in W^{1,1}.
\end{equation*}
\end{enumerate}
\end{defn}

The condition {\bf(A1)} is a mild integrability condition. The condition {\bf(A2)} is a uniform-in-$\ep$ positivity condition, guaranteeing in particular the non-degeneracy of the stochastic differential equation (SDE) \eqref{sde-manifold} when $A\in\mathcal{A}$.

It should be pointed out that for $A\in\mathcal{A}$, the SDE \eqref{sde-manifold} may not be well-posed even in the weak sense, and therefore, transition probabilities and stationary distributions are hardly defined. However, it is quite convenient to work with the Fokker-Planck equation associated with the SDE \eqref{sde-manifold}:
\begin{equation}\label{eqn-fpe-time}
    \pa_t u=\LL^*_{\ep} u,
\end{equation}
where $\LL^*_{\ep}$ is the Fokker-Planck operator, which is the formal  $L^2$-adjoint operator of the generator $\LL_{\ep}$ given by
$$
\LL_{\ep} :=\frac{\ep^2}{2}\sum_{i=1}^m (A^{\ep}_i)^2 + \ep^2 A^{\ep}_0 +B.
$$
This weak formalism of the SDE \eqref{sde-manifold} is often adopted when its coefficients have low regularity (see e.g. \cite{MR3443169}). We are mostly interested in \emph{stationary (probability) measures} of \eqref{sde-manifold} (that is, stationary solutions of \eqref{eqn-fpe-time} in the class of probability measures) that generalize stationary distributions of SDEs. The reader is referred to Appendix \ref{app-stationary-measure} for the basics of stationary measures.

According to Theorem \ref{lem-stationary-measure}, if $A\in\mathcal{A}$, then for each $\ep$, the SDE \eqref{sde-manifold} admits a unique stationary measure $\mu_{\ep}$ having a positive density $u_{\ep}\in W^{1,p}$, where $p>d$ is the same as in {\bf (A1)}. We have suppressed the dependence of $\mu_{\ep}$ and $u_{\ep}$ on $A$; this shall cause no trouble. Moreover, the set $\MM_{A}$ of all the limiting measures of $\mu_{\ep}$ as $\ep\to0$ under the weak*-topology is non-empty thanks to the compactness of $M$ and Prokhorov's theorem, and each element in $\MM_{A}$ must be an invariant measure of \eqref{eqn-ode} \cite{MR1015933}. 

The stochastic stability of an invariant measure of \eqref{eqn-ode} is defined as follows.

\begin{defn}[Stochastic stability]
    An invariant measure $\mu$ of \eqref{eqn-ode} is said to be \emph{stochastically stable} with respect to $\mathcal{A}$ if $\MM_{A}=\{\mu\}$ for all $A\in\mathcal{A}$.
\end{defn}

In the present paper, we mainly focus on addressing the stochastic stability issue of the invariant measure $\mu_0$. The \emph{strong admissible class} of vector fields is considered in pursuit of stronger results regarding the stochastic stability.

\begin{defn}[Strong admissible class]
    A collection $A$ of vector fields on $M$ is said to be in the {\em strong admissible class $\mathcal{SA}$} if $A=\left\{A^{\ep}_{i},i\in\{0,\dots,m\},\,\,\ep\right\}$ for some $m\geq d$ and it satisfies {\bf(A2)} and
\begin{enumerate}
\item [{\bf(SA1)}] there exists $p>d$ such that 
$$
\|A^{\ep}_0\|_{1,p}+\max_i\|A^{\ep}_i\|_{2,p}\lesssim1.
$$
\end{enumerate}
\end{defn}

%In order to investigate the stochastic stability of invariant measures of $\vp^t$, we study the limiting behaviors of $\mu_{\ep}$ as $\ep\to 0$. Clearly, the limiting measure(s) of $\mu_{\ep}$ is (are) invariant measure(s) of $\vp^t$. It is interesting to know what kind of invariant measures could survive from noise perturbations. Are they physical measures? Can some physical measure lose stability under noise perturbations? Particularly, how about the persistence of $\mu_0$ under noise perturbations?

%%%%%%%%%%%%%%%%%%%%%%%%%

\subsection{Statement of main results}\label{subsec-statement-main-results}

Our first result addresses uniform estimates of stationary measures $\{u_{\ep}\}_{\ep}$ as well as the stochastic stability of $\mu_0$.

Recall that an invariant measure $\mu$ of \eqref{eqn-ode} is called a \emph{physical measure} (or \emph{physical}) if ${\rm Vol}(B_{\mu})>0$, where 
$$
B_{\mu}:=\left\{x\in M:\lim_{t\to\infty}\frac{1}{t}\int_{0}^{t}\de_{\vp^{s}(x)}ds=\mu\quad\text{under the weak*-topology}\right\}
$$
is referred to as the basin of $\mu$, where $\vp^t$ denotes the flow generated by solutions of \eqref{eqn-ode}.

\begin{thmx}\label{thm-uniform-estimate-stability}
The following statements hold.
\begin{itemize}
    \item [(1)] For any $A\in \mathcal{A}$, there hold
$$
\|u_{\ep}\|_{1,2}\lesssim1\quad\text{and}\quad 1\lesssim\min u_{\ep}\leq \max u_{\ep}\lesssim1.
$$
In particular, any $\mu\in\MM_{A}$ has a density $u$ belonging to $W^{1,2}$ and satisfying $u,\frac{1}{u}\in L^{\infty}$.

\item [(2)] If $\mu_0$ is physical, then it is stochastically stable with respect to $\mathcal{A}$. 

\item [(3)] If $\mu_0$ is the only invariant measure of \eqref{eqn-ode} with a density in $W^{1,2}$, then it is stochastically stable with respect to $\mathcal{A}$. 
\end{itemize}
In the case of either (2) or (3), the limit $\lim_{\ep\to0}u_{\ep}=u_{0}$ holds weakly in $W^{1,2}$ and strongly in $L^{p}$ for any $p\in[1,\infty)$.
\end{thmx}

\begin{rem}
We make some comments about Theorem \ref{thm-uniform-estimate-stability}.
\begin{itemize}
    \item[(i)] The uniform estimates of $\{u_{\ep}\}_{\ep}$, having conclusions in (2) and (3) as immediate consequences, are natural yet far-reaching consequences of the conservativeness of the unperturbed system \eqref{eqn-ode}. Our approach to establishing these uniform estimates builds on analyzing the stationary Fokker-Planck equation satisfied by $u_{\ep}$, namely, $\LL_{\ep}^{*}u_{\ep}=0$ in the weak sense. Developing uniform-in-$\ep$ Harnack's estimate and Moser iteration for $u_{\ep}$ plays a crucial role in the proof. 

    In addition to the implications for stochastic stability, these uniform estimates are key to the justification of the sub-exponential large deviation principle of $\{u_{\ep}\}_{\ep}$ (see Remark \ref{rem-ldp}) and the derivation of more than just an optimal lower bound on the exponential convergence rate to the stationary measure (see Theorem \ref{thm-chi-convergence} and the subsequent discussion).

    \item[(ii)] It is shown in Corollary \ref{cor-2024-03-06} that $\mu_0$ is physical if and only if it is ergodic. The ergodicity of $\mu_0$ clearly implies that it is the only invariant measure of \eqref{eqn-ode} with a density in $W^{1,2}$, yielding that (2) is stronger than (3). In Theorem \ref{thm-construction}, we construct examples that satisfy (3) while not fulfilling (2). It turns out that constructing such examples is highly non-trivial, asserting the sharpness of (2).

    \item[(iii)] The limit $\lim_{\ep\to0}u_{\ep}=u_{0}$ holds in $L^{1}$ for each $A\in\mathcal{A}$. In literature, such a result is often called the \emph{strong} stochastic stability with respect to $\mathcal{A}$.

   \item[(iv)] In the case where $B$ is divergence free, the stochastic stability of the volume {\rm Vol} with respect to homogeneous noises (that is, $A^{\ep}_0=0$ and $A^{\ep}_{i}=A_i$ is a constant vector field for each $i\in\{1,\dots,m\}$) was obtained in \cite{zbMATH05852930}. This result is straightforward as the unique stationary measure $\mu_{\ep}$ coincides with {\rm Vol}. It is generalized in Theorem \ref{thm-converse} that sheds light on the issues of \emph{invariant measure selection by noise} and \emph{stochastic instability} when \eqref{eqn-ode} admits multiple invariant measures like $\mu_0$.

%It is worth noting that these uniform estimates are in sharp contrast to limiting behaviours of stationary measures of randomly perturbed {\em dissipative} systems. In the latter scenario,  stationary measures tend to concentrate on attractors of the unperturbed system, resulting in the singularity of their densities as noise vanishes (see e.g. \cite{MR1652127,MR3771195,MR4575665}), and hence, the failure of uniform estimates of densities as in Theorem \ref{thm-uniform-estimate-stability}.
\end{itemize}
\end{rem}

Introduced by Kolmogorov and Sinai \cite{MR0399421,MR1009437}, the stochastic stability of invariant measures, especially physical measures, has attracted considerable attention and received affirmative results in many instances. The reader is referred to \cite{MR1015933} and references therein for earlier investigations concerning uniquely ergodic systems or systems having relatively simple dynamics. Significant advancements have been achieved in various settings, including uniformly hyperbolic systems \cite{MR0388452,MR0857204}, systems that feature spectral gaps applicable to many exponentially mixing piecewise expanding maps \cite{MR0685377,MR1233850,MR1386223,MR1430741,MR1644405,MR1679080}, and non-uniformly hyperbolic maps \cite{MR0874051,MR1162396,MR1772247,MR2052296,MR2069701,MR2134078,MR2259614,MR2351040,MR3062897,MR3126392,MR3170624}.

A distinguishing feature of our results (Theorem \ref{thm-uniform-estimate-stability} and the subsequent Theorems \ref{thm-construction} and \ref{thm-stronger-uniform-bounds}) is that they do not rely on uniform, non-uniform, or partial hyperbolicity assumptions, which are typically required in most previous works. While these hyperbolicity assumptions have become standard in the geometric or statistical theory of smooth dynamical systems, they are often not met or challenging to verify in the case of complex systems such as those related to complex fluids. Therefore, our study opens up a new avenue for the exploration of stochastic stability and related issues.

When $\mu_0$ is a physical measure, Theorem \ref{thm-uniform-estimate-stability} (2) asserts its stochastic stability with respect to $\mathcal{A}$. While there are an abundance of examples asserting the stochastic instability of $\mu_0$ with respect to $\mathcal{A}$ when it fails to be physical (see Remark \ref{rem-selection-by-noise}), the converse is generally wrong as shown in the following result. In which, we construct a volume-preserving system where the normalized volume is stochastically stable but non-physical.

\begin{thmx}\label{thm-construction}
There exist a three-dimensional smooth, connected, and closed Riemannian manifold $M$ and a smooth divergence-free vector field $B:M\to TM$ satisfying the following conditions:
\begin{enumerate}
    \item $M=M_1\cup M_2$ and $M_{1}\cap M_{2}=\partial$, where $M_1$ and $M_2$ are three-dimensional smooth, connected and compact manifolds with a common boundary $\pa$,
    \item the interior $M_i^{\mathrm{o}}$ of $M_i$ for each $i=1,2$ and the boundary $\pa$ are invariant under $\vp^t$,
    \item ${\rm Vol}|_{M^{\mathrm{o}}_{i}}$ is strongly mixing for each $i=1,2$,
\end{enumerate}
where $\vp^t$ is the flow generated by $B$. Then, the following hold.
\begin{itemize}
    \item[\rm(i)] Any invariant measure of $\vp^t$ has the form of $\nu_1+\nu_2+\nu_{\pa}$, where $\nu_1$, $\nu_2$ and $\nu_{\pa}$ are invariant measures of $\vp^t$ when restricted to $M_1^{\mathrm{o}}$, $M_2^{\mathrm{o}}$ and $\pa$, respectively.

    \item[\rm(ii)] $\text{\rm Vol}$ is neither ergodic nor physical, but it is the only invariant measure of $\vp^t$ admitting a density in $W^{1,2}$. 

    \item[\rm(iii)] $\text{\rm Vol}$ is stochastically stable with respect to $\mathcal{A}$.
\end{itemize}    
\end{thmx}

Note that, in the context of the flow $\vp^{t}$ described in Theorem \ref{thm-construction}, both $\frac{\text{\rm Vol}|_{M_{1}}}{\text{\rm Vol}(M_{1})}$ and $\frac{\text{\rm Vol}|_{M_{2}}}{\text{\rm Vol}(M_{2})}$ are physical measures of $\vp^{t}$. However, it is important to note that neither of them is stochastically stable with respect to $\mathcal{A}$. This observation prompts an intriguing question since there is a prevailing belief that physical measures should demonstrate stability under small random perturbations to some extent. This issue can be attributed to a timescale problem, implying that a physical measure may exhibit stochastic stability only over certain finite timescales. The exploration of this matter as well as the same issue for dissipative systems will be the subject of future investigations.

The conclusions in Theorem \ref{thm-uniform-estimate-stability} can be enhanced if additional conditions on $u_0$, $A$ and $B$ are imposed.

%Our last result is devoted to {\rd obtaining} stronger results about   $\{u_{\ep}\}_{\ep}$ and stochastic stability under stronger conditions on $u_0$ and $A$.

%When the manifold $M$ has some particular structure, we are able to obtain uniform $W^{1,p}$-estimates for $\{u_{\ep}\}_{\ep}$ for some $p>d$ when the noise perturbations belong to some smaller class $\mathcal{A}_1\subset \mathcal{A}$, which satisfies the following additional condition except for {\bf (A1)}-{\bf (A2)}:

\begin{thmx}\label{thm-stronger-uniform-bounds}
Assume $u_0\in W^{2,p_0}$. Suppose there are vector fields $X_i:M\to TM$, $i\in \{1,\dots, n\}$ with $n\geq d$, belonging to $W^{1,p}$ for some $p>d$ such that 
\begin{itemize}
    \item [(i)] for each $i$, $\div (u_0 X_i)=0$ and $[u_0X_i, u_0B]=0$, where $[\cdot,\cdot]$ denotes the Lie bracket;
    
    %\item [(ii)] for each $i$, $\mu_0$ is an invariant measure of the flow generated by $X_i$;
    
    \item [(ii)] $\{X_i\}_{i=1}^n$ spans the tangent bundle $TM$.
\end{itemize}

Then, for each $A\in \mathcal{SA}$, 
\begin{itemize}
    \item [(1)] $\|u_{\ep}\|_{1,q}\lesssim1$ for all $q\geq1$;

    \item[(2)] if $\mu_0$ is the only invariant measure of \eqref{eqn-ode} with a density in $W^{1,2}$, then
$$
\lim_{\ep\to0}u_{\ep}=u_{0}\,\,\text{in}\,\, C^{\al},\quad\forall \al\in(0,1).
$$
\end{itemize}

\end{thmx}

Theorem \ref{thm-stronger-uniform-bounds} applies particularly to rotations on tori. In fact, for a rotation on $\T^d:=\R^d/\Z^d$, we can choose $\mu_0$ as the normalized Lebesgue measure, set $n=d$, and define $X_i:=\frac{\pa}{\pa x_i}$ for $i\in \{1,\dots, d\}$. Under these choices, the conditions in Theorem \ref{thm-stronger-uniform-bounds} are satisfied.

Consider the case when $M$ is the 2-sphere $\S^2$ equipped with the spherical coordinate $(\phi,\theta)$, where $\phi$ and $\theta$ stand for the longitude and latitude, respectively. Let $B:=\frac{\pa }{\pa \phi}$ be the vector field generating a volume-preserving flow. Straightforward calculations show that if $X$ is a divergence-free vector field satisfying $[X,B]=0$, and $B$ and $X$ are linearly independent, then $X$ must have a component $\frac{C}{\sin \phi}\frac{\pa}{\pa \theta}$ for some $C\in \R\sm \{0\}$. Consequently, $X$ does not belong to $W^{1,p}$ for any $p>2$. Hence, the conditions in Theorem \ref{thm-stronger-uniform-bounds} cannot be satisfied in this case.

\subsection{Some consequences}

We discuss two consequences about uniform estimates of $\{u_{\ep}\}_{\ep}$ derived in Theorem \ref{thm-uniform-estimate-stability}.
One concerns an optimal lower bound on the exponential convergence rate of the distribution of $X^{\ep}_t$ to the corresponding stationary measure $\mu_{\ep}$ with a good control on the sub-exponential factor. More precisely, we consider the Cauchy problem for \eqref{eqn-fpe-time} with the initial condition:
\begin{equation}\label{initial-condition}
    u(\cdot, 0)=\nu_0,
\end{equation}
where $\nu_0$ is a probability measure on $M$. Recall that for probability measures $\mu_{1},\mu_{2}$ on $M$ with $\mu_{1}\ll\mu_{2}$, the $\chi^{2}$-distance between $\mu_{1}$ and $\mu_{2}$ is defined by 
$$
\chi^2(\mu_1, \mu_{2}):=\int \left(\frac{d\mu_{1}}{d\mu_2}-1 \right)^2 d\mu_2. 
$$
We prove the following result. 

\begin{thmx}\label{thm-chi-convergence}
Let $A\in \mathcal{A}$ satisfy {\bf (A1)} with $p>d+2$. Assume $\nu_0$ has a density $v_0\in L^2$. Then, the following hold.
\begin{itemize}
    \item[(1)] For each $\ep$, \eqref{eqn-fpe-time}-\eqref{initial-condition} admits a unique solution $v^{\ep}\in V_{2,loc}^{1,0}(M\times [0,\infty))\cap \HH^{1,p}_{loc}(M\times (0,\infty))$, which is continuous and positive in $M\times (0,\infty)$ and satisfies $\int v^{\ep}(\cdot, t)= 1$ for all $t\in(0,\infty)$. 
    
    \item[(2)] There exists $C>0$ (independent of $\nu_0$) such that 
    $$
    \chi^2(\nu^{\ep}_t, \mu_{\ep})\leq e^{-C\ep^2 t}\chi^2(\nu_0, \mu_{\ep}),\quad\forall t\geq 0, \,\,0<\ep\ll1,
    $$
    where $d\nu^{\ep}_t:=v^{\ep}(\cdot, t)dx$.
\end{itemize}
\end{thmx}

The constant $C$ can be implicitly determined and depends particularly on uniform lower and upper bounds of $\{u_{\ep}\}_{\ep}$ (see Remark \ref{rem-constant-C-in-convergence-rate} and the proof of Theorem \ref{thm-chi-convergence}). Examining the case of the zero vector field $B\equiv0$, we observe that the exponential convergence rate is optimal up to constant multiplication. 

Let us consider the Fokker-Planck operator $\LL_{\ep}^{*}$ in the weighted space $L^{2}(u_{\ep}^{-1}):=L^{2}(M,u_{\ep}^{-1}d{\rm Vol})$. Then, it has a discrete spectrum and has $0$ as the principal eigenvalue with eigenspace ${\rm span}\{u_{\ep}\}$. Note that the estimate in Theorem \ref{thm-chi-convergence} (2) can be rephrased as 
$$
\|v^{\ep}(\cdot,t)-u_{\ep}\|^{2}_{L^{2}(u_{\ep}^{-1})}\leq e^{-C\ep^2 t} \|v_0-u_{\ep}\|^{2}_{L^{2}(u_{\ep}^{-1})},\quad\forall t\geq 0, \,\,0<\ep\ll1,
$$
which says in particular that the spectral gap of $\LL_{\ep}^{*}$ has a lower bound $\frac{C\ep^{2}}{2}$, but is generally much stronger. In fact, the information about the spectral gap only asserts
$$
\|v^{\ep}(\cdot,t)-u_{\ep}\|^{2}_{L^{2}(u_{\ep}^{-1})}\leq C_{\ep} e^{-C\ep^2 t} \|v_0-u_{\ep}\|^{2}_{L^{2}(u_{\ep}^{-1})},\quad\forall t\geq 0, \,\,0<\ep\ll1
$$
without a control on the growth of $C_{\ep}$.

Introducing mixing properties into the vector field $B$ can lead to an improved  convergence rate -- a phenomenon known as \emph{enhanced dissipation/relaxation} in fluid dynamics. This effect is established only in the context of divergence-free vector fields and symmetric diffusions (see e.g. \cite{MR2434887,MR4156602,MR4580539}), that is, the Fokker-Planck equation \eqref{eqn-fpe-time} reads in local coordinates 
$$
\partial_{t}u=\frac{\ep^{2}}{2}\nabla\cdot\left(A\nabla u\right)-\nabla\cdot(Bu),
$$
where $A$ is symmetric and positive definite and $B$ is divergence-free.

\medskip

The other consequence confirms the sub-exponential large deviation principle (LDP), also known as the zeroth-order WKB expansion \cite{Graham_1989}, for randomly perturbed conservative systems.

\begin{rem}\label{rem-ldp}
For randomly perturbed conservative systems, the sub-exponential LDP of stationary measures $\{\mu_{\ep}\}_{\ep}$ concerns the rigorous justification of 
$$
u_{\ep}=R_{\ep}e^{-\frac{2}{\ep^2}V}\quad\text{with}\quad \min V=0\,\,\text{and}\,\,R_{\ep}=R_{0}+o(\ep^{2}),
$$
where $V$ and $R_{\ep}$ are called the quasi-potential function and prefactor, respectively. In the situations of Theorems \ref{thm-uniform-estimate-stability} and \ref{thm-stronger-uniform-bounds}, there holds $V=0$, and hence, $R_{\ep}=u_{\ep}$ has relevant properties.

This sheds light on rigorously justifying the sub-exponential LDP of stationary measures $\{\mu_{\ep}\}_{\ep}$ for randomly perturbed dissipative systems, particularly when the global attractor of the dissipative system is a Riemannian manifold $M$, and the dissipative system, when restricted to $M$, is a conservative system in the same sense as that of \eqref{eqn-ode}. In such instances, one needs to justify 
$$
u_{\ep}=\frac{R_{\ep}}{\ep^{N-d}}e^{-\frac{2}{\ep^2}V} \quad\text{with}\quad \min V=0\,\,\text{and}\,\,R_{\ep}=R_{0}+o(\ep^{2}),
$$
where $N$ is the system's dimension and $d$ is the dimension of $M$. The dynamics of the dissipative system on $M$ for sure plays a crucial role in quantifying the limit of $u_{\ep}$.

The sub-exponential LDP for randomly perturbed dissipative systems has many applications, including classical first exit problems \cite{MR877726}, stochastic bifurcations \cite{MR0928950,Graham_1989}, stochastic populations \cite{Assaf_2017}, and the landscape and flux theory of non-equilibrium systems \cite{doi:10.1080/00018732.2015.1037068}. However, thus far, it has been only established in scenarios where the global attractor is a non-degenerate equilibrium point \cite{MR826705,MR877726,MR1055418}, owing to certain essential difficulties.   
\end{rem}

\medskip

%To see this, consider the simplest case where $B=0$ and $A_i^{\ep}$ for $i\in \{0,\dots, m\}$ are independent of $\ep$. In this setting, the generator $\LL_{\ep}=\ep^2 \left(\frac{1}{2}\sum_{i=1}^m A_i^2+A_0\right)=:\ep^2 \LL$. By rescaling the time variable in the Fokker-Planck equation \eqref{eqn-fpe-time} as $t \mapsto \frac{t}{\ep^2}$, we obtain the equation $\partial_t u = \mathcal{L}^* u$, whose solution converges to its steady state at a constant exponential rate $C > 0$. Reverting the time rescaling implies that $X_t^\varepsilon$ converges to $ \mu_\varepsilon$ at an exponential rate $C\ep^2$.

We conclude this subsection by discussing our study in one dimension and on bounded domains. The one-dimensional case is highly exceptional, as the conservative vector field $B$, which is identified as a function on $M$, must either vanish entirely or be sign-definite. We focus on the sign-definite scenario, which is evidently more interesting. This unique characteristic enables us to provide a significantly more concise proof of the one-dimensional counterpart of Theorem \ref{thm-uniform-estimate-stability}(1) and to establish uniform estimates for stationary measures for a substantially broader range of random perturbations. The corresponding findings are detailed in Section \ref{sec-1D}.

Our study extends to the setting on an open, bounded, and smooth domain $M\subset\R^{d}$, where the vector field $B$ is tangent to $\partial M$ and conservative, and the randomly perturbed dynamical system \eqref{sde-manifold} is equipped with the obliquely reflecting boundary condition. Such a setting is often used in the mathematical modeling of diffusion processes in confined environments and finds wide applications in physics, chemistry, biology and many other disciplines (see e.g. \cite{MR2676235, schuss2015}). Results analogous to Theorems \ref{thm-uniform-estimate-stability} and \ref{thm-chi-convergence} are obtained with details given in Section \ref{section-bounded-domain-case}.

%%%%%%%%%%%%%%

%The study of Kolmogorov's problem can be significantly simplified in the one-dimensional situation. Since any smooth, connected and closed one-dimensional manifold is smoothly diffeomorphic to $\S^1$, it suffices to examine the randomly perturbed conservative systems on $\S^1$. It turns out that the vector field $B$, regarded as a $2\pi$-periodic function, either vanishes everywhere or preserves its sign. Clearly, the latter deserves more attention and is our primary focus. Benefiting from the sign-preserving property of $B$ and the much simplified Sobolev embedding theorem, the proof of Theorem \ref{thm-uniform-estimate-stability} becomes notably straightforward. We include it in Section \ref{sec-1D} for the reader's reference. Apart from that, we apply the method of Berstein estimates to leverage the sign-preserving property of $B$ and achieve uniform bounds for $\{u_{\ep}\}_{\ep}$ as well as their derivatives, as outlined in Theorems \ref{thm-uniform-bounds-u_ep-1d}. A key distinction from Theorem \ref{thm-uniform-estimate-stability} lies in the novelty of this method: it allows for the potential ``blow-up" of noise perturbations, as opposed to the requirement of uniform bounds stipulated in {\bf (A1)}.

\subsection{Organization and notations}\label{subsec-notaion}
The rest of the paper is organized as follows. In Section \ref{sec-div-free}, we establish uniform estimates for the stationary measures of \eqref{sde-manifold} when the vector field $B$ is divergence-free. Section \ref{sec-main} is dedicated to proving Theorems \ref{thm-uniform-estimate-stability}-\ref{thm-chi-convergence}. When the system \eqref{eqn-ode} admits multiple invariant measures like $\mu_0$, the problem of invariant measure selection by noise is studied in Subsection \ref{subsec-selection-by-noise}. In Section \ref{sec-1D}, we address the one-dimensional case. The problems on bounded domains are treated in Section \ref{section-bounded-domain-case}. Appendix \ref{app-stationary-measure} contains results regarding the existence and uniqueness of stationary measures for SDEs with less regular coefficients as well as the well-posedness of corresponding Fokker-Planck equations. Additionally, in Appendix \ref{app-formulas}, we collect some commonly used formulas from calculus on manifolds.

\medskip

{\bf Notation.} The following list of notations are used throughout this paper.
\begin{itemize}

    \item $d_*:=\frac{d}{d-2}$ if $d\geq 3$ and fix any $d_*\in(1,\infty)$ if $d=2$.

   \item $\N$ denotes the set of positive integers, and $\N_0:=\{0\}\cup \N$.

    \item The (vector-valued) spaces $L^{p}(M)$, $W^{k,p}(M)$, and $C^{\al}(M)$ are written as $L^{p}$, $W^{k,p}$, and $C^{\al}$, respectively. 
    
    \item The usual $L^{p}$-norm and $W^{k,p}$-norm are denoted by $\|\cdot\|_{p}$ and $\|\cdot\|_{k,p}$, respectively.

    \item Let $I\subset [0,\infty)$ be an interval.  Denote $\mathbb{H}^{1,p}_0(M\times I)=L^p(I,W_0^{1,p})$ and by $\mathbb{H}^{-1,p'}(M\times I)$ the dual space of $\mathbb{H}^{1,p}_0(M\times I)$.  
    
    \item Denote by $\HH^{1,p}_{loc}(M\times I)$ the space of measurable functions $u$ on $M\times I$ such that $\eta u\in \mathbb{H}_0^{1,p}(M\times I)$ and $\pa_t(\eta u)\in \mathbb{H}^{-1,p}$ for every $\eta \in C_0^{\infty}(M\times I)$.

    \item Denote by $V_{2,loc}^{1,0}(M\times I)$ the space of measurable functions $u$ on $M\times I$ such that $t\mapsto u(\cdot,t)$ is continuous under the $L^2$-norm and $\eta u\in L^2(I,W^{1,2})$ for any $\eta\in C_0^{\infty}(M\times I)$.

    \item Denote by $C^{i,j}(M\times I)$ the space of continuous functions $u$ on $M\times I$ such that $u$ is continuously differentiable up to the $i$-th order in $x$ and up to the $j$-th order in $t$.

    \item For constants $\al(\ep)$ and $\beta(\ep)$ indexed by $\ep$, we write
    \begin{itemize}
        \item $\al(\ep)\lesssim\beta(\ep)$ (resp. $\al(\ep)\gtrsim\beta(\ep)$) if there is a positive constant $C$, independent of $\ep$, such that $\al(\ep)\leq C\beta(\ep)$ for all $\ep$ (resp. $\al(\ep)\geq C\beta(\ep)$ for all $\ep$);

        \item $\al(\ep)\approx\beta(\ep)$ if $\al(\ep)\lesssim\beta(\ep)$ and $\al(\ep)\gtrsim\beta(\ep)$.
    \end{itemize}

    \item For $f\in C^{0}$, $\min_{M}f$ and $\max_{M}f$ are written as $\min f$ and $\max f$, respectively.

    \item The integral $\int_{M} fd{\rm Vol}$ is written as $\int f$, and the integral $\int_{M}fd\mu$, where $\mu$ is a measure on $M$, is written as $\int fd\mu$.

    \item Einstein's summation convention is used throughout the paper unless otherwise specified. 
\end{itemize}

%%%%%%%%%%%%%%%%%%%%%%%%%%%%%%%%%%%%%%%%
%%%%%%%%%%%%%%%%%%%%%%%%%%%%%%%%%%%%%%%%
%%%%%%%%%%%%%%%%%%%%%%%%%%%%%%%%%%%%%%%%

\section{\bf Uniform estimates in the divergence-free case}\label{sec-div-free}

This section is devoted to uniform estimates for stationary measures of \eqref{sde-manifold} under different assumptions on the vector field $B$. The corresponding result serves as a crucial step in the proof of our main results, and is of independent interest. 

We make the following assumptions on the vector field $B$.

\medskip

\begin{itemize}
    \item[\bf(A)$_{B}$] $B\in W^{1,p_0}$ for some $p_0>\frac{d}{2}$ and $\div B=0$ {\rm{Vol}}-a.e.
\end{itemize}

\medskip

That is, $B$ is divergence-free and has weak regularity. It is well known that the flow generated by $B$, if exists, is volume-preserving. Unfortunately, $B\in W^{1,p_0}$ is insufficient for \eqref{eqn-ode} to generate such a flow. But, this causes no trouble at all since the dynamical implication is off the table in this section.

Note that $B\in L^{p_1}$ with some $p_1>d$ thanks to the Sobolev embedding theorem. Hence, Theorem \ref{lem-stationary-measure} applies and yields that if $A\in\mathcal{A}$, then for each $\ep$, the SDE \eqref{sde-manifold} admits a unique stationary measure $\mu_{\ep}$ having a positive density $u_{\ep}\in W^{1,p_*}$, where $p_*:=\min\{p_1, p\}$ and $p$ is as assumed in {\bf (A1)}. Moreover, 
\begin{equation}\label{eqn-fpe}
-\frac{\ep^2}{2}\int A^{\ep}_i f\left[A^{\ep}_i u_{\ep}+(\div A^{\ep}_i) u_{\ep}\right]+\int\left(\ep^{2}A^{\ep}_0f+Bf\right)u_{\ep}=0,\quad \forall f\in W^{1,2}.
\end{equation}

The following two theorems are the main results in this section. 

\begin{thm}\label{thm-lower-bound-div-free}
Assume {\bf(A)}$_{B}$. For any $A\in \mathcal{A}$, there hold
\begin{equation*}
\|u_{\ep}\|_{1,2}\lesssim1\quad\text{and}\quad 1\lesssim \min u_{\ep}\leq \max u_{\ep}\lesssim1.
\end{equation*}
\end{thm}

\begin{thm}\label{thm-stronger-uniform-bounds-div-free}
Assume {\bf (A)$_B$} with $p_0>d$. Suppose there are vector fields $X_i:M\to TM$, $i\in \{1,\dots, n\}$ with $n\geq d$, belonging to $W^{1,p}$ for some $p>d$ such that 
\begin{itemize}
    \item [(i)] for each $i$, $\div X_i=0$ and $[X_i, B]=0$;
    
    \item [(ii)] $\{X_i\}_{i=1}^{n}$ spans the tangent bundle $TM$.
\end{itemize}
Then, for any $A\in \SS\mathcal{A}$, there holds $\|u_{\ep}\|_{1,q}\lesssim1$ for all $q\geq1$.
\end{thm}

%%%%%%%%%%%%%%%%%%%%%%%%%%%%%%%%

The rest of this subsection is devoted to the proof of Theorems \ref{thm-lower-bound-div-free} and \ref{thm-stronger-uniform-bounds-div-free}. To proceed with the proof of Theorem \ref{thm-lower-bound-div-free}, we need two lemmas. From now on, we will frequently use some formulas relevant to calculus on manifolds. The reader is directed to Appendix \ref{app-formulas} for details.

\begin{lem}\label{lem-uniform-L2-gradient}
Assume {\bf(A)}$_{B}$. For any $A\in \mathcal{A}$, we have $\|u_{\ep}\|_{1,2}\lesssim1$.
\end{lem}

\begin{proof}
Let $A\in \mathcal{A}$. Taking $f=u_{\ep}$ in \eqref{eqn-fpe} yields
\begin{equation}\label{eqn-2023-01-08-1}
-\frac{\ep^2}{2} \int A^{\ep}_iu_{\ep}[A^{\ep}_i u_{\ep}+(\div A^{\ep}_i) u_{\ep}] +\int(\ep^2 A^{\ep}_0 u_{\ep}+B u_{\ep})u_{\ep}=0.
\end{equation}

Since $\div B=0$ by {\bf(A)}$_{B}$, we derive from the divergence theorem that 
\begin{equation*}
\int (Bu_{\ep}) u_{\ep}=\frac{1}{2}\int Bu_{\ep}^{2}=-\frac{1}{2}\int \left(\div B \right) u_{\ep}^{2}=0,
\end{equation*}
which together with \eqref{eqn-2023-01-08-1} leads to 
\begin{equation*}%\label{eqn-Aug-1-2}
\begin{split}
    \sum_{i=1}^m\int |A^{\ep}_i u_{\ep}|^2 +\int(A^{\ep}_i u_{\ep})(\div A^{\ep}_i) u_{\ep}-2\int (A^{\ep}_0 u_{\ep})u_{\ep}=0.
\end{split}  
\end{equation*}

The assumption {\bf (A1)} and the Sobolev embedding theorem ensure $\max_i \|A^{\ep}_i\|_{\infty}\lesssim1$ and  $ \|A^{\ep}_0\|_p+\max_i \|\div A^{\ep}_i\|_p\lesssim1$. Hence, we apply H\"older's inequality to deduce that
\begin{equation*}
\begin{split}
\sum_{i=1}^m \int |A^{\ep}_i u_{\ep}|^2% &=-\int(A^{\ep}_i u_{\ep})(\div A^{\ep}_i) u_{\ep}+2\int (A^{\ep}_0 u_{\ep})u_{\ep}\\
   &\lesssim \|\nabla u_{\ep}\|_2\|u_{\ep}\|_r\left( \sum_{i=1}^m\|\div A^{\ep}_i\|_p +\|A^{\ep}_0\|_p\right)\lesssim \|\nabla u_{\ep}\|_2\|u_{\ep}\|_r,
\end{split}
\end{equation*}
where $r:=\left( \frac{1}{2}-\frac{1}{p}\right)^{-1}\in (2, 2d_*)$. As $\sum_{i=1}^m \int |A^{\ep}_i u_{\ep}|^2\gtrsim \|\nabla u_{\ep}\|^2_2$ ensured by {\bf (A2)}, we arrive at $\|\nabla u_{\ep}\|_2\lesssim \|u_{\ep}\|_r$.  The Sobolev embedding theorem then leads to 
\begin{equation}\label{estimate-20230206}
    \|u_{\ep}\|_{2d_*}\lesssim \|u_{\ep}\|_2+\|\nabla u_{\ep}\|_2\lesssim  \|u_{\ep}\|_r,
\end{equation}
where we used $\|u_{\ep}\|_2\lesssim \|u_{\ep}\|_r$ in the second inequality.  

Note that $\|u_{\ep}\|_r\leq \|u_{\ep}\|^{\al}_1\|u_{\ep}\|^{1-\al}_{2d_*}$ (by interpolation) with $\al:=\left(\frac{1}{d}-\frac{1}{p}\right)\left( \frac{1}{2}+\frac{1}{d}\right)^{-1}$. The fact $\|u_{\ep}\|_{1}=1$ results in $\|u_{\ep}\|_r\leq \|u_{\ep}\|^{1-\al}_{2d_*}$. It then follows from \eqref{estimate-20230206} that $\|u_{\ep}\|_{2d_*}\lesssim 1$, and thus, $\|u_{\ep}\|_r\lesssim 1$. The desired conclusion follows readily from the second inequality in \eqref{estimate-20230206}.
\end{proof}

\begin{lem}\label{lem-upper-bound}
Assume {\bf(A)}$_{B}$. For any $A\in \mathcal{A}$, the following hold.
\begin{enumerate}
\item[\rm(1)] $\max u_{\ep}\lesssim1$;

\item[\rm(2)] For any $\ga>0$, we have $\min u_{\ep}\gtrsim  \|u^{-1}_{\ep}\|^{-1}_{\ga}$.
\end{enumerate}
\end{lem}
\begin{proof}
We first establish some estimates for $u_{\ep}$.  Recall that $u_{\ep}\in W^{1,p_*}$ and  $u_{\ep},\, u_{\ep}^{-1}\in L^{\infty}$ due to the Sobolev embedding theorem.

Let $q\in\R\setminus\{0,1\}$. In the following, for any two constants $\al(\ep,q)$ and $\be(\ep,q)$ indexed by $\ep$ and $q$, we write $\al(\ep,q) \lesssim \be(\ep,q)$ to imply the existence of $C$ (independent of $\ep$ and $q$) such that $\al(\ep,q) \leq C \be(\ep,q)$ for all $(\ep,q)$. Setting $f=u^{q-1}_{\ep}\in W^{1,2}$ in \eqref{eqn-fpe} gives rise to
\begin{equation*}
-\frac{\ep^2}{2}\int (A^{\ep}_i u^{q-1}_{\ep})\left[A^{\ep}_iu_{\ep}+(\div A^{\ep}_i) u_{\ep}\right]+\int(\ep^2 A^{\ep}_0 u_{\ep}^{q-1}+B u_{\ep}^{q-1}) u_{\ep}=0.
\end{equation*}
Since $\div B=0$ by {\bf(A)}$_{B}$, we apply the divergence theorem to find 
$$
\int (Bu^{q-1}_{\ep}) u_{\ep}=\frac{1}{q}\int B u_{\ep}^{q}=-\frac{1}{q}\int (\div B) u^q_{\ep}=0.
$$
Therefore, 
\begin{equation*}
    \begin{split}
        0&=\int (A^{\ep}_i u^{q-1}_{\ep})\left[A^{\ep}_iu_{\ep}+(\div A^{\ep}_i) u_{\ep}\right]-2\int(A^{\ep}_0 u_{\ep}^{q-1}) u_{\ep}\\
        % &=\int A_i (u_{\ep} A_i u^{p-1}_{\ep})-\int A_i u_{\ep}^{p-1} A_i u_{\ep}+(p-1)\intu_{\ep}^{p-1}A_0 u_{\ep}\\
        % &=-\int (\div A_i) u_{\ep} A_i u_{\ep}^{p-1} -\int A_i u_{\ep}^{p-1} A_i u+(p-1)\intu_{\ep}^{p-1}A_0 u_{\ep}\\
        &=(q-1)\sum_{i=1}^m \int u_{\ep}^{q-2} |A^{\ep}_i u_{\ep}|^2 +(q-1) \int (\div A^{\ep}_i) u_{\ep}^{q-1} A^{\ep}_i u_{\ep}-2(q-1)\int u_{\ep}^{q-1}A^{\ep}_0 u_{\ep}.
    \end{split}
\end{equation*}

Since $A\in \mathcal{A}$, we apply H\"older's inequality to derive 
\begin{equation*}
\begin{split}
\int u_{\ep}^{q-2}|\nabla u_{\ep}|^2\lesssim  \sum_{i=1}^m \int u_{\ep}^{q-2} |A^{\ep}_i u_{\ep}|^2&=-\int (\div A^{\ep}_i) u_{\ep}^{q-1} A^{\ep}_i u_{\ep}+2\int u_{\ep}^{q-1}A^{\ep}_0 u_{\ep}\\
&\lesssim \left(\int u_{\ep}^{q-2}|\nabla u_{\ep}|^2\right)^{\frac{1}{2}} \|u_{\ep}^{\frac{q}{2}}\|_r \left(\sum_{i=1}^m\|\div A^{\ep}_i\|_{p}+\|A^{\ep}_0\|_{p}\right)\\
&\lesssim \left(\int u_{\ep}^{q-2}|\nabla u_{\ep}|^2\right)^{\frac{1}{2}} \|u_{\ep}^{\frac{q}{2}}\|_r,
\end{split}
\end{equation*}
where $r:=\left( \frac{1}{2}-\frac{1}{p}\right)^{-1}\in (2, 2d_*)$. As a result,  $\int |\nabla u_{\ep}^{\frac{q}{2}}|^2 =\frac{q^2}{4}\int  u_{\ep}^{q-2}|\nabla u_{\ep}|^2\lesssim q^2 \|u_{\ep}^{\frac{q}{2}}\|^2_r$, which together with the Sobolev embedding theorem leads to
\begin{equation*}
\begin{split}
\| u_{\ep}^{\frac{q}{2}}\|^2_{2d_*}\lesssim \| u_{\ep}^{\frac{q}{2}}\|^2_{2}+\|\nabla u_{\ep}^{\frac{q}{2}}\|^2_{2}\lesssim \| u_{\ep}^{\frac{q}{2}}\|^2_{2}+ q^2\|u_{\ep}^{\frac{q}{2}}\|^2_{r}.
\end{split}
\end{equation*}

Noting that $\| u_{\ep}^{\frac{q}{2}}\|^2_{2}\lesssim \| u_{\ep}^{\frac{q}{2}}\|^2_{r}$, we find $C_{*}>0$, independent of $\ep$ and $q$, such that $\| u^{\frac{q}{2}}_{\ep}\|^{2}_{2d_*}\leq C_* q^2 \| u_{\ep}^{\frac{q}{2}}\|^2_{r}$. As $\| u^{\frac{q}{2}}_{\ep}\|_{r}\leq \| u^{\frac{q}{2}}_{\ep}\|^{\al}_{2}\| u^{\frac{q}{2}}_{\ep}\|^{1-\al}_{2d_*}$ (by interpolation) with $\al:=1-\frac{d}{p}$, we deduce that 
\begin{equation}\label{eqn-2023-01-10-5}
    \left\| u_{\ep}\right\|_{qd_*}\leq C_*^{\frac{1}{\al q}} q^{\frac{2}{\al q}} \left\| u_{\ep}\right\|_{q} \quad \text{if } q>0,\,\,q\neq1,
\end{equation}
and  
\begin{equation}\label{eqn-2023-01-10-5-1}
    \| u^{-1}_{\ep}\|_{|q|d_*}\leq C^{\frac{1}{\al|q|}}_* |q|^{\frac{2}{\al |q|}}\| u^{-1}_{\ep}\|_{|q|}\quad \text{if } q<0.
\end{equation}
% \begin{equation}\label{eqn-2023-01-10-5}
% \begin{gathered}
% \left\| u_{\ep}\right\|^{p}_{pd_*}\leq C_3\left(\left\| u_{\ep}\right\|^{p}_{p}+\left\|\nabla u_{\ep}^{\frac{p}{2}}\right\|^{2}_{2}\right)\leq C_3\left(1+\frac{C_1 p^2}{4} \right)\left\| u^{\frac{p}{2}}_{\ep}\right\|^2_{q}\quad \text{if } p>0,\,\,p\neq1,\\
% \left\| u^{-1}_{\ep}\right\|^{|p|}_{|p|d_*}\leq C_3\left(\left\| u^{-1}_{\ep}\right\|^{|p|}_{|p|}+\left\|\nabla u_{\ep}^{\frac{p}{2}}\right\|^{2}_{2}\right)\leq C_3\left(1+\frac{C_2 p^2}{4} \right)\left\| u^{-1}_{\ep}\right\|^{|p|}_{|p|}\quad \text{if } p<0.
% \end{gathered}
% \end{equation}

\medskip

Now, we prove the results.

\medskip

(1) Setting $q=2d_*^k$ in \eqref{eqn-2023-01-10-5} for each $k\in \N_{0}$, we see that  
$$
\left\| u_{\ep}\right\|_{2d_*^{k+1}}\leq C_*^{\frac{1}{2\al d_*^k}} (2d_*^k)^{\frac{1}{\al d_*^k}}\left\| u_{\ep}\right\|_{2d_*^{k}},\quad\forall k\in\N_{0}.
$$
By iteration, 
$$
\|u_{\ep}\|_{\infty}=\lim_{k\to \infty}\| u_{\ep}\|_{2d^k_*}\leq (4C_*)^{\frac{1}{2\al}\sum_{k=0}^{\infty} \frac{1}{d^k_*}}d_*^{\frac{1}{\al}\sum_{k=0}^{\infty} \frac{k}{d^k_*}}\| u_{\ep}\|_{2}\lesssim \| u_{\ep}\|_{2}\lesssim \|u_{\ep}\|^{\frac{1}{2}}_{\infty} \|u_{\ep}\|^{\frac{1}{2}}_{1}.
$$ 
The result then follows from the fact $\|u_{\ep}\|_{1}=1$.

\medskip

(2) For $\ga>0$, we set $q=-\ga d_*^k$ in \eqref{eqn-2023-01-10-5-1} for each $k\in \N_{0}$ to find that 
$$
 \left\| u^{-1}_{\ep}\right\|_{\ga d_*^{k+1}}\leq C_*^{\frac{1}{\al\ga d_*^k}}\left(\ga d_*^k\right)^{\frac{2}{\al\ga d_*^k}} \left\| u^{-1}_{\ep}\right\|_{\ga d_*^k},\quad\forall k\in\N_0.
$$
It follows that 
\begin{equation*}
\|u^{-1}_{\ep}\|_{\infty}\leq (C_* \ga^2)^{\frac{1}{\al\ga}\sum_{k=0}^{\infty} \frac{1}{ d^k_*}} d_*^{\frac{2}{\al\ga}\sum_{k=0}^{\infty} \frac{k}{ d^k_*}} \|u^{-1}_{\ep}\|_{\ga},
\end{equation*}
that is, $\min u_{\ep}\gtrsim \|u^{-1}_{\ep}\|^{-1}_{\ga}$. This completes the proof.
\end{proof}

We are ready to prove Theorem \ref{thm-lower-bound-div-free}.

\begin{proof}[Proof of Theorem \ref{thm-lower-bound-div-free}]
%Since Sobolev compact embedding theorem ensures that any bounded sequences $W^{1,2}(\T^d)$ is precompact in $L^1(\T^d)$, we see from Lemma \ref{lem-uniform-L2-gradient} and Lemma \ref{lem-upper-bound} that $\{u_{\ep}\}_{\ep}$ is precompact in $L^1(\T^d)$. The positive lower bound of the limit(s) follows immediately once we prove the positive uniform lower bound of $\{u_{\ep}\}_{\ep}$.

Given Lemmas \ref{lem-uniform-L2-gradient} and \ref{lem-upper-bound}, it remains to prove 
\begin{equation}\label{bounded-u_ep-inverse-norm}
\|u_{\ep}^{-1}\|_{\ga}\lesssim1\quad\text{for some}\quad\ga>0.    
\end{equation}

To verify \eqref{bounded-u_ep-inverse-norm}, we set $v_{\ep}:=\ln u_{\ep}-\int  \ln u_{\ep} $ and break the proof into four steps.

\medskip

\paragraph{\bf Step 1} We show that for each $\eta\in C^2(\R)$,
\begin{equation}\label{eqn-log u}
\sum_{i=1}^m\int \left[\eta'(v_{\ep})-\eta''(v_{\ep})\right] |A^{\ep}_i v_{\ep}|^2 =\int \left[\eta''(v_{\ep})-\eta'(v_{\ep})\right] \left[(\div A^{\ep}_i)A^{\ep}_i v_{\ep} -2A^{\ep}_0v_{\ep}\right].
\end{equation}

Recall that $u_{\ep}\in W^{1,p_*}$ and  $u_{\ep},\, u^{-1}_{\ep}\in L^{\infty}$ due to the Sobolev embedding theorem.  Then, for any $g\in W^{1,2}$, we are able to set $f=\frac{g}{u_{\ep}}\in W^{1,2}$ in \eqref{eqn-fpe} to derive
\begin{equation}\label{eqn-Aug-29}
    \begin{split}
        0&=-\frac{\ep^2}{2}\int  A^{\ep}_i \frac{g}{u_{\ep}}[A^{\ep}_i u_{\ep}+(\div A^{\ep}_i) u_{\ep}]+\int \left(\ep^2 A^{\ep}_0\frac{g}{u_{\ep}}+B\frac{g}{u_{\ep}}\right) u_{\ep}\\
        &=-\frac{\ep^2}{2}\int \left(\frac{A^{\ep}_i g}{u_{\ep}}-\frac{g}{u^2_{\ep}}A^{\ep}_i u_{\ep}\right)[A^{\ep}_i u_{\ep}+(\div A^{\ep}_i) u_{\ep}]+\ep^2\int \left(\frac{A^{\ep}_0 g}{u_{\ep}}-\frac{g}{u^2_{\ep}}A^{\ep}_0 u_{\ep}\right) u_{\ep}\\
        &\quad+\int  \left(\frac{B g}{u_{\ep}}-\frac{g}{u^2_{\ep}}B u_{\ep}\right) u_{\ep}\\
        &=-\frac{\ep^2}{2}\int  \left(A^{\ep}_i g A^{\ep}_i v_{\ep} -g \sum_{i=1}^d|A^{\ep}_i v_{\ep}|^2\right)-\frac{\ep^2}{2}\int (\div A^{\ep}_i)\left(A^{\ep}_i g- gA^{\ep}_i v_{\ep}\right)\\
        &\quad +\ep^2\int  A^{\ep}_0 g-gA^{\ep}_0v_{\ep}+\int B g-gBv_{\ep}.
    \end{split}
\end{equation}

Clearly, $v_{\ep}\in W^{1,p_*}$. For each $\eta\in C^2(\R)$, we take $g=\eta'(v_{\ep})\in W^{1,2}$ and then apply the divergence theorem to find from $\div B=0$ (by {\bf (A)$_B$}) that $\int B\eta'(v_{\ep})=0$
and $\int  \eta'(v_{\ep}) Bv_{\ep}  =\int  B \eta(v_{\ep})=0$. As a result, the equality \eqref{eqn-log u} follows from \eqref{eqn-Aug-29}.

\medskip

\paragraph{\bf Step 2} We show 
\begin{equation}\label{eqn-v_ep-uniform-estimates}
D_0:=\sup_{\ep}\left( \|v_{\ep}\|^2_2+\|\nabla v_{\ep}\|^2_2\right)<\infty.
\end{equation}

Taking $\eta(t)=t$ in \eqref{eqn-log u} yields $\sum_{i=1}^m\int  |A^{\ep}_i v_{\ep}|^2=-\int (\div A^{\ep}_i) A^{\ep}_i v_{\ep}-2A^{\ep}_0 v_{\ep}$. By {\bf (A1)} and the Sobolev embedding theorem, there hold $\max_i \|A^{\ep}_i\|_{\infty}\lesssim1$ and $\|A^{\ep}_0\|_p+\max_i \|\div A^{\ep}_i\|_p\lesssim1$. We apply H\"{o}lder's inequality to find for $r:=\left( \frac{1}{2}-\frac{1}{p}\right)^{-1}\in (2, 2d_*)$ that
\begin{equation*}
    \sum_{i=1}^m\int  |A^{\ep}_i v_{\ep}|^2\lesssim \|\nabla v_{\ep}\|_2 \|1\|_r  \left(\sum_{i=1}^m\|\div A^{\ep}_i\|_p + \|A^{\ep}_0\|_p \right)\lesssim \|\nabla v_{\ep}\|_2.
\end{equation*}
By {\bf (A2)}, $\|\nabla v_{\ep}\|^2\lesssim \sum_{i=1}^m\int  |A^{\ep}_i v_{\ep}|^2\lesssim \|\nabla v_{\ep}\|_2$, and thus,  $\|\nabla v_{\ep}\|_2\lesssim1$. As $\int  v_{\ep}=0$, we apply Poincar\'e inequality to find $\|v_{\ep}\|_2\lesssim1$, and hence, \eqref{eqn-v_ep-uniform-estimates}.

\medskip 

Before proceeding with the proof, we mention that constants $C_1,\dots,C_6$, which appear in the following {\bf Step 3} and {\bf Step 4}, are independent of $\ep$ and $q\geq1$ (to be included in {\bf Step 3}). Moreover, for any two constants $\al(\ep,q)$ and $\be(\ep,q)$ indexed by $\ep$ and $q$, we write $\al(\ep,q) \lesssim \be(\ep,q)$ to imply the existence of $C$ (independent of $\ep$ and $q$) such that $\al(\ep,q) \leq C \be(\ep,q)$ for all $(\ep,q)$. 

\medskip

\paragraph{\bf Step 3} We claim the existence of $D_1,\,D_2>0$, independent of $\ep$,  such that 
\begin{equation}\label{eqn-pre-iteration}
\|v_{\ep}\|_{2qd_*}\leq D_1 q D_2^{\frac{1}{2q}} + D_2^{\frac{1}{2q}} q^{\frac{2}{\al q}} \|v_{\ep}\|_{2q}, \quad \forall q\geq 1,
\end{equation}
where $\al:=1-\frac{d}{p}\in (0,1)$.

Let $q\geq 1$. Setting $\eta(t)=\frac{1}{2q+1}t^{2q+1}$ in \eqref{eqn-log u} gives rise to 
\begin{equation}\label{eqn-2023-01-10-1}
\begin{split}
\sum_{i=1}^m\int  v_{\ep}^{2q} |A^{\ep}_i v_{\ep}|^2&=2q \sum_{i=1}^m\int  v_{\ep}^{2q-1} |A^{\ep}_i v_{\ep}|^2 + 2q\int  v_{\ep}^{2q-1} \left[(\div A^{\ep}_i) A^{\ep}_i v_{\ep} -2A^{\ep}_0 v_{\ep}\right]\\
&\quad-\int v_{\ep}^{2q} \left[(\div A^{\ep}_i) A^{\ep}_i v_{\ep} -2A^{\ep}_0 v_{\ep}\right].
\end{split}
\end{equation}

Applying Young's inequality leads to 
\begin{equation*}
|v_{\ep}|^{2q-1}\leq \frac{2q-1}{2q} v_{\ep}^{2q} \de^{\frac{2q}{2q-1}}+\frac{1}{2q} \de^{-2q},\quad\forall \de>0.
\end{equation*}
Setting $\de=[2(2q-1)]^{\frac{1-2q}{2q}}$ yields $|v_{\ep}|^{2q-1}\leq \frac{1}{4q} v_{\ep}^{2q}+\frac{1}{2q} [2(2q-1)]^{2q-1}$, and hence, 
\begin{equation*}
\begin{split}
2q \sum_{i=1}^m\int  v_{\ep}^{2q-1} |A^{\ep}_i  v_{\ep}|^2 
&\leq \frac{1}{2}\sum_{i=1}^m\int  v_{\ep}^{2q} |A^{\ep}_i  v_{\ep}|^2 +[2(2q-1)]^{2q-1}\|A^{\ep}_i\|_{\infty}\int |\nabla v_{\ep}|^2\\
&\leq \frac{1}{2}\sum_{i=1}^m\int  v_{\ep}^{2q} |A^{\ep}_i  v_{\ep}|^2+C_1 (4q)^{2q-1},
\end{split}
\end{equation*}
where we used \eqref{eqn-v_ep-uniform-estimates} in the second inequality. Inserting this into \eqref{eqn-2023-01-10-1} results in 
\begin{equation}\label{eqn-2023-01-10-2}
\begin{split}
    \frac{1}{2}\sum_{i=1}^m\int  v_{\ep}^{2q} |A^{\ep}_i v_{\ep}|^2 &\leq C_1 (4q)^{2q-1}+2q\int  v_{\ep}^{2q-1} \left|(\div A^{\ep}_i) A^{\ep}_i v_{\ep} -2A^{\ep}_0 v_{\ep}\right|\\
    &\quad-\int v_{\ep}^{2q} \left|(\div A^{\ep}_i) A^{\ep}_i v_{\ep} -2A^{\ep}_0 v_{\ep}\right|.
\end{split}
\end{equation}

For the second and third terms on the right-hand side of \eqref{eqn-2023-01-10-2}, we apply H\"older's inequality and then Young's inequality to find
\begin{equation}\label{eqn-2023-01-10-3}
\begin{split}
&2q\int  v_{\ep}^{2q-1} \left|(\div A^{\ep}_i) A^{\ep}_iv_{\ep} -2A^{\ep}_0 v_{\ep}\right|\\
&\qquad\qquad\leq 2q (\max_i \|A^{\ep}_i\|_{\infty}+1) \left(\sum_{i=1}^m\|\div A^{\ep}_i\|_p +  2\|A^{\ep}_0\|_p\right)\|v_{\ep}^{q-1}\|_{r} \|v_{\ep}^q \nabla v_{\ep}\|_2\\
&\qquad\qquad \leq \frac{\la}{8}\|v_{\ep}^q \nabla v_{\ep}\|_2^2+C_2 q^2 \|v_{\ep}^{q-1}\|_r^2,
\end{split}
\end{equation}
where $r$ is defined in {\bf Step 2}, and 
\begin{equation}\label{eqn-2023-01-10-4}
\begin{split}
&\int v_{\ep}^{2q}\left|(\div A^{\ep}_i) A^{\ep}_i v_{\ep} -2A^{\ep}_0 v_{\ep}\right|\\
& \qquad\qquad \leq (1+\max_i \|A^{\ep}_i\|_{\infty})\left(\sum_{i=1}^m\|\div A^{\ep}_i\|_p +2\|A^{\ep}_0\|_p\right) \|v_{\ep}^{q}\|_{r} \|v_{\ep}^q \nabla v_{\ep}\|_2\\
&\qquad\qquad \leq \frac{\la}{8}\|v_{\ep}^q \nabla v_{\ep}\|^2_2+ C_3\|v_{\ep}^{q}\|^2_{r},
\end{split}
\end{equation}
where $\la>0$ is the constant appearing in {\bf (A2)}. Inserting \eqref{eqn-2023-01-10-3} and \eqref{eqn-2023-01-10-4} back into \eqref{eqn-2023-01-10-2}, we derive from {\bf (A2)} that
\begin{equation*}
\begin{split}
    &\|v_{\ep}^{q}\nabla v_{\ep}\|_{2}^{2}\lesssim (4q)^{2q-1}+q^2 \|v_{\ep}^{q-1}\|_r^2+\|v_{\ep}^{q}\|_r^2.
\end{split}
\end{equation*}
According to H\"older's inequality and Young's inequality, we find
$$
\|v_{\ep}^{q-1}\|^2_r\lesssim \|v_{\ep}^{q}\|^{\frac{2(q-1)}{q}}_r\lesssim \frac{q-1}{q}\|v_{\ep}^q\|_r^2+\frac{1}{q},
$$
and hence,
\begin{equation}\label{eqn-2023-03-15-1}
\begin{split}
\|v_{\ep}^{q}\nabla v_{\ep}\|_{2}^{2} \lesssim (4q)^{2q-1}+(q^2-q+1)\|v_{\ep}^q\|^2_r+q. 
\end{split}
\end{equation}
As $v_{\ep}^{2q-2}\leq \frac{q-1}{q} v_{\ep}^{2q}+\frac{1}{q}$ by Young's inequality, we deduce
\begin{equation*}
        \|v_{\ep}^{q-1}\nabla v_{\ep}\|_{2}^{2}\leq \frac{q-1}{q}\|v_{\ep}^{q}\nabla v_{\ep}\|_{2}^{2} +\frac{D_0}{q},
\end{equation*}
where $D_0$ is given in \eqref{eqn-v_ep-uniform-estimates}. This together with \eqref{eqn-2023-03-15-1} leads to 
$$
\|v_{\ep}^{q-1}\nabla v_{\ep}\|_{2}^{2}\lesssim (4q)^{2q}+ q^2 \|v_{\ep}\|^{2q}_{qr}.
$$
It then follows from the equality $\nabla |v_{\ep}|^{q}=q|v_{\ep}|^{q-2} v_{\ep} \nabla v_{\ep}$ that
$$
\|\nabla |v_{\ep}|^{q}\|_{2}^{2}\lesssim(4q)^{2q+2}+ q^4 \|v_{\ep}\|^{2q}_{qr}.
$$
Since $\|v_{\ep}\|_{2qd_*}^{2q}\lesssim \|\nabla |v_{\ep}|^q\|_{2}^{2}+ \|v_{\ep}\|_{2q}^{2q}$ by the Sobolev embedding theorem, we arrive at
\begin{equation}\label{eqn-Aug-2-1}
\begin{split}
\|v_{\ep}\|_{2qd_*}^{2q} \leq C_4 (4q)^{2q+2}+C_4 q^4\|v_{\ep}\|_{qr}^{2q}+C_4\|v_{\ep}\|_{2q}^{2q}.
\end{split}
\end{equation}

Noting that applications of the interpolation inequality and Young's inequality yield
\begin{equation*}
\begin{split}
    C_4 q^4\|v_{\ep}\|_{qr}^{2q}\leq C_4 q^4\|v_{\ep}\|_{2q}^{2q\al} \|v_{\ep}\|_{2qd_*}^{2q(1-\al)}\leq \al \left(C_4 q^4\right)^{\frac{1}{\al}}\|v_{\ep}\|_{2q}^{2q}+(1-\al)\|v_{\ep}\|_{2qd_*}^{2q},
\end{split}
\end{equation*}
we substitute the above inequality into \eqref{eqn-Aug-2-1} to find
\begin{equation*}
\begin{split}
   \|v_{\ep}\|_{2qd_*}^{2q}\leq \frac{1}{\al}C_4 (4q)^{2q+2} +\frac{1}{\al}\left[C_4+\al \left(C_4 q^4\right)^{\frac{1}{\al}}\right]\|v_{\ep}\|_{2q}^{2q}\leq C_5 (4q)^{2q+2} +C_5 q^{\frac{4}{\al}} \|v_{\ep}\|_{2q}^{2q}.
\end{split}
\end{equation*}
Taking the $2q$-th root of both sides results in
$$
\|v_{\ep}\|_{2qd_*}\leq (C_5)^{\frac{1}{2q}} (4q)^{\frac{q+1}{q}}+(C_5)^{\frac{1}{2q}} q^{\frac{2}{\al q}}\|v_{\ep}\|_{2q}, 
$$
giving rise to \eqref{eqn-pre-iteration}. 

\medskip 

\paragraph{\bf Step 4} We finish the proof. Set $q=d_*^k$ and $I^{\ep}_k:=\|v_{\ep}\|_{2d_*^k}$ for $k\in\N_{0}$. It follows from \eqref{eqn-pre-iteration} that
\begin{equation*}
I^{\ep}_{k+1}\leq D_1 d_*^k D_2^{\frac{1}{2d_*^k}} +d_*^{\frac{2k}{\al d_*^k}} D_2^{\frac{1}{2d_*^k}} I^{\ep}_k.
\end{equation*}
By iteration, 
\begin{equation*}
\begin{split}
I^{\ep}_{k+1}&\leq D_1 d_*^k D_2^{\frac{1}{2d_*^k}} +D_1 d_*^{k-1} d_*^{\frac{2k}{\al d_*^k}} D_2^{\frac{1}{2d_*^k}+\frac{1}{2d_*^{k-1}}} +d_*^{\frac{2k}{\al d_*^k}+\frac{2(k-1)}{\al d_*^{k-1}}} D_2^{\frac{1}{2d_*^k}+\frac{1}{2d_*^{k-1}}}I^{\ep}_{k-1}\\
&\leq \cdots\leq  D_1 D_2^{\frac{1}{2}\sum_{k=0}^{\infty}\frac{1}{d_*^k}} d_*^{\frac{2}{\al}\sum_{k=0}^{\infty} \frac{k}{d_*^k}} \left(\sum_{i=0}^{k} d_*^i +I^{\ep}_0\right)\lesssim \sum_{i=0}^{k} d_*^i +1,
\end{split}
\end{equation*}
where we used the fact $\sum_{k=0}^{\infty}\left(\frac{1}{d^k_*}+\frac{k}{d^k_*}\right)<\infty$ due to $d_*>1$ and $I_0^{\ep}\lesssim 1$ due to \eqref{eqn-v_ep-uniform-estimates}. As a result, $I^{\ep}_{k+1}\lesssim  d_*^{k}$ for each $k\in\N_0$. Let $n\geq 2$. Then, there exists a unique $k=k(n)\in \N_0$ such that $2d_*^{k}\leq n< 2d_*^{k+1}$. Applying H\"older's inequality yields 
$\|v_{\ep}\|_{n}\lesssim \|v_{\ep}\|_{2d_*^{k+1}}\lesssim d_*^k\leq \frac{n}{2}$, and hence, $\|v_{\ep}\|_{n}\leq\frac{C_6n}{2}$. Setting $\ga:=(eC_6)^{-1}$ results in 
\begin{equation*}
\frac{\ga^n}{n!}\int|v_{\ep}|^n\leq \frac{\ga^n}{n!} \left(\frac{C_6n}{2}\right)^n\leq \frac{\ga^n}{n!} \left(\frac{C_6}{2}\right)^n e^n n! =\frac{1}{2^n},\quad\forall n\geq 2,
\end{equation*}
where we used $n^n\leq e^n n!$ in the second inequality. In addition, we apply H\"older's inequality to find from \eqref{eqn-v_ep-uniform-estimates} that $\|v_{\ep}\|_{1}\lesssim D_0$. As a result, 
\begin{equation*}
\begin{split}
\int  e^{\ga |v_{\ep}|}=1+\ga \int |v_{\ep}|+ \sum_{n=2}^{\infty} \frac{\ga^n}{n!}\int  |v_{\ep}|^n  \lesssim 1+\ga D_0+\frac{1}{2}\approx 1.
\end{split}
\end{equation*}
In particular, 
\begin{equation}\label{eqn-lower-upper-connection}
\int  u_{\ep}^{-\ga}  \int  u_{\ep}^{\ga}=\int  e^{\ga v_{\ep}} \int  e^{-\ga v_{\ep}}\lesssim 1.
\end{equation}

Taking $C_6$ so large that $\ga\in(0,1)$, we deduce $\int  u_{\ep}^{\ga} =\int\frac{u_{\ep}}{u_{\ep}^{1-\ga}} \geq  \frac{1}{\|u_{\ep}\|_{\infty}^{1-\ga}}  \int u_{\ep}\gtrsim1$, where we used the fact $\int u_{\ep}=1$ and Lemma \ref{lem-upper-bound} (1) in the last inequality. It follows from \eqref{eqn-lower-upper-connection} that $\int  u_{\ep}^{-\ga} \lesssim1$, proving \eqref{bounded-u_ep-inverse-norm}, and hence, completing the proof.  
\end{proof}

%%%%%%%%%%%%%%%%%%%%%
Now, we prove Theorem \ref{thm-stronger-uniform-bounds-div-free}.

\begin{proof}[Proof of Theorem \ref{thm-stronger-uniform-bounds-div-free}]
Let $A\in \mathcal{S}\mathcal{A}$.  Without loss of generality, we may assume $p$ is such that $A_i^{\ep}\in W^{2,p}$, for $i\in \{1,\dots, m\}$, $A_0^{\ep}\in W^{1,p}$ and $B\in W^{1,p}$. By Theorem \ref{lem-stationary-measure},  $u_{\ep}\in W^{1,p}$ and \eqref{eqn-fpe} is satisfied. Moreover, the classical regularity theory for elliptic equations ensures $u_{\ep}\in W^{2,p}$.  

The proof is divided into four steps.

\medskip

\paragraph{\bf Step 1} We show that for each $k\in\{1,\dots,n\}$ and $g\in W^{1,2}$, 
\begin{equation}\label{eqn-2023-03-09-5-1}
    \begin{split}
        \frac{\ep^2}{2}\int (A^{\ep}_ig)A^{\ep}_iX_k u_{\ep}-\int (Bg) X_k u_{\ep}=\frac{\ep^2}{2}(\RN{1}+\RN{2}),
    \end{split}
\end{equation}
where 
\begin{equation*}
 \RN{1}:=\int A_i^{\ep}g\left(-(\div A^{\ep}_i) X_k u_{\ep} -[X_k,A^{\ep}_i]u_{\ep}\right)+([A^{\ep}_i, X_k]g) A^{\ep}_i u_{\ep}+2(X_k g) A_0^{\ep} u_{\ep}   
\end{equation*}
% \begin{equation*}
%  \RN{1}:=\int A_i^{\ep}g\left(-(\div A^{\ep}_i) X_k u_{\ep} +2[A^{\ep}_i,X_k]u_{\ep}\right)+2(X_k g) A_0^{\ep} u_{\ep}   
% \end{equation*}
and 
\begin{equation*}
\RN{2}:=\int -(A^{\ep}_i g) (X_k\div A^{\ep}_i) u_{\ep}+([A^{\ep}_i,X_k]g) (\div A^{\ep}_i) u_{\ep}
+2(X_k g)(\div A^{\ep}_0) u_{\ep}.
\end{equation*}
% \begin{equation*}
% \RN{2}:=\int A^{\ep}_i g \left(-X_k(\div A^{\ep}_i)+\div[A^{\ep}_i,X_k]\right)u_{\ep}
% +2(X_k g)(\div A^{\ep}_0) u_{\ep}-([A^{\ep}_i, [A^{\ep}_i,X_k]]g)u_{\ep}.
% \end{equation*}
% \begin{equation}\label{eqn-2023-03-09-5-1}
%     \begin{split}
%         -\frac{\ep^2}{2}\left( -\int (A^{\ep}_ig)A^{\ep}_iX_k u_{\ep}+\RN{1}+\RN{2}+\RN{3}\right)+\ep^2\int (A^{\ep}_0 X_k g) u_{\ep}-\int (Bg) X_k u_{\ep}=0,
%     \end{split}
% \end{equation}
%where $\RN{1}:=\int (X_k A^{\ep}_i g)(\div A^{\ep}_i)u_{\ep}$, $\RN{2}:=-\int (A^{\ep}_i[A^{\ep}_i,X_k] g) u_{\ep}$ and $\RN{3}:=-\int (A^{\ep}_ig) [X_k,A^{\ep}_i]u_{\ep}$.

To do so, we fix $k\in \{ 1,\dots, n\}$ and take $g\in W^{2,2}$. Setting $f=X_k g\in W^{1,2}$ in \eqref{eqn-fpe}, we find 
\begin{equation}\label{eqn-2023-03-09-2}
    -\frac{\ep^2}{2}\int A^{\ep}_iX_k g\left[A^{\ep}_i u_{\ep}+(\div A^{\ep}_i) u_{\ep}\right] +\ep^2\int (A^{\ep}_0 X_k g) u_{\ep}+\int  (BX_k g) u_{\ep}=0.
\end{equation}

Straightforward calculations yield
\begin{equation}\label{eqn-2023-03-09-3}
    \begin{split}
        &\int A^{\ep}_iX_k g\left[A^{\ep}_i u_{\ep}+(\div A^{\ep}_i) u_{\ep}\right]-2\int (A^{\ep}_0 X_k g) u_{\ep}\\
        &\qquad =\int (X_k A^{\ep}_i g+[A^{\ep}_i,X_k]g)\left[A^{\ep}_i u_{\ep}+(\div A^{\ep}_i)u_{\ep}\right]-2\int (A^{\ep}_0 X_k g) u_{\ep}\\
        &\qquad =\int (X_k A^{\ep}_i g)A^{\ep}_i u_{\ep}+\int (X_k A^{\ep}_i g)(\div A^{\ep}_i) u_{\ep}+\int ([A^{\ep}_i,X_k]g) \left[A^{\ep}_i u_{\ep}+(\div A^{\ep}_i) u_{\ep}\right]-2\int (A^{\ep}_0 X_k g) u_{\ep} \\
        %&\qquad =\int (X_k A^{\ep}_i g)A^{\ep}_i u_{\ep}+\int (X_k A^{\ep}_i g)(\div A^{\ep}_i)u_{\ep}-\int (A^{\ep}_i[A^{\ep}_i,X_k] g) u_{\ep} -2\int (A^{\ep}_0 X_k g) u_{\ep}\\
        &\qquad=:\int (X_k A^{\ep}_i g)A^{\ep}_i u_{\ep}+\RN{1}'+\RN{2}'+\RN{3}'.
    \end{split}
\end{equation}

Since $\div X_k=0$ and $[X_k,B]=0$, we derive 
\begin{equation*}
    \begin{split}
        \int (X_k A^{\ep}_i g)A^{\ep}_i u_{\ep}&=-\int (\div X_k)A^{\ep}_igA^{\ep}_i u_{\ep}-\int (A^{\ep}_ig)X_kA^{\ep}_i u_{\ep} \\
        &=-\int (A^{\ep}_ig)A^{\ep}_iX_k u_{\ep}-\int (A^{\ep}_ig) [X_k,A^{\ep}_i]u_{\ep}
    \end{split}
\end{equation*}
% we derive 
% \begin{equation}\label{eqn-2023-03-09-3}
%     \begin{split}
%        &\int A^{\ep}_iX_k g\left[A^{\ep}_i u_{\ep}+(\div A^{\ep}_i) u_{\ep}\right]-2\int (A^{\ep}_0 X_k g) u_{\ep}
%     =-\int (A^{\ep}_ig)A^{\ep}_iX_k u_{\ep}-\int (A^{\ep}_ig) [X_k,A^{\ep}_i]u_{\ep}+\RN{1}'+\RN{2}'+\RN{3}'.
%     \end{split}
% \end{equation}
% where 
% $$
% \RN{3}:=-\int (\div X_k)A_igA_i u_{\ep} -\int (A_ig) [X_k,A_i]u_{\ep} .
% $$
and 
\begin{equation*}%\label{eqn-2023-03-09-1}
\begin{split}
    \int  (BX_k g) u_{\ep}&=\int  (X_kBg+[B,X_k]g) u_{\ep}=\int (X_k Bg) u_{\ep}\\
    &=-\int (Bg) u_{\ep}(\div X_k)-\int  (Bg)X_k u_{\ep}=-\int  (Bg)X_k u_{\ep},
\end{split}
\end{equation*}
which together with \eqref{eqn-2023-03-09-2} and \eqref{eqn-2023-03-09-3} lead to 
\begin{equation}\label{eqn-2023-10-12-1}
\frac{\ep^2}{2}\int (A^{\ep}_ig)A^{\ep}_iX_k u_{\ep}-\int (Bg) X_k u_{\ep}=\frac{\ep^2}{2}(\RN{1}'+\RN{2}'+\RN{3}')- \frac{\ep^2}{2}\int (A^{\ep}_ig) [X_k,A^{\ep}_i]u_{\ep}.
\end{equation}

Note that 
\begin{equation*}
\begin{split}
    \RN{1}'&=-\int (A^{\ep}_i g) (\div A^{\ep}_i) u_{\ep} (\div X_k)-\int (A^{\ep}_ig) X_k(u_{\ep}\div A^{\ep}_i)\\
    &=-\int (A^{\ep}_ig) (\div A^{\ep}_i) X_k u_{\ep}-\int  (A^{\ep}_i g) (X_k\div A^{\ep}_i) u_{\ep},
\end{split}
\end{equation*}
% \begin{equation*}
%     \begin{split}
%     \RN{2}'&=-\int ([A^{\ep}_i, X_k]A^{\ep}_i g) u_{\ep}-\int ([A^{\ep}_i, [A^{\ep}_i,X_k]]g)u_{\ep}\\
%     &=\int (A^{\ep}_i g)u_{\ep} (\div[A^{\ep}_i, X_k]) +\int (A^{\ep}_i g) [A^{\ep}_i,X_k] u_{\ep}-\int ([A^{\ep}_i, [A^{\ep}_i,X_k]]g)u_{\ep},
%     \end{split}
% \end{equation*}
and $\RN{3}'=2\int X_k g\left[(\div A^{\ep}_0) u_{\ep}+A^{\ep}_0 u_{\ep}\right]$. Inserting them into \eqref{eqn-2023-10-12-1} results in \eqref{eqn-2023-03-09-5-1}.
Obviously, each term in \eqref{eqn-2023-03-09-5-1} is well-defined even if $g\in W^{1,2}$. Since $W^{2,2}$ is dense in $W^{1,2}$, it follows from standard approximation arguments that \eqref{eqn-2023-03-09-5-1} holds for any $g\in W^{1,2}$.

\medskip

In {\bf Step 2} and {\bf Step 3}, for any two constants $\al(\ep,q)$ and $\be(\ep,q)$ indexed by $(\ep,q)$ with $q$ to be included in {\bf Step 2}, we write $\al(\ep,q) \lesssim \be(\ep,q)$ to imply the existence of $C$, which is independent of $\ep$ and $q$, such that $\al(\ep,q) \leq C \be(\ep,q)$ for any $(\ep,q)$. 

\medskip

\paragraph{\bf Step 2} We show  
\begin{equation}\label{eqn-2023-03-21-2}
\begin{split}
    \int |X_k u_{\ep}|^{2q}|\nabla X_k u_{\ep}|^2 &\lesssim \sum_{k=1}^{n}\left(\|X_k u_{\ep}\|_{(q+1)r}^{2(q+1)}+\|X_k u_{\ep}\|_{qr}^{2q}\right),\quad \forall q\geq 0,
\end{split}
\end{equation}
where $r:=\left(\frac{1}{2}-\frac{1}{p}\right)^{-1}\in (2, 2d_*)$. We remind the reader of the use of Einstein's summation convention on the index $k$ in the left-hand side of \eqref{eqn-2023-03-21-2}. In the following, we first prove inequalities for fixed $k\in \{1,\dots, n\}$ and then sum them up to achieve \eqref{eqn-2023-03-21-2}. 

Let $q\geq 0$ and $k\in \{1,\dots, n\}$. Set $g:=|X_k u_{\ep}|^{2q}X_k u_{\ep}\in W^{1,2}$ in \eqref{eqn-2023-03-09-5-1}. Note that
\begin{equation}\label{eqn-2023-03-09-4}
Y g=(2q+1)|X_k u_{\ep}|^{2q} Y X_k u_{\ep},    
\end{equation}
for any vector field $Y:M\to TM$. Then, we see from  $\div B=0$ that 
\begin{equation*}
    \begin{split}
        -\int  (Bg)X_k u_{\ep}=-(2q+1)\int  |X_{k} u_{\ep}|^{2q} (B X_k u_{\ep}) X_k u_{\ep} =-\frac{2q+1}{2q+2} \int  B|X_k u_{\ep}|^{2q+2}=0.
    \end{split}
\end{equation*}
Hence, it follows from \eqref{eqn-2023-03-09-5-1} that
\begin{equation}\label{eqn-2023-03-09-5}
    \begin{split}
        \int (A^{\ep}_ig)A^{\ep}_iX_k u_{\ep} =\RN{1}+\RN{2}.
    \end{split}
\end{equation}

Applying \eqref{eqn-2023-03-09-4}, we derive from {\bf (A2)} that 
\begin{equation}\label{eqn-2023-03-09-6}
\begin{split}
    \frac{\text{LHS of }\eqref{eqn-2023-03-09-5}}{2q+1}=\sum_{i=1}^d\int  |X_k u_{\ep}|^{2q} |A_i^{\ep} X_k u_{\ep}|^2 \gtrsim \int  |X_k u_{\ep}|^{2q} |\nabla X_k u_{\ep}|^2 
\end{split}
\end{equation}
and
\begin{equation}\label{eqn-2023-03-09-7}
        \frac{\text{RHS of }\eqref{eqn-2023-03-09-5}}{2q+1}=\frac{\RN{1}}{2q+1}+\frac{\RN{2}}{2q+1},
\end{equation}
with
\begin{equation*}
\begin{split}
    \frac{\RN{1}}{2q+1}&=-\int|X_k u_{\ep}|^{2q}(A^{\ep}_i X_k u_{\ep})\left((\div A^{\ep}_i) X_k u_{\ep} + [X_k, A^{\ep}_i]u_{\ep}\right)\\
    &\quad +\int |X_k u_{\ep}|^{2q}\left(([A^{\ep}_i,X_k]X_k u_{\ep}) A^{\ep}_i u_{\ep}+2(X^2_k u_{\ep}) A^{\ep}_0 u_{\ep}\right)
\end{split}
\end{equation*}
% $$
% \RN{1}':=\frac{\RN{1}}{2q+1}=\int|X_k u_{\ep}|^{2q} A^{\ep}_i X_k u_{\ep} \left(-(\div A^{\ep}_i)X_k u_{\ep}+2[A^{\ep}_i,X_k]u_{\ep}\right)+2|X_k u_{\ep}|^{2q}X_k^2 u_{\ep} A^{\ep}_0 u_{\ep}
% $$
and
\begin{equation*}
\begin{split}
\frac{\RN{2}}{2q+1}=\int |X_k u_{\ep}|^{2q} u_{\ep}\left( - (A^{\ep}_i X_k u_{\ep}) (X_k \div A^{\ep}_i) +([A^{\ep}_i, X_k]X_k u_{\ep})\div A^{\ep}_i +2(X_k^2 u_{\ep}) \div A^{\ep}_0\right).
\end{split}
\end{equation*}

Since $A^{\ep}_i\in W^{2,p}$ for $i\in \{1,\dots,m\}$, and $A^{\ep}_0,\,\, X_k\in W^{1,p}$, we apply the Sobolev embedding theorem to find that
$A^{\ep}_0,\,\, A^{\ep}_i, \,\, \div A^{\ep}_i,\,\, X_k \in C^0$ and $[X_k, A^{\ep}_i], \,\, X_k\div A_i^{\ep},\,\, \div A_0^{\ep}\in L^p$. Hence, an application of H\"older's inequality yields
\begin{equation*}
\begin{split}
    \frac{|\RN{1}|}{2q+1}&\lesssim \left(\int |X_k u_{\ep}|^{2q}|\nabla X_k u_{\ep}|^2 \right)^{\frac{1}{2}} \left\||X_k u_{\ep}|^{q} \nabla u_{\ep}\right\|_r,
\end{split}
\end{equation*}
and
\begin{equation*}
\begin{split}
\frac{|\RN{2}|}{2q+1}&\lesssim \sup_{\ep}\|u_{\ep}\|_{\infty}\times \left(\int |X_k u_{\ep}|^{2q}|\nabla X_k u_{\ep}|^2 \right)^{\frac{1}{2}} \left\| |X_k u_{\ep}|^q \right\|_{r} %\left( \|X_k\|_{2,p}+\|A^{\ep}_0\|_{1,p}+\sum_{i=1}^m\|A^{\ep}_i\|_{2,p}\right) \\
    \lesssim \left(\int |X_k u_{\ep}|^{2q}|\nabla X_k u_{\ep}|^2\right)^{\frac{1}{2}} \left\|X_k u_{\ep} \right\|_{qr}^q, 
\end{split}
\end{equation*}
where we used the fact $\|u_{\ep}\|_{\infty}\lesssim1$ due to Theorem \ref{thm-lower-bound-div-free} in the second inequality. This together with \eqref{eqn-2023-03-09-5}-\eqref{eqn-2023-03-09-7} gives
\begin{equation*}
\begin{split}
\int |X_k u_{\ep}|^{2q}|\nabla X_k u_{\ep}|^2& \lesssim \left\||X_k u_{\ep}|^{q} \nabla u_{\ep}\right\|^2_r +\|X_k u_{\ep}\|_{qr}^{2q}.
\end{split}
\end{equation*}
Since $\{X_j\}_{j=1}^{n}$ spans the tangent bundle $TM$, we see from Lemma \ref{lem-vector-basis} the existence of $D>0$ (independent of $\ep$) such that
$\frac{1}{D}|\nabla u_{\ep}|\leq \sum_{j=1}^{n} |X_j u_{\ep}|\leq D |\nabla u_{\ep}|$ on $M$. That is,  
$|\nabla u_{\ep}|\approx \sum_{j=1}^{n} |X_j u_{\ep}|$ uniformly on $M$. Hence, 
\begin{equation*}
 \int |X_k u_{\ep}|^{2q}|\nabla X_k u_{\ep}|^2\lesssim \sum_{j=1}^{n} \|X_j u_{\ep}\|_{(q+1)r}^{2(q+1)}+\|X_k u_{\ep}\|_{qr}^{2q}.   
\end{equation*}
Summarizing the above inequalities in $k$ yields \eqref{eqn-2023-03-21-2}.

\medskip

\paragraph{\bf Step 3} We prove the existence of $C_*>0$ independent of $\ep$ such that
\begin{equation}\label{eqn-Oct-29}
    \sum_{k=1}^{n}\|X_k u_{\ep}\|_{2(q+1)d_*}\leq C_* (q+1)^{\frac{1}{\al (q+1)}}\sum_{k=1}^{n}\|X_k u_{\ep}\|_{2(q+1)}+C_*(q+1)^{\frac{1}{\al (q+1)}},\quad \forall q\geq 0,
\end{equation}
where $\al=1-\frac{d}{p}\in (0,1)$. 

According to H\"older's inequality and Young's inequality, we find
\begin{equation*}
    \|X_k u_{\ep}\|_{qr}^{2q}\lesssim  \|X_k u_{\ep}\|_{(q+1)r}^{2q}\lesssim \frac{q}{q+1} \|X_k u_{\ep}\|_{(q+1)r}^{2(q+1)} +\frac{1}{q+1}.
\end{equation*}
Inserting this into \eqref{eqn-2023-03-21-2} gives rise to 
\begin{equation}\label{eqn-2023-03-26-1}
    \int |X_k u_{\ep}|^{2q}|\nabla X_k u_{\ep}|^2  \lesssim \sum_{k=1}^{n}\|X_k u_{\ep}\|_{(q+1)r}^{2(q+1)}+1.
\end{equation}

Noting that
$$
(q+1)^2\int |X_k u_{\ep}|^{2q}|\nabla X_k u_{\ep}|^2 =\left\|\nabla |X_k u_{\ep}|^{q+1} \right\|_2^2,
$$
and 
$$
\sum_{k=1}^{n}\|X_k u_{\ep}\|_{2(q+1)d_*}^{2(q+1)}\lesssim \sum_{k=1}^{n}\left(\|X_k u_{\ep}\|_{2(q+1)}^{2(q+1)}+ \left\|\nabla |X_k u_{\ep}|^{q+1}\right\|_2^2 \right)
$$
due to the Sobolev embedding theorem, we derive from \eqref{eqn-2023-03-26-1} that
\begin{equation}\label{eqn-2023-03-26-2}
\begin{split}
    \sum_{k=1}^{n}\|X_k u_{\ep}\|_{2(q+1)d_*}^{2(q+1)}
    &\leq C_1 (q+1)^2\sum_{k=1}^{n} \left(\|X_k u_{\ep}\|_{2(q+1)}^{2(q+1)}+\|X_k u_{\ep}\|_{(q+1)r}^{2(q+1)}\right)+C_1 (q+1)^2.
\end{split}
\end{equation}

Since applications of the interpolation inequality and then Young's inequality  lead to 
\begin{equation*}
    \begin{split}
        \|X_k u_{\ep}\|_{(q+1)r}^{2(q+1)}
        &\leq \|X_k u_{\ep}\|_{2(q+1)}^{2\al(q+1)}\|X_k u_{\ep}\|_{2(q+1)d_*}^{2(1-\al)(q+1)}\\
        &\leq \al \de^{-\frac{1}{\al}} \|X_k u_{\ep}\|_{2(q+1)}^{2(q+1)}+(1-\al)\de^{\frac{1}{1-\al}} \|X_k u_{\ep}\|_{2(q+1)d_*}^{2(q+1)}, 
    \end{split}
\end{equation*}
where $\al:=1-\frac{p}{d}\in (0,1)$, we set $\de:=[2(1-\al)C_1(q+1)^2]^{-(1-\al)}$ to derive from  \eqref{eqn-2023-03-26-2} that 
\begin{equation*}
\begin{split}
     \sum_{k=1}^{n}\|X_k u_{\ep}\|_{2(q+1)d_*}^{2(q+1)}&\leq C_2(q+1)^{\frac{2}{\al}} \sum_{k=1}^{n}\|X_k u_{\ep}\|_{2(q+1)}^{2(q+1)}+C_2(q+1)^{\frac{2}{\al}}.
\end{split}
\end{equation*}
Taking the $2(q+1)$-th root of both sides leads to
\begin{equation*}
\begin{split}
    \sum_{k=1}^{n}\|X_k u_{\ep}\|_{2(q+1)d_*}\lesssim \left( \sum_{k=1}^{n}\|X_k u_{\ep}\|_{2(q+1)d_*}^{2(q+1)}\right)^{\frac{1}{2(q+1)}}\lesssim (q+1)^{\frac{1}{\al(q+1)}}\sum_{k=1}^{n}\|X_k u_{\ep}\|_{2(q+1)}+ (q+1)^{\frac{1}{\al(q+1)}}
\end{split}
\end{equation*}
and thus, \eqref{eqn-Oct-29} holds.
\medskip

\paragraph{\bf Step 4} We finish the proof. Set $q+1=d_*^k$ for $k\in \N_0$ and denote $I^{\ep}_k:=\sum_{k=1}^{n}\|X_k u_{\ep}\|_{2d_*^{k}}$. It follows from the result in {\bf Step 3} that 
$$
I^{\ep}_{k+1}\leq C_* d_*^{\frac{k}{\al d^k_*}}I^{\ep}_k+C_*d_*^{\frac{k}{\al d^k_*}}. 
$$
By iteration, we arrive at
\begin{equation*}
\begin{split}
    I^{\ep}_{k+1}\leq C_*^{2}d_*^{\frac{k}{\al d^k_*}+\frac{k}{\al d^{k-1}_*}}I^{\ep}_{k-1}+C_*^{2}d_*^{\frac{k}{\al d^k_*}+\frac{k}{\al d^{k-1}_*}}\leq \dots\leq C_*^{k+1}d_*^{\frac{1}{\al}\sum_{i=0}^{\infty}\frac{k}{d_*^i}}I^{\ep}_0+kC_*^{k+1}d_*^{\frac{1}{\al}\sum_{i=0}^{\infty}\frac{k}{d_*^i}}.
\end{split}
\end{equation*}

Recall from {\bf Step 2} that  $\sum_{k=1}^{n} |X_k u_{\ep}| \approx |\nabla u_{\ep}|$ uniformly on $M$. Since $\|\nabla u_{\ep}\|_2\lesssim1$ by Lemma \ref{lem-uniform-L2-gradient}, we find $I^{\ep}_0\approx \|\nabla u_{\ep}\|_2\lesssim 1$, and thus, $\|\nabla u_{\ep}\|_{2d^k_*}\approx I^{\ep}_{k}\lesssim 1$ for all $k\in \N$. The interpolation inequality then ensures $\|\nabla u_{\ep}\|_{p'}\lesssim 1$ for any $p'>1$. This completes the proof. 
\end{proof}

%%%%%%%%%%%%%%%%%%%%%%%%%%%%%%

%%%%%%%%%%%%%%%%%%%%%%%%%%%%%%%
\section{\bf Proof of main results}\label{sec-main}

This section is devoted to the proof of our main results. 

\subsection{Converting to divergence-free vector fields}\label{subsec-transform}

Recall that the system \eqref{eqn-ode} is assumed to have an invariant measure $\mu_{0}$ with a positive density $u_{0}\in W^{1,p_0}$ for some $p_0>d$. We introduce the invertible transformations to convert \eqref{sde-manifold} into a system whose unperturbed part is divergence-free so that results obtained in Section \ref{sec-div-free} apply.

Set 
\begin{equation}\label{def-tilde-B}
\Tilde{B}:=u_{0} B, \quad \Tilde{u}_0:= 1 \quad \andd\quad d\Tilde{\mu}_0:=d{\rm Vol},   
\end{equation}
and for each $\ep$,
\begin{equation}\label{def-new-measure-diffusion}
\Tilde{A}^{\ep}_0:=u_{0} A^{\ep}_0-\frac{1}{2}\sqrt{u_{0}}(A^{\ep}_j \sqrt{u_{0}})A^{\ep}_j,\quad \Tilde{A}^{\ep}_i:=\sqrt{u_{0}}A^{\ep}_i\quad \text{for}\quad  i\in \{1,\dots, m\},
\end{equation}

\begin{equation}\label{def-tilde-u}
    \Tilde{u}_{\ep}:=\frac{u_{\ep}}{u_{0}} \quad \andd \quad d\Tilde{\mu}_{\ep}:=\Tilde{u}_{\ep}d{\rm Vol}.
\end{equation}
Denote  $\Tilde{A}:=\left\{\Tilde{A}^{\ep}_i,\,\,i\in \{0,\dots, m\},\,\, \ep\right\}$. The following result is elementary. 

\begin{lem}\label{lem-new-measure-vectors}
Let $A\in \mathcal{A}$. Then, the following hold.
\begin{itemize}
    \item [(1)] $\Tilde{B}\in W^{1,p_0}$ and $\div \Tilde{B}=0$.
    
    %\item [(2)] $\Tilde{\mu}_0$ is the invariant measure of $\Tilde{\vp}^t$.
    
    \item [(2)] $\Tilde{A}\in \mathcal{A}$.

    \item [(3)] $\tilde{\mu}_{\ep}$ is a stationary measure of the SDE \eqref{sde-manifold} with $B$ and $A$ replaced by $\tilde{B}$ and $\Tilde{A}$, respectively. 
    %the Fokker-Planck equation $\pa_t u=\Tilde{\LL}^*_{\ep} u$, where $\Tilde{\LL}_{\ep}:=\frac{\ep^2}{2}\sum_{i=1}^m (\tilde{A}^{\ep}_i)^2 + \ep^2 \tilde{A}^{\ep}_0 +\tilde{B}$.
\end{itemize}
\end{lem}
\begin{proof}
Recall that $u_{0}\in W^{1,p_0}$ and $u_0>0$. It follows from the embedding $W^{1,p_0}\hookrightarrow C^{0}$ that $0<\min u_{0}\leq \max u_{0}<\infty$. 

\medskip

(1) Clearly, $B\in W^{1,\infty}$ implies $\tilde{B}\in W^{1,p_0}$. It remains to prove $\div\Tilde{B}=0$. The fact that $\mu_{0}$ is an invariant measure of \eqref{eqn-ode} ensures $\int fu_{0}=\int (f\circ\vp^{t})u_{0}$ for all $f\in C^{0}$, where $\vp^{t}$ is the flow generated by \eqref{eqn-ode} or the vector field $B$. If $f\in C^{1}$, we differentiate $\int (f\circ\vp^{t})u_{0}$ with respect to $t$ and then set $t=0$ to derive
\begin{equation*}
    0=
    \int (Bf) u_{0}=\int (u_{0}B)f=\int \Tilde{B}f.
\end{equation*}
Applying the divergence theorem then gives $\int f\div \Tilde{B}=0$. Hence, we see from the arbitrariness of $f\in C^1$ that $\div \Tilde{B}=0$. 

% To see (2), we calculate
% \begin{equation}\label{eqn-Aug-30-1}
%     \left. \frac{d}{dt}\right|_{t=0}\int  (g\circ \tilde\vp^t) \tilde u_0=\int  (\Tilde{B} g) \Tilde{u}_0=0,\quad \forall g\in C^1.
% \end{equation}
% Fix $g_1\in C^1$. Note that
% \begin{equation*}
% \begin{split}
%     \left. \frac{d}{dt}\right|_{t=t_0}\int  (g_1\circ \tilde\vp^t) \tilde u_0=\left. \frac{d}{dt}\right|_{t=t_0}\int  (g_1\circ \tilde\vp^{t_0}\circ \tilde\vp^{t-t_0})\tilde u_0=\left. \frac{d}{dt}\right|_{t=0}\int  (g_1\circ \tilde\vp^{t_0}\circ \tilde\vp^{t})\tilde u_0,\quad \forall t_0\geq 0.
% \end{split}
% \end{equation*}
% Applying \eqref{eqn-Aug-30-1} to $g=g_1\circ \tilde\vp^{t_0}$ leads to 
% \begin{equation*}
%     \left. \frac{d}{dt}\right|_{t=t_0}\int  (g_1\circ \tilde\vp^t) \tilde u_0=0,\quad \forall t_0\geq 0.
% \end{equation*}
% Since $g_1$ is arbitrary in $C^1$, it follows that $\Tilde{\mu}_0$ is invariant under $\vp^t$.

\medskip

(2) Since $A\in \mathcal{A}$, there exists $p>d$ such that $\|A^{\ep}_0\|_p+\max_i\|A^{\ep}_i\|_{1,p}\lesssim1$. Since $u_0\in W^{1,p_0}$, straightforward calculations show that for $p_1:=\min \{p_0,p\}$, there holds $\|\Tilde{A}^{\ep}_0\|_{p_1}+\max_i\|\Tilde{A}^{\ep}_i\|_{1, p_1}\lesssim1$. That is, $\Tilde{A}$ satisfies {\bf (A1)}. Obviously, $\Tilde{A}$ satisfies {\bf (A2)}, and hence, $\Tilde{A}\in \mathcal{A}$. 

% Conversely, suppose $\Tilde{A}\in \mathcal{A}$. Note that 
% \begin{equation*}
%     A^{\ep}_0=\frac{1}{u_0}\Tilde{A}^{\ep}_0 +\frac{1}{4 u_0^2} (\Tilde{A}^{\ep}_j u_0)\Tilde{A}^{\ep}_j\quad\andd \quad A^{\ep}_i=\frac{1}{\sqrt{u_0}}\Tilde{A}^{\ep}_i\quad \text{for } i\in \{1,\dots, m\}.
% \end{equation*}
% Following similar calculations yields $A\in \mathcal{A}$. 

\medskip

(3) For $f\in C^2$, we calculate 
\begin{equation*}
    \begin{split}
        \int (\LL_{\ep}f) u_{\ep}&=\int  \left[\frac{\ep^2}{2}u_{0}A^{\ep}_iA^{\ep}_i f+\ep^2 u_{0}A^{\ep}_0f+u_{0}Bf \right]\Tilde{u}_{\ep} \\
        &=\int \left[\frac{\ep^2}{2}\sqrt{u_{0}}A^{\ep}_i(\sqrt{u_{0}}A^{\ep}_if)-\frac{\ep^2}{2}\sqrt{u_{0}}(A^{\ep}_i \sqrt{u_{0}})A^{\ep}_i f+\ep^2 u_{0}A^{\ep}_0f+u_{0}Bf\right] \Tilde{u}_{\ep}=\int (\tilde{\LL}_{\ep}f)\tilde{u}_{\ep}.
    \end{split}
\end{equation*}
Since $\int (\LL_{\ep}f) u_{\ep}=0$, we conclude that $\Tilde{\mu}_{\ep}$ is a stationary measure of the SDE \eqref{sde-manifold} with $B$ and $A$ replaced by $\tilde{B}$ and $\Tilde{A}$, respectively. 
\end{proof}

%%%%%%%%%%%%%%%%%%

\subsection{Proof of Theorems \ref{thm-uniform-estimate-stability}-\ref{thm-chi-convergence}}

For the proof of Theorem \ref{thm-uniform-estimate-stability}, we need the following lemma. Recall the definition of a physical measure and its basin from Subsection \ref{subsec-statement-main-results}.

\begin{lem}\label{lem-properties-physical-measure}
If $\mu$ is a physical measure of \eqref{eqn-ode}, then $\mu=\frac{\mu_0|_{B_{\mu}}}{\mu_{0}(B_{\mu})}$ and it is ergodic.
\end{lem}
\begin{proof}
Fix $f\in C^0$. The fact that $\mu$ is a physical measure implies
\begin{equation}\label{eqn-consequence-physical-measure}
    \lim_{t\to\infty}\frac{1}{t}\int_{0}^{t}f(\vp^{s}(x))ds=\int fd\mu,\quad\forall x\in B_{\mu}.
\end{equation}
It follows from Fubini's theorem and the dominated convergence theorem  that
$$
\lim_{t\to\infty}\frac{1}{t}\int_{0}^{t}\int_{B_{\mu}}f\circ\vp^{s}d\mu_0ds=\lim_{t\to\infty}\int_{B_{\mu}}\left(\frac{1}{t}\int_{0}^{t}f\circ\vp^{s}ds\right)d\mu_0=\mu_{0}(B_{\mu})\int fd\mu.
$$
Note that
$$
\int_{B_{\mu}}f\circ\vp^{s}d\mu_0=\int (1_{B_{\mu}}f)\circ\vp^{s}d\mu_{0}=\int 1_{B_{\mu}}fd\mu_0,
$$
where we used the $\vp^{t}$-invariance of $B_{\mu}$ in the first equality and the fact that $\mu_0$ is an invariant measure of $\vp^{t}$ in the second one. Hence, we arrive at $\int 1_{B_{\mu}}fd\mu_0=\mu_{0}(B_{\mu})\int fd\mu$ for all $f\in C^0$,
leading to $\mu=\frac{\mu_0|_{B_{\mu}}}{\mu_{0}(B_{\mu})}$. In particular, $\mu(B_{\mu})=1$.

For the ergodicity of $\mu$, we note that $\mu(B_{\mu})=1$ implies that $\lim_{t\to\infty}\frac{1}{t}\int_{0}^{t}\de_{\vp^{s}(x)}ds=\mu$ under the weak*-topology for $\mu$-a.e. $x\in M$. This is equivalent to the ergodicity of $\mu$. 
\end{proof}

\begin{cor}\label{cor-2024-03-06}
The following statements are equivalent.
\begin{itemize}
    \item[(1)] $\mu_0$ is physical.

    \item[(2)] $\mu_0$ is ergodic.

\item[(3)] There is a physical measure $\mu$ of \eqref{eqn-ode} with ${\rm Vol}(B_{\mu})=1$.

\end{itemize}
Whenever these statements hold, $\mu_0$ is the unique physical measure.
\end{cor}
\begin{proof}
(1)$\implies$(2). It follows directly from Lemma \ref{lem-properties-physical-measure}.

(2)$\implies$(1) and (3). Since $\mu_0$ has a positive density, Birkhoff's ergodic theorem ensures that $\mu_0$ itself is physical and satisfies ${\rm Vol}(B_{\mu_0})=1$.

(3)$\implies$(2). If there is a physical measure $\mu$ of \eqref{eqn-ode} with ${\rm Vol}(B_{\mu})=1$, then the equivalence between $\mu_0$ and ${\rm Vol}$ yields $\mu_0(B_{\mu})=\mu_0(M)=1$, and hence, $\mu=\mu_0$ thanks to Lemma \ref{lem-properties-physical-measure}.

Now, suppose that (1)-(3) hold. If $\mu$ is a physical measure, then Lemma \ref{lem-properties-physical-measure} yields $\mu\ll\mu_0$, and hence, $\mu=\mu_0$ thanks to the ergodicity of $\mu_0$. 
\end{proof}

\begin{proof}[Proof of Theorem \ref{thm-uniform-estimate-stability}]
(1) Let $A\in \mathcal{A}$ and $\Tilde{A}$, $\Tilde{B}$ and  $\Tilde{\mu}_{\ep}$ be as defined in \eqref{def-tilde-B}-\eqref{def-tilde-u}. Given Lemma \ref{lem-new-measure-vectors}, we apply Theorem \ref{thm-lower-bound-div-free} to find $\|\Tilde{u}_{\ep}\|_{1,2}\lesssim1$ and $1\lesssim\min \Tilde{u}_{\ep}\leq  \max \Tilde{u}_{\ep}\lesssim1$. Since $u_{\ep}=u_0 \Tilde{u}_{\ep}$ and $u_0\in W^{1,p_0}$ for some $p_0>d$, straightforward calculations yield
$\|u_{\ep}\|_{1,2}\lesssim1$ and $1\lesssim\min u_{\ep}\leq \max u_{\ep}\lesssim1$. As a result, $\{u_{\ep}\}_{\ep}$ is precompact in $W^{1,2}$ under the weak topology. The ``In particular" part follows readily.

(2) The conclusions in (1) guarantees that each element of $\MM_{A}$ (must be an invariant measure of \eqref{eqn-ode}) is equivalent to $\mu_0$. Then, the ergodicity of $\mu_0$ (by Corollary \ref{cor-2024-03-06}) asserts $\MM_{A}=\{\mu_0\}$. 

(3) In this case, $\MM_{A}=\{\mu_0\}$ follows immediately from conclusions in (1). 

Whenever either (2) or (3) is true, there holds $\MM_{A}=\{\mu_0\}$. It then follows from the uniform estimates of $\{u_{\ep}\}_{\ep}$ in $W^{1,2}$ (by (1)) and the Rellich–Kondrachov theorem that $\lim_{\ep\to0}u_{\ep}=u_{0}$ weakly in $W^{1,2}$ and strongly in $L^{p}$ for any $p\in[1,\frac{2d}{d-2})$. Since for each $p_2>p_1\geq1$ there holds
$$
\|u_{\ep}-u\|_{p_{2}}\leq \left[\max \left( u_{\ep}+u\right)\right]^{1-\frac{p_{1}}{p_{2}}}\|u_{\ep}-u\|_{p_1}^{\frac{p_{1}}{p_{2}}}\lesssim \|u_{\ep}-u\|_{p_{1}}^{\frac{p_{1}}{p_{2}}},
$$
we conclude that $\lim_{\ep\to0}u_{\ep}=u_{0}$ holds strongly in $L^{p}$ for any $p\in[1,\infty)$. This completes the proof. 
\end{proof}

Theorem \ref{thm-stronger-uniform-bounds} follows readily.

\begin{proof}[Proof of Theorem \ref{thm-stronger-uniform-bounds}]
(1) Let $A\in \mathcal{SA}$ and $\Tilde{A}$, $\Tilde{B}$ and  $\Tilde{\mu}_{\ep}$ be as defined in \eqref{def-tilde-B}-\eqref{def-tilde-u}. Then, all the results in Lemma \ref{lem-new-measure-vectors} hold. Moreover, since $u_0\in W^{2,p_0}$, we can follow the proof of Lemma \ref{lem-new-measure-vectors} to find $\Tilde{B}\in W^{1,\infty}$  and $\Tilde{A}\in \mathcal{SA}$. 

Set $\Tilde{X}_i:=u_0 X_i$. Clearly, $[\Tilde{X}_i,\Tilde{B}]=0$, $\div \Tilde{X}_i=0$, and $\{\Tilde{X}_i\}_{i=1}^n$ spans the tangent bundle $TM$. We then apply Theorem \ref{thm-stronger-uniform-bounds-div-free} to conclude that $\|\Tilde{u}_{\ep}\|_{W^{1,q}}\lesssim1$ for all $q\geq 1$. Since $u_{\ep}=u_0 \Tilde{u}_{\ep}$ and $u_0\in W^{2,p_0}$, the conclusion follows.

(2) Theorem \ref{thm-uniform-estimate-stability}(3) asserts that $\MM_A=\{\mu_0\}$. The conclusion then follows from (1) and the compact Sobolev embedding theorem. 
\end{proof}

Now, we prove Theorem \ref{thm-construction}.

%-----------------------------------------%
%-----------------------------------------%
%-----------------------------------------%

\begin{proof}[Proof of Theorem \ref{thm-construction}]
    The construction of the manifold $M$ and the vector field $B$ satisfying (1)-(3) are done in {\bf Step 1} and {\bf Step 2}. We prove (i)-(iii) in {\bf Step 3}.  

\medskip

\paragraph{\bf Step 1} Set $D:=\left\{x=(x_1,x_2)\in\R^{2}: x_1^2+x_2^2\leq 1\right\}$. It is well-known \cite{MR0554383} that there exists a smooth area-preserving Bernoulli diffeomorphism $F:D\to D$  such that $F-\text{Id}$ is ``infinitely flat" on $\pa D$, that is, both $F-\text{Id}$ and all its derivatives of any order vanish on $\pa D$.

We construct a smooth suspension flow  $\Tilde{\psi}^t$ of $F$ on the smooth, connected Riemannian suspension manifold $(\Tilde{D},\Tilde{g})$ such that $\Tilde{\psi}^t$ is volume-preserving and has a Poincar\'e map smoothly conjugate to $F$. Moreover, both $\Tilde{D}^{\mathrm{o}}$ and $\pa\Tilde{D}$ are invariant under $\Tilde{\psi}^t$.   

To do so, we define a $\Z$-action on $D\times\R$ generated by the map 
\begin{equation}\label{eqn-equi-map}
(x,s)\mapsto (F(x),s-1): D\times \R\to D\times \R,
\end{equation}
which naturally induces an equivalence relation $\sim$ on $D\times \R$:
$$
(x,s)\sim (x',s')\quad \text{iff}\quad (x',s')=(F^n(x),s-n)\text{ for some }n\in\Z.
$$

Set $\Tilde{D}:=D\times \R/\sim$. It is a smooth, connected, and compact manifold with boundary $\pi_D(\pa D\times \R)$, where $\pi_D:D\times\R\to \Tilde{D}$ is the natural projection. To furnish $\Tilde{D}$ with a Riemannian metric $\Tilde{g}$, we first take the following Riemannian metric on $D\times [-1,0]$:
\begin{equation*}
    g=(dx_1,dx_2) \left[[DF(\ga(x,s))]^s\right]^{\top}[DF(\ga(x,s))]^s(dx_1,dx_2)^{\top}+ds^2=:g_s+ds^2,
\end{equation*}
where $\ga:D\times [-1,0]\to D$ is a smooth function satisfying
\begin{equation*}
    \ga(x,s)=\begin{cases}
        F^{-1}(x),& s\in \left[-1,-\frac{3}{4}\right],\\
        x,&s\in \left[-\frac{1}{4},0\right],
    \end{cases}
    \qquad x\in D.
\end{equation*}
It is easy to check that $g_0=dx_1^2+dx_2^2=F^*g_{-1}$, where $F^*g_{-1}$ denotes the pullback of $g_{-1}$ under $F$. Therefore, through the smooth diffeomorphism $(x,s)\mapsto (F^n(x), s-n)$ on $D\times (n-1,n]$ for each $n\in \Z$, we can extend $g$ to be a smooth Riemannian metric on $D\times \R$, stilled denoted by $g$. Since $F=\mathrm{Id}$ on $\pa D$ and $F-\mathrm{Id}$ is infinitely smooth, we see that $DF=\mathrm{Id}$ on $\pa D$ and $\ga(x,s)=x$ for any $(x,s)\in \pa D\times [-1,0]$, yielding
\begin{equation}\label{g-boundary}
\begin{split}
g&=g_0+ds^2\quad \text{ on }\quad  \pa D\times \R.
\end{split}
\end{equation}
Note that the construction $g$ through the extension ensures its invariance under the $\Z$-action. Hence, a Riemannian metric $\tilde{g}$ on $\Tilde{D}$ is naturally induced. 

Now, we construct the suspension flow $\Tilde{\psi}^t$. Let $\psi^t$ be the flow generated by $\frac{\pa}{\pa s}$ on $D\times \R$. Since $F$ is area-preserving, there holds $|\det (DF)|\equiv 1$, and thus, $\det (g_s)\equiv 1$. Straightforward calculations then yield $\div \frac{\pa}{\pa s}=0$, and hence, $\psi^t$ is volume-preserving. Here, we used $g_s$ with a slight abuse of notation to represent the matrix when defining the $2$-form $g_s$. Since $\psi^t$ preserves the equivalence relation $\sim$, it naturally induces a volume-preserving flow $\Tilde{\psi}^t:=\pi_D \circ \psi^t\circ(\pi_D)^{-1}$ on $\Tilde{D}$ (see Figure \ref{figure-suspension-flow-disk} for an illustration). Clearly, $\Tilde{\psi}^t$ is smooth, generated by the pushforward $(\pi_D)_*\frac{\pa}{\pa s}$, and satisfies 
\begin{equation}\label{eqn-flow-boundary}
\Tilde{\psi^t}(\Tilde{D}^\mathrm{o})=\Tilde{D}^\mathrm{o},\quad \andd\quad \Tilde{\psi^t}(\pi_{D}(x,s))=\pi_{D}(x,s+t),\quad \forall (x,s)\in \pa D\times \R.
\end{equation}
In particular, 
\begin{equation*}
\Tilde{\psi}^1(\pi_D(x,0))=\pi_D(\psi^1(x,0))=\pi_D(x,1)=\pi_D(F(x),0),\quad \forall x\in D.
\end{equation*}
Since $\pi_D|_{D\times \{0\}}:D\times \{0\}\to\pi_D(D\times \{0\})$ is a smooth diffeomorphism, the above equality asserts that $\Tilde{\psi}^1|_{\pi_D(D\times\{0\})}$, which is a Poincar\'e map of $\Tilde{\psi}^t$, is smoothly conjugate to $F$. Hence, $\Tilde{\psi}^t$ is strongly mixing. 

\begin{figure}
	\begin{center}
		\begin{tikzpicture}[
			thick,
			nonterminal/.style = {
				%the shape
				rectangle,			
				% the size
				minimum size = 10mm,			
				% the border
				thick,
				% draw=red!50!black!50,  %50% red, 50% black, and that mixed with 50% white			
				draw = black,
				%the filling
				% top color = white,
				% bottom color = red!50!black!20,
				%font
				% font = \itshape
			}
		]
			\matrix [row sep=10mm, column sep=5mm] {
				%first row
				\node (DR1)  {$(D\times \R,g)$};&&
				\node (DR-quo1) {$(\Tilde{D},\Tilde{g})$};&& \\			
                    %second row
                    \node (DR2) {$(D\times \R,g)$};&&
				\node (DR-quo2) {$(\Tilde{D},\Tilde{g})$};&& \\						
			};
			
			% \graph{
			% 	horizontal arrows: (DR1) -> (DR-quo1);
			%				       (DR2) -> (DR-quo2);
			%     vertical arrows: (DR1)->(DR2), (DR-quo1)->(DR-quo2)
			% };
			\path (DR1) edge[->]node[above]{$\pi_D$}  (DR-quo1);
                \path (DR2) edge[->]node[above]{$\pi_D$}  (DR-quo2);
                
                \path (DR1) edge[->]node[right]{$\psi^t$} (DR2);
                \path (DR-quo1) edge[->] node[right]{$\Tilde{\psi}^t$} (DR-quo2);
		\end{tikzpicture}
	\end{center}
	    \caption{Construction of the suspension flow $\Tilde{\psi}^t$.}\label{figure-suspension-flow-disk}
\end{figure}

\medskip
\paragraph{\bf Step 2} We construct $M$ and $B$. 

Let $\S^2:=\{(x_1,x_2,x_3):x_1^2+x_2^2+x_3^2=1 \}$ be the unit sphere and denote by
$$
\S^2_1:=\{(x_1,x_2,x_3)\in \S^2: x_3\geq 0\} \quad\andd \quad \S^2_2:=\{(x_1,x_2,x_3)\in \S^2: x_3\leq 0\}
$$
the upper hemisphere and lower hemisphere, respectively. Let $h_1:\S^2_1\to D$ and $h_2:\S^2_2\to D$ be the standard stereographic projections from the south pole $(0,0,-1)$ and the north pole $(0,0,1)$, respectively.
% Define the map
% $$
% h: \S^2\to D,\quad (x_1,x_2,x_3)\mapsto (x_1,x_2).
% $$
% {\rd Clearly, for $i=1,2$, the map $h_i:=h|_{\S^2_i}:\S^2_i\to D$ is a  smooth diffeomorphism}. 
Through $h_i$, the equivalence relation $\sim$ on $D\times \R$ is naturally lifted to its counterpart $\sim_i$ on $\S^2_{i}\times \R$. Set $M_i:=\S^2_i\times \R/\sim_i$ and let $\pi_i:\S^2_i\times \R\to M_{i}$ be the natural projection. 
%$F_i:=h_i\circ F\circ h^{-1}_i:\S^2_i\to \S^2_i$ is a $C^{\infty}$ Bernorlli diffeomorphism. As in {\bf Step 1}, $F_i$ defines a $\Z$-action on $\S^2_i\times \R$ and further gives an equivalent relation $\sim_i$. Set $M_i:=\S^2_i\times \R/\sim_i$ and let $\pi_i$ be the natural projection. T
Then, 
$$
\Tilde{h}_i:=\pi_D\circ (h_i,\mathrm{Id})\circ\pi_i^{-1}:M_i\to\Tilde{D}
$$ 
defines a smooth diffeomorphism, allowing us to endow $M_i$ with the Riemannian metric $g_i:=(\Tilde{h}_i)^* \Tilde{g}$ and define a flow $\vp^t_i:=(\Tilde{h}_i)^{-1}\circ\Tilde{\psi}^t\circ \Tilde{h}_i$ on $M_i$. Obviously, $(M_i,g_i)$ is a smooth and connected Riemannian manifold with boundary $\pa:=(\Tilde{h}_i)^{-1}\pa \Tilde{D}$ and $\vp^t_i$ is a smooth and volume-preserving flow, generated by the vector field $B_i:=(\Tilde{h}_i^{-1})_*(\pi_D)_* \frac{\pa}{\pa t}$. Since $\vp^t_i$ is smoothly conjugate to $\Tilde{\psi}^t$, it follows that $\vp^t_i$ is as well strongly mixing. We refer the reader to Figure \ref{figure-suspension-flow} for clarity.  

Note that the map \eqref{eqn-equi-map}, which generates the equivalence relation $\sim$, preserves $\pa D\times \R$ since $F-\mathrm{Id}$ is infinitely smooth on $\pa D$. We see from the form \eqref{g-boundary} of $g$ on $\pa D\times \R$ and the construction of $(M_i,g_i)$, $i=1,2$ that $M:=M_1\cup M_2$ is a smooth, connected, and closed manifold and $g_M:=g_i$ on $M_i$ for $i=1,2$ is a well-defined smooth Riemannian metric on $M$.
%how $M_i$ for $i=1,2$ are obtained that $\pa M_1=\pa M_2=\pa$ and $M:=M_1\cup M_2$ is a smooth, connected, and closed manifold. 
%Since the metric $g$ on $D\times \R$ is defined through $\ga$, or more essentially through $F$ and satisfies \eqref{g-boundary} on $\pa D\times \R$, we find $g_1=g_2$ on $\pa$, and $g_M:=g_i$ on $M_i$ for $i=1,2$ is a well-defined smooth Riemannian metric on $M$. 

Now, we construct a smooth flow $\vp^t$ on $(M,g)$. By the definition of the flow $\psi^t$ on $D\times \R$ that results in a particular form of $\Tilde{\psi}^t$ on $\pa \Tilde{D}$, namely, the second equality in \eqref{eqn-flow-boundary}, we see from the construction of $\vp^t_i$, $i=1,2$ that $\vp^t:=\vp^t_i$ on $M_i$ for $i=1,2$ is a well-defined, smooth and volume-preserving flow on $M$. Setting $B:=B_i$ on $M_i$ for $i=1,2$, we conclude that $B$ is smooth and divergence-free, and generates $\vp^t$. The properties (1)-(3) follow readily.  

\begin{figure}
	\begin{center}
		\begin{tikzpicture}[
			thick,
			nonterminal/.style = {
				%the shape
				rectangle,			
				% the size
				minimum size = 10mm,			
				% the border
				thick,
				% draw=red!50!black!50,  %50% red, 50% black, and that mixed with 50% white			
				draw = black,
				%the filling
				% top color = white,
				% bottom color = red!50!black!20,
				%font
				% font = \itshape
			}
		]
			\matrix [row sep=10mm, column sep=5mm] {
				%first row
				\node (DR)  {$D\times \R$};&&
				\node (DR-quo1) {$\Tilde{D}$};&& 
                    \node (DR-quo2) {$\Tilde{D}$};&& \\			
                    %second row
                    \node (SR) {$\S^2_i\times \R$};&&
				\node (SR-quo1) {$M_i$};&& 
                    \node (SR-quo2) {$M_i$};&& \\			
			};
			
			% \graph{
			% 	horizontal arrows: (DR) -> (DR-quo1)->(DR-quo2);
			%				       (SR) -> (SR-quo1)->(SR-quo2);
			%     vertical arrows: (DR)->(SR), (DR-quo1)->(SR-quo1), (DR-quo2)->(SR-quo2)
			% };
			\path (DR) edge[->]node[above]{$\pi_D$}  (DR-quo1);
                \path (DR-quo1) edge[->] node[above]{$\Tilde{\psi}^t$} (DR-quo2);
                \path (SR) edge[->]node[above]{$\pi_i$}  (SR-quo1);
                \path (SR-quo1) edge[->] node[above]{$\vp^t_i$} (SR-quo2);
                
                \path (SR) edge[->]node[right]{$(h_i, \mathrm{Id})$} (DR);
                \path (SR-quo1) edge[->] node[right]{$\Tilde{h}_i$} (DR-quo1);
                \path (SR-quo2) edge[->] node[right]{$\Tilde{h}_i$} (DR-quo2);
		\end{tikzpicture}
	\end{center}
	    \caption{Construction of the  suspension flow $\vp^t_i$ on $M$.}\label{figure-suspension-flow}
\end{figure}

\medskip

\paragraph{\bf Step 3} (i) is an immediate consequence of (1) and (2). (ii) follows readily from (3). (iii) follows directly from (ii) and Theorem \ref{thm-uniform-estimate-stability}(3).
\end{proof}

%--------------------------------%

Finally, we prove Theorem \ref{thm-chi-convergence} that builds on two lemmas. In the first lemma, we prove a uniform-in-$\ep$ Poincar\'e inequality.

\begin{lem}\label{lem-uniform-Poincare-inequality}
Let $A\in \mathcal{A}$. Then, there exists $C_{*}>0$ such that  
\begin{equation*}
    \int |f-\bar{f}_{\ep}|^2 d\mu_{\ep}\leq C_{*} \sum_{i=1}^m\int |A^{\ep}_if|^2 d\mu_{\ep},\quad \forall f\in W^{1,2}\andd 0<\ep\ll 1,
\end{equation*}
where $\bar{f}_{\ep}:=\int fd\mu_{\ep}$. 
\end{lem}

\begin{proof}
Let $f\in W^{1,2}$ and set $\bar{f}:=\int f$. It follows from the fact $\int |f-\bar{f}_{\ep}|^2 d\mu_{\ep}=\min_{t\in \R}\int |f-t|^2 d\mu_{\ep}$ and Theorem \ref{thm-uniform-estimate-stability} that
\begin{equation}\label{eqn-Jan-1-2}
    \int |f-\bar{f}_{\ep}|^2 d\mu_{\ep}\leq \int |f-\bar{f}|^2 d\mu_{\ep}\leq \left(\sup_{\ep}\max u_{\ep}\right)\|f-\bar{f}\|_2^2\leq C_{1}\left(\sup_{\ep}\max u_{\ep}\right)\|\nabla f\|_2^{2}.
\end{equation}
where we used the classical Poincar\'e inequality in the third inequality and the constant $C_1>0$ is independent of $f$. 

Since $A\in \mathcal{A}$, Theorem \ref{thm-uniform-estimate-stability} implies that
\begin{equation*}
\|\nabla f\|_2^2\leq \frac{1}{\la} \sum_{i=1}^m\int |A^{\ep}_if|^2\leq \frac{1}{\la \inf_{\ep}\min u_{\ep}} \sum_{i=1}^m\int |A^{\ep}_if|^2 d\mu_{\ep},
\end{equation*}
where $\la$ is given in {\bf (A2)} in Definition \ref{admissible-class}. This together with \eqref{eqn-Jan-1-2} yields the result. 
\end{proof}

\begin{rem}\label{rem-constant-C-in-convergence-rate}
It is clear from the proof of Lemma \ref{lem-uniform-Poincare-inequality} that the constant $C_{*}$ can be chosen as $C_{*}=\frac{C_1\sup_{\ep}\max u_{\ep}}{\la \inf_{\ep}\min u_{\ep}}$.    
\end{rem}

The second lemma addresses the well-posedness of the Cauchy problem \eqref{eqn-fpe-time}-\eqref{initial-condition}.

\begin{lem}\label{lem-Cauchy-integral-formula}
Let $A\in \mathcal{A}$ satisfy {\bf (A1)} with $p>d+2$. Assume $\nu_0$ has a density $v_0\in L^2$. Then, the following hold for each $\ep$. 
\begin{itemize}
    \item[(1)] \eqref{eqn-fpe-time}-\eqref{initial-condition} admits a unique solution $v^{\ep}\in V_{2,loc}^{1,0}(M\times [0,\infty))\cap \HH^{1,p}_{loc}(M\times (0,\infty))$, which is continuous and positive in $M\times (0,\infty)$ and satisfies $\int v^{\ep}(\cdot, t)=1$ for all $t\in(0,\infty)$. 
    
    \item[(2)] For any $h\in C^2(\R^2)$ and $f\in W^{1,2}\cap L^{\infty}$, there holds
\begin{equation}\label{eqn-Jan-21-1}
\frac{d}{dt}\int h(v^{\ep},f)=-\frac{\ep^2}{2}\int A^{\ep}_i h_1(v^{\ep},f)\left[A^{\ep}_i v^\ep+(\div A^{\ep}_i) v^{\ep}\right]+\int\left[\ep^{2}A^{\ep}_0h_1(v^{\ep},f)+Bh_1(v^{\ep},f)\right]v^{\ep}
\end{equation}
for all $t\in (0,\infty)$, where $h_1:=\pa_1h$.  
\end{itemize}
 
\end{lem}
\begin{proof}
(1) It follows from Theorem \ref{thm-Cauchy-appx}. Moreover, for each $\phi\in C^{1,1}(M\times [0,\infty))$, there holds 
\begin{equation}\label{eqn-Feb-24-1}
\begin{split}
    &\int \phi(\cdot, t)v^{\ep}(\cdot,t)-\int\phi(\cdot, s)v(\cdot, s)-\int_s^t\int v^{\ep}\pa_t \phi \\
    &\quad = -\frac{\ep^2}{2}\int_s^t\int A^{\ep}_i \phi\left[A^{\ep}_i v^\ep+(\div A^{\ep}_i) v^{\ep}\right]+\int_s^t\int\left(\ep^{2}A^{\ep}_0\phi+B\phi\right)v^{\ep},\quad \forall 0<s<t.
\end{split}
\end{equation}

(2) Fix $h\in C^2(\R^2)$. Assume at the moment that $f\in C^1$ and that $A_i^{\ep}$ for $i\in \{0,\dots, m\}$ and $B$ are  $C^3$. Then, the classical regularity theory for solutions of parabolic equations ensures $v^{\ep}\in C^{2,1}(M\times (0,\infty))$. Set $\phi:=h_1(v^{\ep},f)$ with $h_1:=\pa_1 h$ in \eqref{eqn-Feb-24-1}. Since integration by parts yields 
\begin{equation*}
\begin{split}
    &\int h_1(v^{\ep}(\cdot,t),f)v^{\ep}(\cdot,t)-\int h_1(v^{\ep}(\cdot,s),f)v^{\ep}(\cdot,s)-\int_s^t\int v^{\ep}\pa_t h_1(v^{\ep}, f)\\
    &\quad=\int_s^t\int h_1(v^{\ep}, f)\pa_t v^{\ep}
    =\int_s^t \int \pa_t h(v^{\ep}, f) 
    =\int h(v^{\ep}(\cdot, t), f) -\int h(v^{\ep}(\cdot,s), f),
\end{split}
\end{equation*}
we derive 
\begin{equation}\label{eqn-Feb-24-2}
\begin{split}
    &\int h(v^{\ep}(\cdot, t), f)-\int h(v^{\ep}(\cdot,s), f)\\
    &\quad =-\frac{\ep^2}{2}\int_s^t\int A^{\ep}_i h_1(v^{\ep}, f)\left[A^{\ep}_i v^\ep+(\div A^{\ep}_i) v^{\ep}\right]+\int_s^t\int\left[\ep^{2}A^{\ep}_0h_1(v^{\ep}, f)+Bh_1(v^{\ep}, f)\right]v^{\ep}
\end{split}
\end{equation}
for $0<s<t$, leading to \eqref{eqn-Jan-21-1}.

Now, we address the general case where $A\in \mathcal{A}$ and $B$ are not necessarily smooth and $f\in W^{1,2}\cap L^{\infty}$. Recall that $v^{\ep}\in V_{2,loc}^{1,0}(M\times [0,\infty))\cap C(M\times (0,\infty))$. Hence, all the terms in \eqref{eqn-Feb-24-2} remain well defined. Approximating $A^{\ep}_i$ and $B$ by sequences of smooth functions, we deduce from the uniqueness of the solutions of \eqref{eqn-fpe-time}-\eqref{initial-condition} that solutions of the approximating equations converges to $v^{\ep}$ in $\in V_{2,loc}^{1,0}(M\times [0,\infty))\cap C(M\times (0,\infty))$. Then, standard approximation arguments ensure that \eqref{eqn-Feb-24-2} still holds, from which \eqref{eqn-Jan-21-1} follows readily.
\end{proof}

We are ready to prove Theorem \ref{thm-chi-convergence}.

\begin{proof}[Proof of Theorem \ref{thm-chi-convergence}]
(1) has been proven in Lemma \ref{lem-Cauchy-integral-formula} (1). Assuming the following identity
\begin{equation}\label{an-identity-2025-05-04}
\frac{d}{dt}\chi^2(\nu^{\ep}_t, \mu_{\ep})=-\ep^2 \int \left(A^{\ep}_i \frac{v^{\ep}}{u_{\ep}}\right)^2 u_{\ep}, 
\end{equation}
we prove (2). Set $C:=\frac{1}{C_{*}}$, where $C_{*}$ is given in Lemma \ref{lem-uniform-Poincare-inequality}. Applying Lemma \ref{lem-uniform-Poincare-inequality}, we derive from $\int \frac{v^{\ep}}{u_{\ep}}u_{\ep}=1$ that 
\begin{equation*}
   -\int \left(A^{\ep}_i \frac{v^{\ep}}{u_{\ep}}\right)^2 u_{\ep}\leq -C \int \left|\frac{v^{\ep}}{u_{\ep}}-1\right|^2 d\mu_{\ep}=-C\chi^2(\nu^{\ep}_t, \mu_{\ep}),\quad \forall 0<\ep\ll1.
\end{equation*}
Hence, $\frac{d}{dt}\chi^2(\nu^{\ep}_t, \mu_{\ep})\leq -C\ep^2\chi^2(\nu_t^{\ep}, \mu_{\ep})$, which together with the Gr\"ownwall's inequality leads to $\chi^2(\nu^{\ep}_t, \mu_{\ep})\leq e^{-C\ep^2 t}\chi^2(\nu_0, \mu_{\ep})$.

It remains to establish the identity \eqref{an-identity-2025-05-04}. Choosing $h(r_1,r_2)=r_2\left(\frac{r_1}{r_2}-1\right)^2$ for $(r_1,r_2)\in \R^2$ and $f=u_{\ep}\in W^{1,p}\hookrightarrow L^{\infty}$ in Lemma \ref{lem-Cauchy-integral-formula} (2), we achieve 
\begin{equation}\label{eqn-Jan-1-4}
\begin{split}
    \frac{d}{dt}\chi^2(\nu^{\ep}_t, \mu_{\ep})&=-\ep^2\int A^{\ep}_i \left(\frac{v^{\ep}}{u_{\ep}}-1 \right)\left[A^{\ep}_i v^{\ep}+(\div A^{\ep}_i) v^{\ep}\right] +2\int\left[\ep^{2}A^{\ep}_0\left(\frac{v^{\ep}}{u_{\ep}}-1 \right)+B\left(\frac{v^{\ep}}{u_{\ep}}-1 \right)\right]v^{\ep}\\
    &=-\ep^2 \int A^{\ep}_i \frac{v^{\ep}}{u_{\ep}}\left[A^{\ep}_i v^{\ep}+(\div A^{\ep}_i) v^{\ep}\right]+2\int\left(\ep^{2}A^{\ep}_0\frac{v^{\ep}}{u_{\ep}}+B\frac{v^{\ep}}{u_{\ep}}\right)v^{\ep}\\
    &=-\ep^{2}\int \left(A^{\ep}_i \frac{v^{\ep}}{u_{\ep}}\right)^2 u_{\ep}-\frac{\ep^{2}}{2}\int A^{\ep}_i \left(\frac{v^{\ep}}{u_{\ep}}\right)^2A^{\ep}_i u_{\ep}-\ep^2 \int A^{\ep}_i \frac{v^{\ep}}{u_{\ep}}(\div A^{\ep}_i) v^{\ep}\\
    &\quad+2\int\left(\ep^{2}A^{\ep}_0\frac{v^{\ep}}{u_{\ep}}+B\frac{v^{\ep}}{u_{\ep}}\right)v^{\ep},
\end{split}
\end{equation}
where we used in the last equality the following identity
\begin{equation*}
\begin{split}
    \int A^{\ep}_i \frac{v^{\ep}}{u_{\ep}} A^{\ep}_i v^{\ep}=\int A^{\ep}_i \frac{v^{\ep}}{u_{\ep}} A^{\ep}_i \left(\frac{v^{\ep}}{u_{\ep}} u_{\ep}\right) %&=\int A^{\ep}_i \frac{v^{\ep}}{u_{\ep}} \left(u_{\ep}A^{\ep}_i \frac{v^{\ep}}{u_{\ep}}+\frac{v^{\ep}}{u_{\ep}}A^{\ep}_iu_{\ep}\right)
    =\int \left(A^{\ep}_i \frac{v^{\ep}}{u_{\ep}}\right)^2 u_{\ep}+\frac{1}{2}\int A^{\ep}_i \left(\frac{v^{\ep}}{u_{\ep}}\right)^2A^{\ep}_i u_{\ep}.
\end{split}
\end{equation*}

Note that for each vector field $X$, there holds
\begin{equation*}
\left(X\frac{v^{\ep}}{u_{\ep}}\right)v^{\ep}=\left(X\frac{v^{\ep}}{u_{\ep}}\right)\frac{v^{\ep}}{u_{\ep}}u_{\ep}=\frac{1}{2}X\left(\frac{v^{\ep}}{u_{\ep}}\right)^2 u_{\ep}.
\end{equation*}
Applying this to \eqref{eqn-Jan-1-4} yields
\begin{equation}\label{eqn-Jan-1-5}
\begin{split}
\frac{d}{dt}\chi^2(\nu^{\ep}_t, \mu_{\ep})&=-\ep^2 \int \left(A^{\ep}_i \frac{v^{\ep}}{u_{\ep}}\right)^2 u_{\ep}-\frac{\ep^2}{2}\int A^{\ep}_i \left(\frac{v^{\ep}}{u_{\ep}}\right)^2A^{\ep}_i u_{\ep}-\frac{\ep^2}{2}\int A^{\ep}_i \left(\frac{v^{\ep}}{u_{\ep}}\right)^2\div(A^{\ep}_i) u_{\ep} \\
&\quad+\int\left[\ep^{2}A^{\ep}_0\left(\frac{v^{\ep}}{u_{\ep}}\right)^2+B\left(\frac{v^{\ep}}{u_{\ep}}\right)^2\right]u_{\ep}.
\end{split}
\end{equation}

Recalling that $u_{\ep}\in W^{1,p}$ is the density of the stationary measure $\mu_{\ep}$ and satisfies \eqref{eqn-fpe}, we find 
\begin{equation*}
    -\frac{\ep^2}{2}\int A^{\ep}_i \left(\frac{v^{\ep}}{u_{\ep}}\right)^2 [A^{\ep}_i u_{\ep}+\div(A^{\ep}_i) u_{\ep}]+\int\left[\ep^{2}A^{\ep}_0\left(\frac{v^{\ep}}{u_{\ep}}\right)^2+B\left(\frac{v^{\ep}}{u_{\ep}}\right)^2\right]u_{\ep}=0,
\end{equation*}
which together with \eqref{eqn-Jan-1-5} implies \eqref{an-identity-2025-05-04}, completing the proof. 
\end{proof}

\subsection{Invariant measure selection by noise}\label{subsec-selection-by-noise}

Recall from Corollary \ref{cor-2024-03-06} that $\mu_0$ is physical if and only if it is ergodic. When $\mu_0$ fails to be ergodic or physical, the system \eqref{eqn-ode} could admit multiple invariant measures like $\mu_0$. In this subsection, we show that any of them can be selected by noise. 

\begin{defn} [Symmetric admissible class]
Let $u:M\to\R$ be positive. A collection of vector fields $A=\left\{A^{\ep}_{i},i\in\{0,\dots,m\},\,\,\ep\right\}$ is said to be in the {\em symmetric admissible class} $\mathcal{A}^{sym}_u$ if for each $\ep$, there exists a symmetric positive definite endomorphism $A^{sym}_{\ep}: M\to \mathrm{End}(TM)$ belonging to $W^{1,1}$ such that \begin{equation*}\label{condition-A}
\frac{1}{2}\sum_{i=1}^m(A_i^{\ep})^2 f+A^{\ep}_0f=\frac{1}{u}\div (A^{sym}_{\ep} \nabla f),\quad \forall f\in C^2(M).
\end{equation*}
\end{defn}

\begin{thm}\label{thm-converse}
Let $\mu$ be an invariant measure of \eqref{eqn-ode} and have a positive density $u\in W^{1,p}$ for some $p>d$. Then, $\mathcal{A}\cap \mathcal{A}^{sym}_u\neq \emptyset$, and for any $A\in \mathcal{A}^{sym}_u$ there holds $\MM_{A} = \{\mu\}$. 
\end{thm}
\begin{proof}[Proof of Theorem \ref{thm-converse}]

First, we show $\mathcal{A}^{sym}_u\neq \emptyset$. Note that the Nash embedding theorem asserts $M$ is isometrically embedded into $\R^{m}$ for some $m\in \N$. Thus, there exist $m$ smooth vector fields $\Tilde{A}_i$, $i\in \{1,\dots,m\}$, on $M$ so that $\sum_{i=1}^{m} (\Tilde{A}_i)^2=\Delta_M$, where $\Delta_M$ denotes the Laplace-Beltrami operator on $M$ (see e.g. \cite{Hsu}). For each $\ep$, we define 
$$
A^{\ep}_0:=\frac{1}{4u^2} (\Tilde{A}_i u)\Tilde{A}_i \quad\andd\quad A^{\ep}_i:=\frac{1}{\sqrt{u}} \Tilde{A}_i,\quad i\in \{1,\dots,m\}.
$$
As $u$ is positive and belongs to $W^{1,p}$ for some $p>d$, it is easy to see that  $A:=\left\{A^{\ep}_i, i\in \{0,\dots, m\},\,\,\ep\right\}\in \mathcal{A}$.  Straightforward calculations yield 
\begin{equation*}
    \frac{1}{2}\sum_{i=1}^m(A_i^{\ep})^2 f+A^{\ep}_0f=\frac{1}{2u}\sum_{i=1}^m(\Tilde{A}_i)^2f=\frac{1}{2u}\Delta_M f,\quad \forall f\in C^2(M),
\end{equation*}
resulting in $A\in \mathcal{A}^{sym}_u\cap \mathcal{A}$. 

Now, we prove $\MM_{A}=\{\mu\}$ for any $A\in \mathcal{A}^{sym}_u$. Fix such an $A$. Then, for each $\ep$ there exists a symmetric positive definite endomorphism $A^{sym}_{\ep}: M\to \mathrm{End}(TM)$ belonging to $ W^{1,1}$ such that 
$$
\frac{1}{2}\sum_{i=1}^m(A_i^{\ep})^2 f+A^{\ep}_0f=\frac{1}{u}\div (A^{sym}_{\ep} \nabla f),\quad \forall f\in C^2(M).
$$
Since $\mu$ is an invariant measure of \eqref{eqn-ode} and $u\in W^{1,p}$, we see from Lemma \ref{lem-new-measure-vectors} (1) with $u_0$ replaced by $u$ that $\div(uB)=0$. Consequently, an application of the divergence theorem yields 
\begin{equation*}
    \int (\LL_{\ep}f) u=\int \left(\frac{\ep^2}{u}\div(A^{sym}_{\ep}\nabla f)+Bf\right) u=0,\quad \forall f\in C^2(M). 
\end{equation*}
This together with Theorem \ref{lem-stationary-measure} implies $\mu_{\ep}=\mu$, and hence, $\MM_{A}=\{\mu\}$.
\end{proof}

\begin{rem}\label{rem-selection-by-noise}
If the system \eqref{eqn-ode} has multiple invariant measures similar to $\mu_0$, then Theorem \ref{thm-converse} says that they are all selectable, and therefore, none of them (and none of the invariant measures of \eqref{eqn-ode}) are stochastically stable with respect to $\mathcal{A}$. Hence, Theorem \ref{thm-converse} can be regarded as a result towards stochastic instability.
\end{rem}

\section{\bf The one-dimensional case}\label{sec-1D}

We study the one-dimensional case in this section. Given that any smooth, connected, and closed one-dimensional manifold is smoothly diffeomorphic to a circle (see e.g. \cite{MR0348781}), we consider $M$ to be the circle $\S^1$ for the sake of clarity.

%%%%%%%%%%%%%

\subsection{Setup and results}

Consider the following one-dimensional ODE over $\S^1$:
\begin{equation}\label{eqn-ode-1D}
    \dot x=B(x),
\end{equation}
where $B:\S^1\to T\S^1$ is a Lipchitz continuous vector field. Note that a vector field over $\S^{1}$ is naturally identified with a function on $\S^{1}$. In the sequel, the same notation is used for a vector field over $\S^{1}$ and its identification as a function on $\S^{1}$. This shall cause no trouble. It is assumed that \eqref{eqn-ode-1D} is conservative or generalized volume-preserving in the sense that it admits an invariant measure $\mu_0$ with a positive density function $u_0$. Consequently, it must hold that $B$ is either equal to zero ($B\equiv0$), positive ($B>0$), or negative ($B<0$). The $B\equiv0$ case is of no interest, and the other two cases are essentially the same. To maintain the clarity, we here focus on the $B>0$ case. It is then easy to see that \eqref{eqn-ode-1D} is uniquely ergodic and $u_0=\frac{\ga}{B}\in W^{1,\infty}$, where $\ga=(\int \frac{1}{B})^{-1}$ is the normalizing constant.

Consider \eqref{eqn-ode-1D} under small random perturbations: 
\begin{equation}\label{eqn-sde-1D}
    dX^{\ep}_t=B(X^{\ep}_t) dt+\ep^2 A^{\ep}_0(X^{\ep}_t)dt+\ep \sum_{i=1}^m A^{\ep}_i(X^{\ep}_t)\circ dW^i_t,
\end{equation}
where $0<\ep\ll 1$ is the noise intensity, $m\geq 1$, $A=\{A_i^{\ep},\,\,i\in \{0,1,\dots, m\}, \,\,\ep\}$ is a collection of vector fields on $\S^1$, and $\{W^i_t\}$ are $m$ independent and standard one-dimensional Brownian motions on some probability space. The stochastic integrals are understood in the sense of Stratonovich. The collection of vector fields $A$ is taken from the admissible class defined as follows.
\begin{defn}[Admissible class-$1$D]\label{admissible-class-1D}
    A collection $A$ of vector fields on $\S^1$ is said to be in the {\em admissible class $\mathcal{A}_{1D}$} if $A=\left\{A^{\ep}_{i},i\in\{0,\dots,m\},\,\,\ep\right\}$ for some $m\geq 1$ and the following conditions are satisfied:
    \begin{enumerate}
\item [\bf(A1)$_{1D}$] there exists $p>2$ such that $A^{\ep}_0\in L^p$, $A^{\ep}_i\in W^{1,p}$ for $i\in\{1,\dots,m\}$, and  
$$
\|A^{\ep}_0\|_p+ \max_i\|A^{\ep}_i\|_{1,p}\lesssim1;
$$

\item [\bf(A2)$_{1D}$] $\min\sum_{i=1}^m|A^{\ep}_i|^2\gtrsim1$.
\end{enumerate}
\end{defn}

\begin{rem}
    The only point that Definition \ref{admissible-class-1D} is not consistent with Definition \ref{admissible-class} lies in the requirement $p>2$ instead of $p>d=1$ in {\bf(A1)$_{1D}$}.
\end{rem}

% The Fokker-Planck equation associated with \eqref{eqn-sde-1D} reads
% \begin{equation}\label{eqn-fpe-1D}
%     \pa_t u=\frac{\ep^2}{2}\sum_{i=1}^m\left((A^{\ep}_i)^2 u\right)''-\left[\left(B+\ep^2 A_0^{\ep}+ \frac{\ep^2}{2}\sum_{i=1}^m A_i^{\ep}(A_i^{\ep})'\right)u\right]' \quad \text{on}\quad \S^1. 
% \end{equation}

If $A\in\mathcal{A}_{1D}$, Theorem \ref{lem-stationary-measure} is applied to yield that for each $\ep$, \eqref{eqn-sde-1D} admits a unique stationary measure $\mu_{\ep}$. Moreover, $\mu_{\ep}$ admits a positive density $u_{\ep}\in W^{1,p}$, where $p>2$ is the same as in {\bf (A1)$_{1D}$}. In addition, $u_{\ep}$ satisfies the stationary Fokker-Planck equation:
\begin{equation}\label{eqn-fpe-1D-stationary}
    \frac{\ep^2}{2}(a_{\ep} u_{\ep})''-[(B+\ep^2 b_{\ep}) u_{\ep}]'=0\quad  \text{in the weak sense,}
\end{equation}
or equivalently, 
\begin{equation*}
    -\frac{\ep^2}{2}\int  a_{\ep}u'_{\ep}f'+ \int \left(-\frac{\ep^2}{2}a'_{\ep} +B+\ep^2 b_{\ep}\right) u_{\ep}f'=0,\quad \forall f\in W^{1,2},
\end{equation*}
where $a_{\ep}:=\sum_{i=1}^m |A_i^{\ep}|^2$ and $b_{\ep}:=A_0^{\ep}+ \frac{1}{2}\sum_{i=1}^m A_i^{\ep}(A_i^{\ep})'$.

Since $\mu_0$ is the unique invariant measure of \eqref{eqn-ode-1D}, its stochastic stability with respect to $\mathcal{A}_{1D}$ follows readily. We are more interested in enhanced results, which are stated in the following two theorems.

%Recall that $\mu_0$ is uniquely ergodic under $\vp^t$. Hence, it is stochastically stable with respect to any noise perturbation. Nevertheless, one novelty of Theorem \ref{thm-uniform-estimate-stability} is the convergence of the densities $\lim_{\ep\to 0}u_{\ep}=u_0$ in some sense, leading to the strong convergence of $\lim_{\ep\to 0}\mu_{\ep}=\mu_0$ in the space of probability measures. There is no doubt that the proof of Theorem \ref{thm-uniform-estimate-stability} applies to the one-dimensional case. But, thanks to the significantly improved embedding properties of Sobolev spaces in the one-dimensional situation, the proof can be greatly simplified. We include it here for the reader's reference.   

\begin{thm}\label{thm-unifrom-bounds-1d}
For any $A\in \mathcal{A}_{1D}$, there hold
\begin{equation*}\label{uniform-estimates-1D}
\|u_{\ep}\|_{1,2}\lesssim1\quad\text{and}\quad 1\lesssim\min u_{\ep}\leq \max u_{\ep}\lesssim1.
\end{equation*}
%In addition, $\mu_0$ is uniquely ergodic under $\vp^t$ and stochastically stable with respect to $\mathcal{A}_{1D}$. 
Consequently, the limit $\lim_{\ep\to0}u_{\ep}=u_{0}$ holds weakly in $W^{1,2}$ and strongly in $C^{\al}$ for any $\al\in(0,\frac{1}{2})$.
\end{thm}

%A system's stochastic stability is in fact the outcome of the interplay between the deterministic dynamics and noise perturbations. In the proofs of Theorems \ref{thm-lower-bound-div-free} and \ref{thm-stronger-uniform-bounds-div-free}, which lie at the cores in proving Theorems \ref{thm-uniform-estimate-stability} and \ref{thm-stronger-uniform-bounds}, we only applied the dynamical property $\div B=0$ at the beginning when deriving those important estimates. As a comparison, when it comes to the one-dimensional situation, both the deterministic dynamics and the Fokker-Planck equation become rather simple. Leveraging the sign-preserving property of $B$ to investigate Kolmogorov's problem, we achieve the following two theorems.

\begin{thm}\label{thm-uniform-bounds-u_ep-1d}
Let $A=\{A_i^{\ep},\,\,i\in \{0,1,\dots, m\}, \,\,\ep\}$ be a collection of vector fields on $\S^{1}$ and satisfy 
\begin{itemize}
    \item $A^{\ep}_0\in C^0$, $A^{\ep}_i\in C^1$, $i\in \{1,\dots, m\}$, and $\lim_{\ep\to 0} \ep^2 \left(\|A^{\ep}_0\|_{\infty}+\sum_{i=1}^m\|A^{\ep}_i\|_{1,\infty}\right)=0$;  
    
\item $\min\sum_{i=1}^m |A^{\ep}_i|^2>0$ for each $\ep$.
\end{itemize}
Then, $1\lesssim \min u_{\ep}\leq \max u_{\ep}\lesssim1$ and the limit $\lim_{\ep\to 0}u_{\ep}=u_0$ holds weakly in $L^p$ for any $p>1$.

If, in addition, 
$$
A^{\ep}_0\in C^1, \,\, A^{\ep}_i\in C^2,\,\,i\in \{1,\dots, m\},\,\,\text{and}\,\,\lim_{\ep\to 0} \ep^2\left(\|A^{\ep}_0\|_{1,\infty}+\sum_{i=1}^m \|A^{\ep}_i\|_{2,\infty}\right)=0,
$$
then $\|u'_{\ep}\|_{\infty}\lesssim1$ and the limit $\lim_{\ep\to0}u_{\ep}=u_0$ holds in $C^{\al}$ for any $\al\in (0,1)$.
\end{thm}

\begin{rem}
    Theorem \ref{thm-unifrom-bounds-1d} is a one-dimensional counterpart of Theorem \ref{thm-uniform-estimate-stability}(1). Uniform-in-$\ep$ bounds of $A_0^{\ep}$ and $A_{i}^{\ep}$ in {\bf(A1)} and {\bf(A1)$_{1D}$} play important roles in the respective proofs of these two theorems.
    
    The aim of Theorem \ref{thm-uniform-bounds-u_ep-1d} is to extend the scope of Theorem \ref{thm-unifrom-bounds-1d} by relaxing the requirement for uniform-in-$\ep$ bounds of $A_0^{\ep}$ and $A_{i}^{\ep}$. An essential factor that enables this relaxation is the positivity of $B$, which is a unique characteristic in dimension one.
    \end{rem}

%%%%%%%%%%%%%%%%%%%%%%%%%%%

\subsection{Proof of Theorems \ref{thm-unifrom-bounds-1d} and \ref{thm-uniform-bounds-u_ep-1d}}

First, we prove Theorem \ref{thm-unifrom-bounds-1d}, which naturally follows from the reasoning leading to Theorem \ref{thm-uniform-estimate-stability}(A). Given the one-dimensional nature of the problem, we present a significantly more concise proof.

\begin{proof}[Proof of Theorem \ref{thm-unifrom-bounds-1d}]
It suffices to establish the uniform estimates for $\{u_{\ep}\}_{\ep}$. Set $w_{\ep}:=Bu_{\ep}$. Since $B$ is Lipschitz continuous and positive, it is equivalent to show
\begin{equation}\label{uniform-estimate-for-w}
\|w_{\ep}\|_{1,2}\lesssim1\quad\text{and}\quad  1\lesssim\min w_{\ep}\leq \max w_{\ep}\lesssim1.    
\end{equation}

It is seen from \eqref{eqn-fpe-1D-stationary} that $w_{\ep}$ satisfies $\frac {\ep^2}{2} \left( \frac{a_{\ep}}{B} w_{\ep}\right)'' -\left[ w_{\ep}+\ep^2 \frac{b_{\ep}}{B}w_{\ep}\right]'=0$ in the weak sense. That is, 
\begin{equation}\label{eqn-Sep-4-1}
    -\frac{\ep^2}{2}\int  \left(\frac{a_{\ep}}{B} w_{\ep}\right)' f'+\int  \left[ w_{\ep}+\ep^2 \frac{b_{\ep}}{B}w_{\ep}\right] f'=0,\quad \forall f\in W^{1,2}.
\end{equation}
The proof of \eqref{uniform-estimate-for-w} is broken into three steps.

\medskip

\paragraph{\bf Step 1} We claim  $\|w_{\ep}\|_{1,2}\lesssim1$ and $\max w_{\ep}\lesssim1$.

Taking $f:=w_{\ep}$ in \eqref{eqn-Sep-4-1} yields
\begin{equation*}
     -\frac{\ep^2}{2}\int  \left(\frac{a_{\ep}}{B} w_{\ep}\right)' w'_{\ep}+\int  \left[ w_{\ep}+\ep^2 \frac{b_{\ep}}{B}w_{\ep}\right] w'_{\ep}=0.
\end{equation*}
Using the fact $\int w'_{\ep}w_{\ep}=\frac{1}{2}\int (w^2_{\ep})'=0$, we find
\begin{equation}\label{eqn-Sep-4-2}
    \int \frac{a_{\ep}}{B} |w'_{\ep}|^2=\int  \left[-\left(\frac{a_{\ep}}{B}\right)'+2 \frac{b_{\ep}}{B}\right]w_{\ep}w'_{\ep}.
\end{equation}

It follows from {\bf (A1)$_{1D}$}, {\bf (A2)$_{1D}$} and the Lipschitz continuity and positivity of $B$ that
\begin{equation}\label{eqn-Sep-4-3}
    \left\|\frac{a_{\ep}}{B}\right\|_{1,p}+\left\|\frac{b_{\ep}}{B}\right\|_{p}\lesssim1\quad \andd\quad \min \frac{a_{\ep}}{B}\gtrsim 1,
\end{equation}
for some $p>2$. Applying H\"older's inequality to the right-hand side of  \eqref{eqn-Sep-4-2} results in
\begin{equation*}
    \|w'_{\ep}\|^2_2\lesssim\int \frac{a_{\ep}}{B} |w'_{\ep}|^2\leq \left\|-\left(\frac{a_{\ep}}{B}\right)'+2\frac{b_{\ep}}{B} \right\|_{p} \|w_{\ep}\|_{r}\|w'_{\ep}\|_2\lesssim\|w_{\ep}\|_{r}\|w'_{\ep}\|_2,
\end{equation*}
where $r:=\left(\frac{1}{2}-\frac{1}{p}\right)^{-1}$. Therefore, 
\begin{equation}\label{estimate-2023-10-10}
    \|w'_{\ep}\|_2\lesssim \|w_{\ep}\|_r
\end{equation}
An application of the Sobolev embedding theorem and the interpolation inequality then leads to 
\begin{equation*}
    \|w_{\ep}\|_{\infty}\lesssim \|w_{\ep}\|_2+\|w'_{\ep}\|_2\lesssim  \|w_{\ep}\|_r\lesssim \|w_{\ep}\|^{1-\frac{1}{r}}_{\infty}\|w_{\ep}\|_1^{\frac{1}{r}},
\end{equation*}
that is, $\|w_{\ep}\|_{\infty}\lesssim \|w_{\ep}\|_1$. Since $\|w_{\ep}\|_1\leq \|u_{\ep}\|_1\max B=\max B$, we find $\|w_{\ep}\|_{\infty}\lesssim1$, which together with \eqref{estimate-2023-10-10} yields $\|w'_{\ep}\|_2\lesssim1$.

\medskip

\paragraph{\bf Step 2} Set $v_{\ep}:=\ln w_{\ep}-\frac{1}{2\pi}\int  \ln w_{\ep}$. We prove $\|v_{\ep}\|_{\infty}\lesssim1$. 

Setting $f:=\frac{1}{w_{\ep}}$ in \eqref{eqn-Sep-4-1} results in
\begin{equation*}
    \frac{\ep^2}{2}\int  \left(\frac{a_{\ep}}{B} w_{\ep}\right)' \frac{w'_{\ep}}{w_{\ep}^{2}}-\int  \left[ w_{\ep}+\ep^2 \frac{b_{\ep}}{B}w_{\ep}\right]\frac{w'_{\ep}}{w_{\ep}^{2}}=0,
\end{equation*}
which together with $\int \frac{w'_{\ep}}{w_{\ep}}=\int  (\ln w_{\ep})'=0$ yields
\begin{equation*}
    \int  \frac{a_{\ep}}{B}\frac{|w'_{\ep}|^2}{w_{\ep}^{2}}=\int \left[-\left(\frac{a_{\ep}}{B} \right)'+2 \frac{b_{\ep}}{B}\right]\frac{w'_{\ep}}{w_{\ep}}.
\end{equation*}
Noting that $v'_{\ep}=\frac{w'_{\ep}}{w_{\ep}}$, we apply H\"older's inequality to find 
\begin{equation*}
    \int \frac{a_{\ep}}{B}|v'_{\ep}|^2=\int \left[-\left(\frac{a_{\ep}}{B} \right)'+2 \frac{b_{\ep}}{B}\right]v'_{\ep}\leq \left\|-\left(\frac{a_{\ep}}{B} \right)'+2 \frac{b_{\ep}}{B}\right\|_p \|1\|_r \|v'_{\ep}\|_2.
\end{equation*}
This together with \eqref{eqn-Sep-4-3}  leads to $\|v'_{\ep}\|_2\lesssim1$. A further application of the Sobolev embedding theorem and Poincar\'e inequality results in $\|v_{\ep}\|_{\infty}\lesssim \|v_{\ep}\|_2+\|v'_{\ep}\|_2\lesssim \|v'_{\ep}\|_2\lesssim1$.

\medskip

\paragraph{\bf Step 3} We show $\min w_{\ep}\gtrsim1$. 

Letting $f:=\frac{1}{w_{\ep}^{3}}$ in \eqref{eqn-Sep-4-1} yields 
\begin{equation*}
     \frac{\ep^2}{2}\int  \left(\frac{a_{\ep}}{B} w_{\ep}\right)' \frac{3w'_{\ep}}{w_{\ep}^{4}}-\int  \left[ w_{\ep}+\ep^2 \frac{b_{\ep}}{B}w_{\ep}\right]\frac{3w'_{\ep}}{w_{\ep}^{4}}=0.
\end{equation*}
Since $\int \frac{w'_{\ep}}{w_{\ep}^{3}}=-\frac{1}{2}\int  (\frac{1}{w_{\ep}^{2}})'=0$, we arrive at 
\begin{equation*}
    \int \frac{a_{\ep}}{B}\frac{|w'_{\ep}|^2}{w_{\ep}^{4}}=\int \left[-\left(\frac{a_{\ep}}{B} \right)'+2 \frac{b_{\ep}}{B}\right]\frac{w'_{\ep}}{w_{\ep}^{3}},
\end{equation*}
and hence, the fact $(\frac{1}{w_{\ep}})'=-\frac{w'_{\ep}}{w_{\ep}^{2}}$ and H\"{o}lder's inequality imply
\begin{equation*}
    \int  \frac{a_{\ep}}{B} \left|\left( \frac{1}{w_{\ep}}\right)'\right|^2=\int \left[\left(\frac{a_{\ep}}{B} \right)'-2 \frac{b_{\ep}}{B}\right]\frac{1}{w_{\ep}}\left(\frac{1}{w_{\ep}}\right)'\leq \left\|\left(\frac{a_{\ep}}{B} \right)'-2 \frac{b_{\ep}}{B}\right\|_p \left\|\frac{1}{w_{\ep}}\right\|_r \left\|\left(\frac{1}{w_{\ep}}\right)'\right\|_2.
\end{equation*}
It then follows from \eqref{eqn-Sep-4-3} that $\left\|\left(\frac{1}{w_{\ep}}\right)'\right\|\lesssim\left\|\frac{1}{w_{\ep}}\right\|_r$, which together with the Sobolev embedding theorem and interpolation inequality yields
\begin{equation}\label{eqn-Sep-4-4}
    \left\|\frac{1}{w_{\ep}}\right\|_{\infty}\lesssim \left\|\frac{1}{w_{\ep}}\right\|_{2}+\left\|\left(\frac{1}{w_{\ep}}\right)'\right\|_2\lesssim \left\|\frac{1}{w_{\ep}}\right\|_r.
\end{equation}

Since $\|w_{\ep}\|_r\left\|\frac{1}{w_{\ep}}\right\|_r=\left(\int  e^{r v_{\ep}}\int  e^{-r v_{\ep}}\right)^{\frac{1}{r}}\lesssim1$ by {\bf Step 2}, we find 
\begin{equation*}
    \left\|\frac{1}{w_{\ep}}\right\|_r\lesssim \frac{1}{\|w_{\ep}\|_r}\lesssim \frac{1}{(\min w_{\ep})^{1-\frac{1}{r}} \|w_{\ep}\|_1^{\frac{1}{r}}}\lesssim \left\|\frac{1}{w_{\ep}}\right\|_{\infty}^{1-\frac{1}{r}},
\end{equation*}
where we used $\frac{1}{(\min w_{\ep})^{1-\frac{1}{r}}}=\|\frac{1}{w_{\ep}}\|_{\infty}^{1-\frac{1}{r}}$ and $\|w_{\ep}\|_{1}\geq \|u_{\ep}\|_1\min B= \min B$ in the third inequality. 
Applying \eqref{eqn-Sep-4-4} then results in $\|\frac{1}{w_{\ep}}\|_{\infty}\lesssim1$. Hence, $\min w_{\ep}=\|\frac{1}{w_{\ep}}\|_{\infty}^{-1}\gtrsim1$.

Consequently, \eqref{uniform-estimate-for-w} follows from {\bf Step 1} and {\bf Step 3}, completing the proof. 
\end{proof}

Next, we prove Theorem \ref{thm-uniform-bounds-u_ep-1d} by means of Berstein-type estimates.

\begin{proof}[Proof of Theorem \ref{thm-uniform-bounds-u_ep-1d}]
Thanks to Theorem \ref{lem-stationary-measure}, $u_{\ep}$ belongs to $W^{1,p}$ for any $p>2$. Hence, $u_{\ep}\in C^0$ by the Sobolev embedding theorem.

We claim $u_{\ep}\in C^1$. Indeed, \eqref{eqn-fpe-1D-stationary} says particularly that $\frac{\ep^2}{2}(a_{\ep} u_{\ep})'-(B+\ep^2 b_{\ep}) u_{\ep}$ has weak derivative $0$. Thus, $\frac{\ep^2}{2}(a_{\ep} u_{\ep})'-(B+\ep^2 b_{\ep}) u_{\ep}$ is absolutely continuous (up to a set of zero Lebesgue measure) and there exists $C_{\ep}\in \R$ such that 
\begin{equation}\label{eqn-Aug-31-1}
    \frac{\ep^2}{2}(a_{\ep} u_{\ep})'-(B+\ep^2 b_{\ep}) u_{\ep}=C_{\ep}\quad \text{a.e. on}\quad \S^1.
\end{equation}
Obviously, $a_{\ep}\in C^1$ and $b_{\ep}\in C^0$. Given $\min a_{\ep}>0$ (by assumption), the continuity of $B$ and $u_{\ep}\in C^0$, we find from  \eqref{eqn-Aug-31-1} that $u'_{\ep}$ is a.e. equal to a continuous function. Thus, we may assume, without loss of generality, that $u_{\ep}\in C^1$. 

Integrating \eqref{eqn-Aug-31-1} over $\S^1$ yields 
\begin{equation}\label{eqn-Sep-4-6}
    -\int  (B+\ep^2 b_{\ep})u_{\ep}=C_{\ep}.
\end{equation}
Note that the assumptions ensure that 
\begin{equation}\label{eqn-Sep-4-5}
\lim_{\ep\to0}\ep^2\left(\|a_{\ep}\|_{\infty}+\|b_{\ep}\|_{\infty}\right)=0.  
\end{equation}
Then, we see from \eqref{eqn-Sep-4-6}, $B>0$ and $\int  u_{\ep}=1$ that
\begin{equation}\label{limsupinf-2023-10-10}
C_{\ep}\approx -1  
\end{equation}

Suppose $u_{\ep}$ attains its maximum and minimum on $\S^1$ at $x_{\ep}$ and $y_{\ep}$, respectively. Then, $u'_{\ep}(x_{\ep})=0=u'_{\ep}(y_{\ep})$. Substituting them into \eqref{eqn-Aug-31-1} yields that the equality $\frac{\ep^2}{2}a'_{\ep} u_{\ep}-(B+\ep^2 b_{\ep}) u_{\ep}=C_{\ep}$ holds at $x_{\ep}$ and $y_{\ep}$, which together with \eqref{eqn-Sep-4-5} and \eqref{limsupinf-2023-10-10} gives 
\begin{equation*}\label{eqn-Aug-31-2}
    u_{\ep}(x_{\ep})= \frac{C_{\ep}}{\frac{\ep^2}{2}a'_{\ep}(x_{\ep})-(B+\ep^2 b_{\ep})(x_{\ep})}\lesssim1\quad\text{and}\quad u_{\ep}(y_{\ep})= \frac{C_{\ep}}{\frac{\ep^2}{2}a'_{\ep}(x_{\ep})-(B+\ep^2 b_{\ep})(x_{\ep})}\gtrsim1.
\end{equation*}

It remains to show the convergence result for $\{u_{\ep}\}_{\ep}$. Recall that $\mu_0$ is the unique invariant measure of \eqref{eqn-ode-1D}. Thanks to the compactness of $\S^1$ and Prokhorov's theorem, the limit $\lim_{\ep\to 0}\mu_{\ep}=\mu_0$ holds under the weak*-topology. This together with $\|u_{\ep}\|_{\infty}\lesssim1$ yields $\lim_{\ep\to 0}u_{\ep}=u_0$ weakly in $L^p$ for any $p>1$. 

\medskip

Now, we prove the results under additional assumptions. Note that $u_{\ep}\in C^{1}$ in this case. Since $\|u_{\ep}\|_{\infty}\lesssim1$, the classical theory for elliptic equations guarantees that $u_{\ep}\in W^{2,p}$ for any $p>2$. Obviously, $a_{\ep}\in C^2$, $b_{\ep}\in C^1$ and $\min a_{\ep}>0$. Since $B$ is Lipschitz continuous and $u_{\ep}\in C^1$, we see from \eqref{eqn-fpe-1D-stationary} that $u''_{\ep}\in C^0$. 

Suppose $(u'_{\ep})^2$ attains its maximum at $x_{\ep}\in \S^1$. Then, $(u'_{\ep}u''_{\ep})(x_{\ep})=0$. Multiplying \eqref{eqn-fpe-1D-stationary} by $u'_{\ep}$ and evaluating at $x_{\ep}$, we deduce 
\begin{equation}\label{eqn-Sep-4-7}
\begin{split}
    0&=\frac{\ep^2}{2}(a_{\ep}u_{\ep})'' u'_{\ep}-[(B+\ep^2 b_{\ep})u_{\ep}]' u'_{\ep}\\
    &=\frac{\ep^2}{2}\left(a_{\ep} u_{\ep}''u'_{\ep}+2a'_{\ep}(u'_{\ep})^2+a''_{\ep}u_{\ep}u'_{\ep}\right)-(B+\ep^2 b_{\ep})(u_{\ep}')^2-(B+\ep^2 b_{\ep})' u_{\ep} u'_{\ep}\\
    &=\left[\ep^2 a'_{\ep} -(B+\ep^2 b_{\ep})\right](u'_{\ep})^2+\left(\frac{\ep^2}{2}a''_{\ep}-(B+\ep^2 b_{\ep})'\right)u_{\ep}u'_{\ep} \quad \text{at }x_{\ep}.
\end{split}
\end{equation}

Since the assumptions ensure $\lim_{\ep\to0}\ep^{2}(\|a_{\ep}'\|_{\infty}+\|a_{\ep}''\|_{\infty}+\|b_{\ep}\|_{\infty}+\|b'_{\ep}\|_{\infty})=0$,
we conclude from $B>0$, \eqref{eqn-Sep-4-7} that $\|u'_{\ep}\|_{\infty}=|u'_{\ep}(x_{\ep})|\lesssim\|u_{\ep}\|_{\infty}\lesssim1$. It then follows from the Arzel\`a-Ascoli theorem that the limit $\lim_{\ep\to0}u_{\ep}=u_0$ holds in $C^{\al}$ for any $\al\in (0,1)$. 
\end{proof}

\section{\bf The bounded domain case}\label{section-bounded-domain-case}

Let $M\subset \R^d$ be an open, bounded, and smooth domain and consider the following ODE over $\ol{M}$:
\begin{equation}\label{ode-bounded-domain}
    \dot{x}=B(x),\quad x\in \ol{M},
\end{equation}
where $B:\ol{M}\to\R^{d}$ is a Lipschitz continuous vector field and satisfies $B\cdot \vec{n}=0$ on $\pa M$, where $\vec{n}$ is the inward unit normal vector field along $\partial M$. It is assumed that \eqref{ode-bounded-domain} is conservative (or generalized volume-preserving) in the sense that it admits an invariant measure $\mu_0$ with a uniformly positive density function $u_0\in W^{1,p_0}(M)$ for some $p_0>d$. In the current setting, it is more common to rewrite $B$ as a vector-valued function $(B_i)$. However, for notation compatibility, we retain  the tensor form and do not distinguish the two representations. This shall cause no confusion. 

Adding small random perturbations to \eqref{ode-bounded-domain} and assuming the resulted stochastic process is obliquely reflected at the boundary, we arrive at the following reflected SDE:
\begin{equation}\label{sde-bounded-domain}
    dX^{\ep}_t=B(X^{\ep}_t)dt+\ep^2 A^{\ep}_0(X^{\ep}_t)dt+\ep\sum_{i=1}^m A^{\ep}_i(X^{\ep}_t)\circ dW^i_t+N^{\ep}(X^{\ep}_t)dl^{\ep}_t,
\end{equation}
where $0<\ep\ll 1$ is the noise intensity, $m\geq d$, $A=\{A^{\ep}_i=(A^{\ep}_{ij}):\ol{M}\to \R^d,\,\, i\in \{0,1\dots, m\},\,\,\ep\}$ is a collection of vector fields on $\ol{M}$, and $\{W^i_t\}$ are $m$ independent and standard one-dimensional Brownian motions on some probability space, $N^{\ep}:=\frac{1}{2}(A^{\ep}_{ij}A^{\ep}_{ik}n_j):\pa M\to \R^d$ is the conormal vector field along $\pa M$, $n_j$ is the $j$-th component of $\vec{n}$, and $l^{\ep}_t$ is an increasing, non-negative, and continuous process that increases only when $X^{\ep}_t\in \pa M$. The stochastic integrals with respect to $W^i_t$ are interpreted in the sense of Stratonovich, while $\int_0^t N^{\ep}(X^{\ep}_s)dl^{\ep}_s$ is understood as a Stieltjes integral. We mention that  a solution of \eqref{sde-bounded-domain} refers to the pair $(X^{\ep}_t, l^{\ep}_t)$. 

The distributions of $X^{\ep}_t$ satisfies the following Fokker-Planck equation with reflecting boundary condition:
\begin{equation}\label{fpe-bounded-domain}
\begin{cases}
    \pa_t v =\LL^*_{\ep} v \quad \text{in}\quad M,\\
    \left(\frac{\ep^2}{2}\pa_k(A^{\ep}_{ij}A^{\ep}_{ik}v)-\frac{\ep^2}{2}A^{\ep}_{ik}\pa_k A^{\ep}_{ij} v-\ep^2 A^{\ep}_{0j}v-B_jv\right)n_j=0 \quad \text{on}\quad \pa M,
\end{cases}
\end{equation}
where $\LL^*_{\ep}$ is the Fokker-Planck operator and has the same form as in \eqref{eqn-fpe-time}. We refer the reader to Appendix \ref{app-subsec-A2} for more details about reflected SDEs and associated Fokker-Planck equations.

We consider less regular coefficients and still assume $A\in \mathcal{A}$ except that $M$ now is an open, bounded, and smooth domain rather than a closed Riemannian manifold. As in the closed Riemannian manifold case, the SDE \eqref{sde-bounded-domain} is hardly well-posed, and therefore, we work with the Fokker-Planck equation \eqref{fpe-bounded-domain}. It follows from Theorem \ref{lem-stationary-measure-on-bounded-domain} that for each $\ep>0$, \eqref{sde-bounded-domain} admits a unique stationary measure $\mu_{\ep}$ with a positive density $u_{\ep}\in W^{1,q}$ for any $q\in[1,p)$, where $p$ is as assumed in {\bf (A1)}. 

Regarding $\{\mu_{\ep}\}_{\ep}$ as probability measures on the compact set $\ol{M}$, the tightness of $\{\mu_{\ep}\}_{\ep}$ holds automatically. We still denote by $\MM_A$ the set of all limiting measures of $\mu_{\ep}$ as $\ep\to0$ under the weak*-topology. Then, all the elements in $\MM_{A}$ are supported on $\ol{M}$. The following theorem collects results parallel to those for the case of a closed Riemannian manifold. 

\begin{thm}
The following hold.
\begin{itemize}
    \item [(1)] For any $A\in \mathcal{A}$, there hold
$$
\|u_{\ep}\|_{1,2}\lesssim1\quad\text{and}\quad 1\lesssim\min u_{\ep}\leq \max u_{\ep}\lesssim1.
$$
In particular, any $\mu\in\MM_{A}$ has a density $u$ belonging to $W^{1,2}$ and satisfying $u,\frac{1}{u}\in L^{\infty}$.

\item [(2)] If $\mu_0$ is physical, then it is stochastically stable with respect to $\mathcal{A}$. 

\item [(3)] If $\mu_0$ is the only invariant measure of \eqref{ode-bounded-domain} with a density in $W^{1,2}$, then it is stochastically stable with respect to $\mathcal{A}$. 

\item [(4)] Let $A\in \mathcal{A}$ satisfy {\bf (A1)} with $p>d+2$ and $v_0\in L^2$ be a probability density. Then, for each $\ep>0$, the problem \eqref{fpe-bounded-domain} with initial condition $v(\cdot, 0)=v_0$ admits a unique solution $v^{\ep}\in V_{2,loc}^{1,0}(M\times [0,\infty))$, which is positive a.e. in $M\times (0,\infty)$, and satisfies $\int v^{\ep}(\cdot,t)=1$ for all $t\in(0,\infty)$ and $\|v^{\ep}\|_{L^{\infty}(M\times I)}<\infty$ for any closed interval $I\subset (0,\infty)$. Moreover, there exists $C>0$ (independent of $v_0$) such that $\chi^2(\nu^{\ep}_t, \mu_{\ep})\leq e^{-C\ep^2 t}\chi^2(\nu_0, \mu_{\ep})$ for all $t\geq 0$ and $0<\ep\ll1$, where $d\nu^{\ep}_t:=v^{\ep}(\cdot, t)dx$. 
\end{itemize}
In the case of either (2) or (3), the limit $\lim_{\ep\to0}u_{\ep}=u_{0}$ holds weakly in $W^{1,2}$ and strongly in $L^{p}$ for any $p\in[1,\infty)$.

\end{thm}
\begin{proof}
Once (1) is established, (2)-(3) follow from arguments analogous to those used in proving Theorem \ref{thm-uniform-estimate-stability}. Although the regularity of $v^{\ep}$ is a little weaker than that in Theorem \ref{thm-chi-convergence}, the arguments in proving Lemma \ref{lem-Cauchy-integral-formula} and Theorem \ref{thm-chi-convergence} can be easily adapted to verify (4). 

It suffices to prove (1). For clarity, we here only deal with the case with $\div B=0$ so that $u_0\equiv 1$; the general case follows in the same way as in the proof of Theorem \ref{thm-uniform-estimate-stability} (1) for the reason that the transform in Section \ref{subsec-transform} ensures $\tilde{B}\cdot \vec{n}=0$ on $\pa M$.

Note that a key step in the proofs of Lemmas \ref{lem-uniform-L2-gradient} and \ref{lem-upper-bound} and Theorem \ref{thm-lower-bound-div-free} is to show that 
\begin{equation}\label{eqn-B-int-by-parts}
\int [Bf(w)] w=0,\quad\forall f\in C^1(\R)
\end{equation}
with $w=u_{\ep}$ in the proof Lemmas \ref{lem-uniform-L2-gradient} and \ref{lem-upper-bound} and $w=v_{\ep}$ in the proof of Theorem \ref{thm-lower-bound-div-free}. In the current setting where $M\subset \R^d$ is an open, bounded, and smooth domain, we apply the divergence theorem and the facts that $\div B=0$ in $M$ and $B\cdot \vec{n}=0$ on $\pa M$ to see that \eqref{eqn-B-int-by-parts} still holds. In fact, for $f\in C^1(\R)$,
\begin{equation*}
\begin{split}
   \int_M [Bf(w)] w&=\int_{\pa M} B\cdot \vec{n} f(w)w-\int_M (\div B) f(w) w-\int_{M} f(w)Bw\\
   &=\int_M B_j\pa_jF(w)= \int_{\pa M} B\cdot \vec{n} F(w)-\int_{M} (\div B) F(w)=0,
\end{split}
\end{equation*}
where $F\in C^2(\R)$ is an antiderivative of $f$. Given the Sobolev embedding theorem, Poincar\'e inequality, and $u_{\ep}>0$ on $\ol{M}$, it is not hard to follow the proofs of Lemmas \ref{lem-uniform-L2-gradient} and \ref{lem-upper-bound} and Theorem \ref{thm-lower-bound-div-free} to derive (1). This completes the proof. 
\end{proof}

%%%%%%%%%%%%%%%%%%%%%%%%%%%%%%%%%%%

\section*{Acknowledgement}

We would like to express gratitude to Professor Jianyu Chen at Soochow University for insightful discussions regarding Theorem \ref{thm-construction}.

%%%%%%%%%%%%%%%%%%%%%%%%%%%%%%%%%%%
\appendix

\section{\bf Fokker-Planck equations with less regular coefficients}\label{app-stationary-measure}

In this appendix, we present some results about Fokker-Planck equations (FPEs) with less regular coefficients on compact manifolds and bounded domains with reflecting boundary conditions. They include the existence, uniqueness, and regularity of stationary measures, as well as the well-posedness of Cauchy problems.

\subsection{FPEs on closed Riemannian manifolds}\label{app-subsec-A1}

Let $d\geq1$ be an integer and $M$ be a $d$-dimensional smooth, connected, and closed Riemannian manifold. Consider the following SDE on $M$:
 \begin{equation}\label{sde-manifold-appendix}
 dX_t=A_0(X_t)dt+\sum_{i=1}^{m} A_{i}(X_t)\circ dW^i_t,
 \end{equation}
 where $m\geq d$, $A_{i}$, $i\in\{0,\dots,m\}$ are vector fields on $M$, $W^{i}_{t}$, $i\in\{1,\dots,m\}$ are independent and standard one-dimensional Brownian motions on some probability space, and the stochastic integrals are understood in the sense of Stratonovich. Denote by $\LL:=\frac{1}{2}\sum_{i=1}^m A_i^2 + A_0$ the generator of \eqref{sde-manifold-appendix}. Its formal $L^{2}$-adjoint operator, namely, the Fokker-Planck operator, is denoted by $\LL^{*}$.

Assume that $A_0\in L^{p}$ and $A_{i}\in W^{1,p}$, $i\in\{1,\dots,m\}$ for some $p>\max\{d,2\}$, and there exists $\la>0$ such that 
 \begin{equation*}
 \sum_{i=1}^m|A_i f|^2\geq \la |\nabla f|^2\quad {\rm Vol}\text{-a.e.}, \quad \forall f\in W^{1,1}.
 \end{equation*}

 \begin{defn}[Stationary measure]\label{def-stationary-dist}
 A probability measure $\mu$ on $M$ is called a {\em stationary measure} of \eqref{sde-manifold-appendix} if $A_0$, $\nabla A_i\in L^1(M,\mu)$ for $i\in \{1,\dots, m\}$, and $\LL^{*}\mu=0$ in the sense that $\int \LL f d\mu=0$ for all $f\in C^{2}$.

% Furthermore, $\mu_{\ep}$ is called {\em regular} if it has a density. %with respect to $dV$. Namely, there exists a measurable function $u_{\ep} \geq 0$ $\mu_g$-a.e. such that $d\mu_{\ep}=u_{\ep}dV$. 
 \end{defn}

 \begin{thm}[\cite{MR0928950,MR1876411}]\label{lem-stationary-measure}
     The SDE \eqref{sde-manifold-appendix} admits a unique stationary measure $\mu$, which has a positive density $u\in W^{1,p}$. Moreover, there holds
     \begin{equation*}
 -\frac{1}{2}\int A_i f\left[A_i u+(\div A_i) u\right]+\int\left(A_0f\right)u=0,\quad \forall f\in W^{1,2}.
 \end{equation*}
 \end{thm}

\begin{thm}[\cite{MR241822,MR1876411}]\label{thm-Cauchy-appx-manifold}
Assume $p>d+2$. For each initial probability density $v_0\in L^2$, the equation $\pa_t v =\LL^* v$ admits a unique solution $v\in V_{2,loc}^{1,0}(M\times [0,\infty))\cap \HH^{1,p}_{loc}(M\times (0,\infty))$ satisfying 
\begin{equation*}
    \int \phi(\cdot, t)v(\cdot,t)=\int \phi(\cdot,0)v_0+\int_0^t\int\left( \pa_t \phi v-\frac{1}{2}A_i\phi[A_i v+(\div A_i) v]+A_0\phi v\right)
\end{equation*} 
for any $\phi\in C^{1,1}(M\times [0,\infty))$ and $t\in (0,\infty)$. Moreover, $v$ is continuous and positive in $M\times (0,\infty)$ and satisfies $\int v(\cdot,t)=1$ for all $t\in (0,\infty)$.
\end{thm}

%%%%%%%%%%%%%%

\subsection{FPEs on bounded domains}\label{app-subsec-A2}

Let $M\subset\R^d$ be an open, bounded, and smooth domain and consider the following obliquely reflected SDE:
\begin{equation}\label{app-reflected-sde}
    dX_t=A_0(X_t)dt+\sum_{i=1}^m A_i(X_t)\circ d W^i_t+N(X_t) dl_t,
\end{equation}
where $m$, $A_i$, and $W_i$ are as those in Subsection \ref{app-subsec-A1}, $N=\frac{1}{2}(A_{ij}A_{ik}n_j):\pa M\to \R^d$ is the conormal vector field along $\pa M$, $A_{ij}$ for $j\in \{1,\dots, d\}$ is the $j$-th component of the vector $A_i$, namely, $A_i=A_{ij}\pa_j$, $n_j$ is the $j$-th component of the inward unit normal vector $\vec{n}$ along $\pa M$, and $l_t$ is an increasing, non-negative, and continuous process that increases only when $X_t\in \pa M$. Intuitively, $X_t$ evolves according to \eqref{app-reflected-sde} without the term $N(X_t) dl_t$ when it is inside $M$. Once touching $\pa M$, $X_t$ is reflected along the direction $N$ into $M$.

A solution of \eqref{app-reflected-sde} refers to the pair $(X_t, l_t)$. When all the coefficients in \eqref{app-reflected-sde} are smooth, the classical theory (see e.g. \cite[Theorem 2.4.1]{pilipenko2014}) ensures the existence and uniqueness of solutions. The conservation of probability is then a simple consequence of the reflection on $\partial M$. Moreover, an application of It\^o's formula indicates that for $f\in C^2(\ol{M})$ obeying $N\cdot \nabla f=0$, the average $(x,t)\mapsto\E_{x}[f(X_t)]$ satisfies (see e.g. \cite[Theorem 3.1.1]{pilipenko2014})
\begin{equation*}
\begin{cases}
    \pa_t w=\LL w &\text{in}\quad M\times (0,\infty),\\
    N\cdot \nabla w=0 &\text{on}\quad \pa M,
\end{cases}
\end{equation*}
where $\LL$ is the generator of $X_t$ having the same form as in Subsection \ref{app-subsec-A1}. Standard duality arguments yield that the distribution of $X_t$ satisfies
\begin{equation}\label{app-fpe-reflected}
    \begin{cases}
        \pa_t v=\LL^* v &\text{in}\quad M\times(0,\infty),\\
        \BB^* v=0 &\text{on}\quad \pa M,
    \end{cases}
\end{equation}
where $\LL^*$ is the Fokker-Planck operator as in Subsection \ref{app-subsec-A1} and \begin{equation}\label{zero-flux-boundary-appendix}
    \BB^*v:=\left(\frac{1}{2}\pa_k (A_{ij}A_{ik} v)-\frac{1}{2}A_{ik}\pa_k A_{ij}v-A_{0j}v\right)n_j.
\end{equation}
%We mention that it is straightforward to verify analytically that the reflecting boundary condition in \eqref{app-fpe-reflected} yields no loss of probability. In fact, integrating \eqref{fpe-appendix} over $M$, we find from integration by parts formula and the reflecing boundary condition that $\frac{d}{dt}\int v=0$. 

We consider \eqref{app-reflected-sde} with less regular coefficients, that is,  $A_i$ for $i\in \{0,\dots, m\}$ satisfies the same conditions as in Subsection \ref{app-subsec-A1} except that $M$ here is an open, bounded, and smooth domain rather than a closed Riemannian manifold. The following theorem establishes the well-posedness of \eqref{app-fpe-reflected}. 

\begin{thm}\label{thm-Cauchy-appx}
Assume $p>d+2$. For each initial probability density $v_0\in L^2$, the problem \eqref{app-fpe-reflected} admits a unique solution $v\in V_{2,loc}^{1,0}(M\times [0,\infty))$ satisfying 
\begin{equation}\label{eqn-Mar-8-1}
    \int \phi(\cdot, t)v(\cdot,t)=\int \phi(\cdot,0)v_0+\int_0^t\int\left( \pa_t \phi v-\frac{1}{2}A_i\phi[A_i v+(\div A_i) v]+A_0\phi v\right)
\end{equation} 
for any $\phi\in C^{1,1}(M\times [0,\infty))$ and $t\in (0,\infty)$. Moreover, $v$ is positive a.e. in $M\times (0,\infty)$, and satisfies $\int v(\cdot,t)=1$ for all $t\in (0,\infty)$ and $\|v\|_{L^{\infty}(M\times I)}<\infty$ for any closed interval $I\subset (0,\infty)$.
\end{thm}

\begin{proof}
The unique existence of a solution $v\in V_{2,loc}^{1,0}(M\times [0,\infty))$ satisfying \eqref{eqn-Mar-8-1} follows from the classical theory for weak solutions of parabolic equations (see e.g. \cite{MR241822, MR1465184}). Setting $\phi\equiv 1$ in \eqref{eqn-Mar-8-1} yields $\int v(\cdot, t)=\int v_0=1$ for all $t\in(0,\infty)$. 

Since $v_0\geq 0$, we apply \cite[Theorem 6.40]{MR1465184} to $-v$ to find that $-v\leq 0$ a.e. in $M\times [0,\infty)$, leading to the non-negativeness of $v$. Thus, an application of the weak Harnack inequality \cite[Theorem 6.42]{MR1465184} yields that $v$ is positive a.e. in $M\times (0,\infty)$. Let $I\in (0,\infty)$ be a closed interval. It follows from the combination of \cite[Theorems 6.41 and 6.38]{MR1465184} that $\|v\|_{L^{\infty}(M\times I)}<\infty$. The proof is finished. 
\end{proof}
\begin{rem}
It is assumed in \cite{MR1465184} that all the lower order coefficients in the parabolic equations are bounded. This is not essential since one may follow \cite[Chapter III]{MR241822} to generalize the results to the case where the coefficients are merely $L^p$-integrable for suitably large $p$. In the situation of Theorem \ref{thm-Cauchy-appx}, $p>d+2$ is sufficient.  
\end{rem}

\begin{defn}[Stationary measure]
A probability measure $\mu$ with density $u\in W^{1,2}$ is called a {\em stationary measure} of \eqref{app-reflected-sde} if $\mu$ is a stationary solution of \eqref{app-fpe-reflected} in the sense that
\begin{equation*}
-\frac{1}{2}\int A_i \phi\left[A_i u+(\div A_i) u\right]+\int\left(A_0\phi\right)u=0,\quad \forall \phi\in C^{\infty}.
\end{equation*}
%Furthermore, $\mu_{\ep}$ is called {\em regular} if it has a density. %with respect to $dV$. Namely, there exists a measurable function $u_{\ep} \geq 0$ $\mu_g$-a.e. such that $d\mu_{\ep}=u_{\ep}dV$. 
\end{defn}

\begin{rem}
We point out that Definition  \ref{def-stationary-dist} is not directly applicable in the present setting. In fact, to account for the reflecting boundary condition in \eqref{app-fpe-reflected}, it is necessary to impose suitable boundary conditions on test functions. To avoid such technical complications, we adopt a stronger formulation to define stationary measures. Nonetheless, this is sufficient for our purposes. 
\end{rem}

The following result summarizes the existence, uniqueness, and regularity of stationary measures.  

\begin{thm}[\cite{MR3623550}]\label{lem-stationary-measure-on-bounded-domain}
The SDE \eqref{app-reflected-sde} admits a unique stationary measure $\mu$ with a positive density $u\in W^{1,q}$ for any $q\in[1,p)$. 
\end{thm}

%%%%%%%%%%%%%%%%%%%%%%%%%%%%%%
\section{\bf Calculus on manifolds}\label{app-formulas}

We collect some of the frequently used formulas for calculus on manifolds. 
\begin{lem}\label{lem-formulas}
Let $f,h:M\to\R$ belong to $W^{1,2}$ and $X:M\to TM$ be a vector field in $W^{1,2}$. Then, the following hold:
\begin{enumerate}
    \item $X(fh)=(Xf) h+fXh$;
    \item $\div (fX)=Xf+f\div X$;
    \item$\int Xf=-\int f\div X$;
    \item if $X\in W^{1,p}$ for some $p>\frac{d}{2}$, then $\int h Xf =-\int fh (\div X)-\int fXh$.
\end{enumerate}
\end{lem}
\begin{proof}
(1)-(2) are obvious.
Note that $\int Xf=-\int \div(fX) -\int f \div X$. Since $M$ has no boundary, we apply the divergence theorem to find $\int \div (fX)=0$. Hence, (3) holds. Similarly, (4) also holds. 
\end{proof}

The next elementary result concerns vector fields spanning the tangent bundle $TM$.

\begin{lem}\label{lem-vector-basis}
Let $X_i:M\to TM$, $i\in \{1,\dots, n\}$ for some $n\geq d$, be continuous vectors fields on $M$. Assume $\{X_i\}_{i=1}^n$ spans the tangent bundle $TM$. Then, there is $C>1$ such that for any $f\in W^{1,1}$, there holds 
$$
\frac{1}{C}|\nabla f|\leq \sum_{i=1}^n |X_i f|\leq C|\nabla f|\quad \text{\rm Vol-a.e. on}\,\,M.
$$
\end{lem}
\begin{proof}
We only prove the lemma when $f\in C^1$; the general case follows from standard approximation arguments. Let $(U, (x_1,\dots,x_d))$ be any smooth local chart of $M$. It is known that $\nabla f=g^{ij} \frac{\pa f}{\pa x_j} \frac{\pa}{\pa x_i}$ in $U$, where $(g^{ij})$ is the inverse of the Riemannian metric $g:=(g_{ij})$. Hence, there exist positive constants $C_1$ and $C_2$, independent of $f$, such that 
\begin{equation}\label{eqn-Jan-25-1}
|\nabla f|^2\leq C_1 \sum_{i=1}^d \left| g^{ij}\frac{\pa f}{\pa x_j}\right|^2\leq C_2 \sum_{i=1}^d \left|\frac{\pa f}{\pa x_j}\right|^2\quad \text{in}\quad U.
\end{equation}

Since $\{X_i\}_{i=1}^n$ spans the tangent bundle $TM$, we find a continuous function $(h_{ij}): U\to \R^{d\times n}$ with full rank such that $\frac{\pa}{\pa x_i}=h_{ij}X_j$ in $U$ for all $i\in \{1,\dots, d\}$, and thus, there is $C_3>0$ such that
$$
\left|\frac{\pa f}{\pa x_i}\right|^2=\left| h_{ij}X_j f\right|^2\leq C_3 \sum_{j=1}^n |X_j f|^2\quad \text{in}\quad U,\quad\forall i\in \{1,\dots, d\}.
$$
This together with \eqref{eqn-Jan-25-1} leads to $|\nabla f|^2 \leq d C_2 C_3\sum_{j=1}^n |X_j f|^2$ in $U$. 

Conversely, as $X_i=X_{ij}\frac{\pa}{\pa x_j}$ for continuous functions $X_{ij}:U\to \R$, $i\in \{1,\dots, n\}$ and $j\in \{1,\dots, d\}$, we see that $\sum_{i=1}^n |X_i f|^2\leq C_4|\nabla f|^2$ in $U$ for some $C_4>0$. As a result, there is $C_5>1$ such that 
$$
\frac{1}{C_5}|\nabla f|\leq \sum_{i=1}^n |X_i f|\leq C_5|\nabla f|\quad \text{in} \quad U.
$$

Since $(U, (x_1,\dots, x_d))$ is any local chart, the desired result follows immediately.
\end{proof}

%-----------------------------------------------------------------%
%\bibliographystyle{amsplain}

\bibliographystyle{plain}
\bibliography{references.bib}

% \begin{thebibliography}{10}

% \bibitem{Zeeman88} E. C. Zeeman,  Stability of dynamical systems. \emph{Nonlinearity} 1 (1988), no. 1, 115-155.

% \bibitem{BKR01} V. I. Bogachev, N. V. Krylov and M. R\"{o}ckner, 
% On regularity of transition probabilities and invariant measures of singular diffusions under minimal conditions. (English summary) \emph{ Comm. Partial Differential Equations} 26 (2001), no. 11-12, 2037-2080.

% \bibitem{Hsu} Elton P. Hsu, {\em A brief introduction to Brownian motion on a Riemannian manifold}. Lecture notes.

% \end{thebibliography}
\end{document}